\documentclass[11pt]{amsart}
\usepackage{amssymb, amsmath, amsthm}
\usepackage[utf8]{inputenc}
\usepackage{graphicx}
\usepackage[final]{hyperref}
\usepackage{color}
\usepackage{epstopdf}
\usepackage{cancel}
\usepackage{yhmath}
 \usepackage{ulem}

\usepackage[a4paper, centering]{geometry}
\usepackage{color}
\usepackage{graphicx}
\graphicspath{{./pics/}}
\usepackage{scalerel,stackengine}

\usepackage{mathtools}

\newtheorem{theorem}{Theorem}
\newtheorem{proposition}[theorem]{Proposition}
\newtheorem{lemma}[theorem]{Lemma}
\newtheorem{corollary}[theorem]{Corollary}

\theoremstyle{definition}
\newtheorem{definition}[theorem]{Definition}
\theoremstyle{remark} 
\newtheorem{remark}[theorem]{Remark}

\definecolor{verde}{RGB}{20,150,100}
\definecolor{purple}{RGB}{200,30,200}

\newcommand{\EEE}{\color{black}}

\stackMath
\newcommand\reallywidecheck[1]{%
\savestack{\tmpbox}{\stretchto{%
  \scaleto{%
    \scalerel*[\widthof{\ensuremath{#1}}]{\kern-.6pt\bigwedge\kern-.6pt}%
    {\rule[-\textheight/2]{1ex}{\textheight}}
  }{\textheight}%
}{0.5ex}}%
\stackon[1pt]{#1}{\scalebox{-1}{\tmpbox}}%
}

\DeclareMathOperator{\dv}{div}

\def\R{\mathbb{R}}
\def\C{\mathbb{C}}

\def\E{\mathcal{E}}

\newcommand{\vps}{\varepsilon}
\newcommand{\Om}{\Omega}

\newcommand{\ra}{\rightarrow}

\newcommand{\lra}{\longrightarrow}
\newcommand{\sr}{\stackrel}
\def\Pr{P_r}

\DeclareMathOperator{\diam}{diam}

\def \e{\varepsilon}

\newcommand{\di}{\operatorname{div}}

\DeclareMathOperator{\Area}{Area}

\begin{document}
\title[]{ 
  The nonlocal isoperimetric problem for polygons:  \\ 
Hardy-Littlewood and Riesz inequalities 
}

\author[]{Beniamin Bogosel, Dorin Bucur, Ilaria Fragal\`a}
\thanks{}

\address [Beniamin Bogosel]{Centre de Math\'ematiques Appliqu\'ees, CNRS,
	\'Ecole polytechnique, Institut Polytechnique de Paris,
	91120 Palaiseau (France)}
\email {beniamin.bogosel@polytechnique.edu}

\address[Dorin Bucur]{
Universit\'e  Savoie Mont Blanc, CNRS  UMR 5127 \\
 Laboratoire de Math\'ematiques, Campus Scientifique \\
73376 Le-Bourget-Du-Lac (France)
}
\email{dorin.bucur@univ-savoie.fr}

\address[Ilaria Fragal\`a]{
Dipartimento di Matematica \\ Politecnico  di Milano \\
Piazza Leonardo da Vinci, 32 \\
20133 Milano (Italy)
}
\email{ilaria.fragala@polimi.it}


\date{\today}

\begin{abstract} Given a non-increasing and radially symmetric kernel in $L ^ 1 _{\rm loc} (\R ^ 2 ; \R _+)$, we investigate counterparts of the classical Hardy-Littlewood and Riesz inequalities  
when the class of admissible domains is the family  of polygons 
 with given  area  and  $N$  sides. The latter corresponds to study the polygonal  isoperimetric problem in nonlocal version. We prove that, for every $N \geq 3$,  the regular $N$-gon is optimal for Hardy-Littlewood inequality. 
Things go differently for  Riesz inequality: while  for $N = 3$ and $N = 4$ it is known that the regular triangle and the square are optimal, for  $N\geq 5$ we prove that symmetry or symmetry breaking may occur (i.e.\ the regular $N$-gon may be optimal or not), 
depending on the value of $N$ and on the choice of  the kernel.   

  \end{abstract}

\maketitle



 

\section{Introduction}

Given a non-negative and non-increasing radially symmetric function $h$ in $L ^ 1 _{\rm loc} (\R ^d)$, called in the sequel an {\it admissible kernel}, 
for any measurable set $\Om \subset \R ^d$ let 
 \begin{equation}\label{f:Jh} 
J _ h (\Om) := \int _\Om \int _\Om h( x-y)   \, dx \, dy\, .
\end{equation} 
 The classical Riesz rearrangement inequality \cite{riesz} states in particular that,
denoting by $\Om ^*$ the ball with the same volume as $\Om$, it holds
$ J _h (\Om) \leq J _h (\Om ^*)$. 
 It is a natural question to ask whether symmetry is preserved when passing 
 to the polygonal setting (in dimension $d=2$): denoting by $\mathcal P _N$  the class of polygons with $N$ sides and area for definiteness equal to $\pi$,     
 and by $\Om ^* _N$  the regular gon in $\mathcal P_N$, this amounts to ask whether 
  \begin{equation}\label{f:pbmax}
  \max \big \{ J _ h (\Om) \ : \ \Om \in \mathcal P _N \big \} = J _ h ( \Om_N ^*) \,.
  \end{equation}


As a general fact, transposing  isoperimetric-type inequalities with balls as optimal domains 
into the setting of polygons with a fixed number of sides is indeed a very natural problem. 
It has been investigated both in the field of geometric measure theory  and in the field of mathematical physics, 
with drastically different levels of difficulties.

For the  isoperimetric inequality in geometric measure theory solved by De Giorgi in \cite{degiorgi}, i.e.\ the minimization of perimeter under volume constraint, the polygonal version
 \begin{equation}\label{f:DG}  {\rm Per} \,  (\Om) \geq  {\rm Per} (\Om ^* _N)\qquad \forall \Om \in \mathcal P _N
\end{equation}
 is  an elementary result, which can be found in several textbooks in convex geometry.

On the other hand, for isoperimetric-type inequalities in mathematical physics, such as Saint-Venant of Faber-Krahn inequalities (for which we refer to the classical monograph \cite{PS} and to the recent book \cite{HdG}), 
the polygonal versions  \begin{equation}\label{f:PS} 
 \lambda _ 1 (\Om) \geq \lambda _ 1 (\Om ^* _N) \qquad \text{ and } \qquad \tau (\Om) \leq \tau (\Om _N ^* )
 \end{equation}
are conjectures formulated several decades ago by P\'olya-Szeg\"o, who also proved them for  $N=3$ and $N=4$. 
Here  $\lambda _1 (\Om)$ is the principal eigenvalue of the Dirichlet Laplacian in $\Om$, while $\tau (\Om)$ is the torsional rigidity of $\Omega$ (namely the $L ^ 1(\Omega)$-norm of the unique solution 
 to the equation $-\Delta u=1$    in $H ^ 1 _0 (\Omega)$). 
The inequalities analogue to \eqref{f:PS}  have been proved for every $N\geq 3$ for the logarithmic capacity \cite{SZ} and for the Cheeger constant \cite{BF16} (for related results, see also \cite{BB22, BF21, Laug21}). 
At present,  \eqref{f:PS}  are open for any $N \geq 5$, and they can be included among the major open problems in shape optimization. Their  validity is also related to a conjecture by Caffarelli and Lin  \cite{CaffLin},  about the  asymptotical optimality of the hexagonal honeycomb for related optimal partition problems, see \cite{BFVV}.

%
To some extent,  Riesz inequality can be viewed as a kind of ``bridge'' 
between the classical isoperimetric inequality \eqref{f:DG} and  the physical inequalities \eqref{f:PS}, as it is intimately connected to each  of them.  It is useful to  briefly explain these connections,  before introducing our results.

The relation between Riesz inequality and the classical  isoperimetric inequality is easily individuated. Indeed,  
the functional $J_h$ in \eqref{f:Jh} differs just by a change of sign and a translation from the 
{\it nonlocal $h$-perimeter}, defined by 
\begin{equation}\label{f:hper} P _ h (\Omega): = \int _\Omega \int _{\Omega^c} h (x-y) \, dx \, dy\,,
\end{equation}
where the quantity $h ( x-y)$ is interpreted as an interaction density between two 
points $x \in \Omega$ and $y \in \Omega ^ c: = \R ^ d \setminus \Omega$. A suitable scaling of $h$ in \eqref{f:hper} allows to recover the usual notion of perimeter via an asymptotic formula.
The  concept of nonlocal perimeter has been first introduced in \cite{BBM01}, and in recent times it has been widely developed, especially
concerning the fractional kernel
$h( x) = {|x| ^ {-(d+s)} }$, $s \in (0, 1)$, and  concerning
bounded integrable kernels (see respectively the seminal papers \cite{CS08, CRS10} and the recent monograph \cite{MRT}).  In particular, the optimality of balls in the nonlocal isoperimetric inequality  has been proved in the fractional case \cite{FS08, FLS08, FFMMM, L14},  for  kernels which are not radially symmetric and decreasing \cite{CN18}, and for some Minkowski type nonlocal perimeters  \cite{CDNV}.  
In this perspective, 
\eqref{f:pbmax} is equivalent to 
  $$\min \big \{ P _ h (\Om) \ : \ \Om \in \mathcal P _N \big \} = P _ h ( \Om_N ^*) \,.$$

To sketch the relation between the polygonal Riesz inequality  and the physical isoperimetric inequalities \eqref{f:PS},
let us focus for simplicity  on the Saint-Venant inequality. 
It is known that  the torsional rigidity satisfies (see \cite[Section 4]{Ca01}) 
 \begin{equation}\label{f:caroll}  \tau (\Om)= \int_0^{+\infty} \Big (\lim_{m \ra +\infty} \int_\Om \dots \int_\Om \prod _{i=1}^m p_{t/m} (x_i-x_{i-1}) dx_0\,\dots dx_m\Big) dt\, , 
 \end{equation} 
 where $p _t$ denotes the heat kernel 
\begin{equation}\label{f:pt}   p_t(x)=  \frac{ e^{- \frac{ | x|^2}{4t}}}{4\pi t}\,.
\end{equation} 
A probabilistic reformulation in the same vein can be given also for the principal frequency, see again \cite[Section 4]{Ca01}. 
Then the classical Saint-Venant  and Faber-Krahn inequalities, with the ball as optimal domain, can be obtained as a consequence of the general rearrangement inequality for multiple integrals due to Brascamp-Lieb-Luttinger \cite[Theorem 1.2]{BLL}, which in $2$-dimensions reads
\begin{equation}\label{f:BLL}\int _\Om \!\! \dots\!\!  \int _\Om \prod _{ i = 1} ^ k h _ i \Big ( \sum_ {j=1} ^m a _{ij} x _j \Big )  \, d x_1 \!  \dots \!  d x_m \leq 
\int _{\Om^*} \!\! \dots \!\! \int _{\Om^*} \prod _{ i = 1} ^ k h _ i ^* \Big ( \sum _{j=1} ^m a _{ij} x _j \Big )  \, d x_1\!   \dots \!  d x_m  \,;
\end{equation} here $h_i$ are measurable non-negative  functions on $\R ^2$ vanishing at infinity,  $h _ i ^*$ are their symmetric decreasing rearrangements,  $\{ a _{ij} \}$ are real numbers, $\Omega$ is a  set of finite Lebesgue measure in $\R^2$, and $\Omega ^*$ is the ball with the same area as  $\Omega$.  

Hence, a possible approach to P\'olya-Sezg\"o conjectures \eqref{f:PS} would be to prove a polygonal  version of  
the result by  Brascamp-Lieb-Luttinger, stating that the inequality  \eqref{f:BLL}  remains true 
 when the integrals  at the left hand side   are extended to a polygon $\Omega  \in \mathcal P _N$,  and the integrals at the right hand side are extended to $\Omega ^* _N$. Of course one has to start  from  small values of $m$. 
 For   $m = 2$  (and a suitable choice of $h _i$ and $\{ a_{ij} \}$), this corresponds exactly to study the Riesz inequality \eqref{q:isop}.  But the problem turns out to be nontrival even for $m = 1$. 
 Indeed in this case it amounts investigate the validity of the following  polygonal version of the classical Hardy-Littlewood inequality (see \cite[Chapter 10]{HLP})  
 \begin{equation}  \label{f:hardy0}  
 \int _\Om h (x) \, dx  \leq  \int _{\Om_N^* }\!   h ( x) \, dx   \qquad \forall \Om \in \mathcal P _N \,  . 
 \end{equation}  
 
 Aim of this paper is to attack  the polygonal Hardy-Littlewood inequality \eqref{f:hardy0} and 
 the polygonal Riesz inequality  
 \begin{equation}\label{q:isop}  J _h (\Om) \leq J _h (\Om ^* _N)\, , \text { or equivalently}  \ P _h (\Om) \geq P _h (\Om ^* _N)   \qquad \forall \Om \in \mathcal P _N\,.
\end{equation}

 To the best of our knowledge, very few results are available in this respect in the literature. 
Concerning the inequality \eqref{f:hardy0},  within the restricted setting of convex polygons, and for particular kernels,  is mentioned as an open question by Fejes-T\'oth in \cite{FT}. Still in the restricted setting of convex polygons, and  for the kernel $h ( x) = |x|$, the inequality is proved  in \cite{morbol}. 
Concerning the inequality \eqref{q:isop}, for $N = 3, 4$ (triangle and quadrilaterals), though not explicitly stated it can be deduced via Steiner symmetrization from Lemma 3.2 in \cite{BLL}. 
 More recently,   in \cite{BCT22}  Bonacini, Cristoferi and Topaloglu have considered the strictly related problem of 
 characterizing ``critical''  triangles and quadrilaterals.  More precisely, a polygon $\Om \in \mathcal P _N$ is said to be critical for $J _h$ if,  
for some positive constant $c$, it holds
\begin{equation}\label{f:area-constrained} 
\frac{d}{d \e}  J _ h (\Om _\e) \Big | _{ \e = 0 } =  c 
\frac{d}{d \e}  | \Om _\e| \Big | _{ \e = 0 }  
\end{equation} 
whenever 
$\{\Om _\e\}$  are  obtained from $\Om$   by one of the following elementary movements:
either the rotation of one side with respect to its midpoint, or the parallel translation of one side with respect to itself
 (see  the Appendix in Section \ref{sec:appendix} for the detailed definitions).  
Clearly, a  polygon maximizing $J _h$ over $\mathcal P _N$ must be a critical polygon.  
In \cite{BCT22} it is proved that, under some weak assumptions on $h$, the regular triangle and the square are respectively the unique critical triangle and quadrilateral. 
It is also conjectured that the same rigidity property holds true for any $N \geq 5$ and that, consequently, the inequality \eqref{q:isop} remains true for every $N \geq 5$.

In this paper we prove that the polygonal Hardy-Littlewood inequality \eqref{f:hardy0} holds true for every admissible kernel $h$ and without any convexity assumption on the admissible polygons (see Theorem \ref{t:hardy}). The situation concerning   
the polygonal Riesz inequality \eqref{q:isop} is more delicate, because a key point turns out to be the choice of the kernel.
In this respect,   the above discussion motivates the assertion that the heat kernel is of special relevance.
We point out that, for $h = p_t$, the corresponding functional $J _h(\Omega)$, namely 
the   quantity
\begin{equation}\label{f:heatc} Q_ {\Om} ( t) := \int_\Om \int_\Om p_t(x-y)dx dy
\end{equation}
is the so-called {\it heat content} of the set $\Om$ at time $t$. It represents the quantity of heat kept by the set $\Om$ once it is warmed at constant temperature $1$ and its heat is  left to diffuse  in the plane. We refer to \cite{Pr04} for related 
isoperimetric properties and to \cite{BeGi16, BeSr90} for related asymptotic expansions in the polygonal setting.
A natural  heuristic way to investigate  the validity of the polygonal inequality \eqref{q:isop}  for $h = p_t$, with $t$ sufficiently small or $t$ sufficiently large,
consists in 
looking at the limiting behaviours of $Q _\Om (t)$ as $t \to 0 ^+ $ and as $t \to + \infty$. In the limit as $t \to  0 ^+$, 
by \cite[Theorem 1]{BeSr90} for any polygon $\Omega$ it holds
$$Q_\Om(t)=|\Om|- \frac{2}{\sqrt \pi} |\partial \Om|t^\frac 12+ O(t)\,,$$
where the $O(t)$ term can be expressed in terms of the inner angles of $\Omega$. Such  asymptotic expansion, combined with the classical isoperimetric inequality \eqref{f:DG}, suggests that inequality \eqref{q:isop} should hold when $h=p _t$ with $t$ sufficiently small. 
On the other hand, 
in the limit as $t \to + \infty$, starting from the asymptotics  of $p _t$ 
and looking at the leading term, we obtain 
$$ 4 \pi t Q _\Om ( t) = |\Om | ^ 2 - \frac{1}{4 t} \int _\Om \int _\Om |x-y| ^ 2 \, dx \, dy  + O \Big (\frac{1}{t^ 2 } \Big ) \,, $$  
so that the inequality  \eqref{q:isop} should hold for $h = p _t$ with $t$ sufficiently large provided it holds for the quadratic kernel $h ( x) = | x| ^ 2$. 

This observation drew our attention to study the inequalities \eqref{q:isop} for the quadratic kernel $h ( x) = | x| ^ 2$, and more generally for power type kernels 
\begin{equation}\label{f:powerh} h ( x) = |x| ^ k\,,  \qquad \text{ with  } k >0 \,.
\end{equation}
Notice carefully that such kernels are not admissible in the sense specified at the beginning of the paper, since they are increasing.   Thus one should write  inequalities \eqref{q:isop} for $h= M - |x| ^ k$, with the constant $M= M (\Omega)$ chosen  so large that $h \geq 0$ on $\Om$.  Equivalently, this amounts to rewrite the reverse inequalities of  \eqref{q:isop}  for $h = |x| ^ k$.  Another family of ``simple'' kernels, which are of natural interest since their linear combinations can serve to approximate any smooth admissible kernel, are 
those of the form
 \begin{equation}\label{f:char}
h(x)=\chi _{B_r(0)}(x)\,
\end{equation}
that we call in the sequel {\it characteristic kernels}.  

Our results about Riesz polygonal inequality are mainly focused on the two families of kernels in \eqref{f:char} and in  \eqref{f:powerh}.
For characteristic kernels 
we prove that, when $r$ is small enough,  inequality \eqref{q:isop} holds  (see Theorem \ref{t:final}); indeed, 
 in this case we also have a rigidity result  characterizing the regular $N$-gon as the unique critical polygon (see Theorem \ref{t:polygons}).
On the other hand,  
when  $r$ is large enough, a big surprise is coming: for  $N$ even, $N \geq 6$, the inequality turns out to be false!  See Theorem \ref{t:false}. 
In particular,   the above mentioned conjecture made in  \cite{BCT22} is, in general, false.
The heuristic reason is that, for $r$ large enough, our problem  is equivalent to the problem 
of finding so-called ``largest small $N$-gons'', namely polygons with fixed diameter and maximal area. 
This is a challenging problem in discrete geometry, for which it is known that symmetry breaking occurs for any $N \geq 6$ even, see Section \ref{sec:results} for more details. 

For power-type kernels, by using our polygonal Hardy-Littlewood inequality,  we prove that the inequality \eqref{q:isop} holds true for $k = 2$ and $k = 4$ (see Theorem \ref{t:powers}). The same strategy fails for  non-integers $k$, for odd integers $k$, as well as for higher even integers $k$. In fact, as a consequence of Theorem \ref{t:false}, the inequality \eqref{q:isop} turns out to be false also for power-type kernels with sufficiently high exponent (still for $N \geq 6$ even).

The conclusion which stems from our analysis is that optimal polygons for  the nonlocal isoperimetric inequality  turn out to be  sensitive to the choice of the kernel $h$, as it may produce either symmetry or symmetry breaking. 
 This  is a highly unexpected phenomenon, which makes the study of the nonlocal isoperimetric inequality quite intriguing: indeed, 
in the light  of our results,  the problem becomes  to understand which  
are specifically the kernels  yielding symmetry for every $N \geq 3$.  
We suspect that this is the case for the heat kernel;  we don't have a proof of this fact, but we give 
some  affirmative numerical results for polynomial approximations of $p_t$ (see Section \ref{sec:numerics}).

\bigskip
{\bf Outline of the paper.} The paper is organized as follows: in Section \ref{sec:results} we state our main results, complemented with some comments 
on their proofs and a list of related open problems; in Section \ref{sec:numerics} we explore numerically some of these open problems;  in the subsequent sections we give  the proofs of the results stated in Section \ref{sec:results}; finally in the Appendix we provide  some first and second order shape derivatives which are used in the proofs.

\section{Main results} \label{sec:results} 

Our main results about the nonlocal polygonal isoperimetric inequality for characteristic kernels 
read as follows.  Below, 
when $h = \chi _{B _ r (0)}$, we set for brevity 
\begin{equation}\label{f:hperr} P _ r (\Omega): = \int _\Omega \int _{\Omega^c} \chi _{B _ r (0)} (x-y) \, dx \, dy
\end{equation}

\begin{theorem} [symmetry]\label{t:final}
For every $N \geq 3$, 
there exists $ r'  =  r' (N) >0$  such that 
\begin{equation}\label{f:newsur}  \forall r \leq r'\,,\quad  \min \Big \{ \Pr (\Om)  \ :\ \Om \in \mathcal P _N \Big \} =   \Pr    ( \Om^*_N  )   \,. \end{equation}   
\end{theorem} 

 {The fact that $r'$ above is in general finite is related to the result of Reinhardt asserting that, when $N \geq 6$ is even, the regular $N$-gon is {\it not} a minimizer of the diameter under an area constraint \cite{R22, S58}. In fact we have the following:}
\begin{theorem}[symmetry breaking] \label{t:false}
For every $N$ even, $N \geq 6$,  setting $ r '' = r'' (N):  = \min \{ {\rm diam}\,  \Omega \, :\,  \Omega \in \mathcal P _N \} $, we have 
  \begin{equation}\label{f:surprise}  
\forall r \in  [r'',   {\rm diam}\,  \Omega ^*_N)   \, , \quad \min \Big \{ \Pr (\Om)  \ :\ \Om \in \mathcal P _N \Big \} <   \Pr    ( \Om ^*_N  )  \,.
\end{equation} 
Moreover,   still for 
$r \in [r'', {\rm diam}\,  \Omega ^*_N)$
the minimum in \eqref{f:surprise} equals $\pi ^ 2 ( r ^ 2 -1)$, and it is attained at a 
polygon $\Om ^\sharp _N$ with diameter $r''$, which  is not the regular polygon.
\end{theorem}

 

Comparing Theorems \ref{t:final} and \ref{t:false} shows  that, for characteristic kernels,
the optimality of the regular $N$-gon does depend on the value of $r$ (at least for $N\geq 6$ even). 
This brings the study of  the nonlocal isoperimetric inequality for polygons into the more complex perspective 
of understanding for which $N$ and for which kernels symmetry or symmetry breaking occurs.

 \smallskip
Before stating some partial results in that direction, let us give some short comments about the proofs of Theorems \ref{t:final} and \ref{t:false}. 

\smallskip 
The proof of Theorem \ref{t:final} is performed in a first stage for convex polygons and then it is extended to the general case. 
For convex polygons, we argue by contradiction.  The idea is that, in the regime of a small $r$, if a convex polygon $\Omega$ minimizes the $r$-perimeter  over $\mathcal P _N$  
and it has $r$-perimeter  strictly 
smaller than $\Omega ^*_N$, then $\Omega$ must be close to $\Omega ^*_N$
(this follows from
a uniform asymptotic estimate for the $r$-perimeter  as $r \to 0$, where the classical perimeter appears in the leading term). 
 In particular, close to $\Omega ^*_N$ there would be a $r$-critical polygon, that is 
a polygon satisfying the stationarity condition  \eqref{f:area-constrained}, for $h = \chi _{B _ r (0)}$. When $r$ is sufficiently small, 
this is not possible thanks to a symmetry result for critical polygons  that we state separately in Theorem \ref{t:polygons} below, 
since  it may have an independent interest.  The second part of the proof dealing with  arbitrary polygons requires some more refined arguments, in particular since 
minimizing sequences may converge to  a ``generalized polygon'' (precisely in the sense of Definition \ref{d:genpol}), possibly containing self-intersections in its boundary.   Roughly speaking, 
the idea is to reduce the problem to a situation similar to the convex setting: this is achieved by exploiting triangulations in order to identify local concentrations of mass, and  by localizing our estimates near the sides of the limit polygon 
where there is no accumulation of vertices.  We refer to Section \ref{sec:final} for the detailed proof.

\smallskip  Let us now state the afore mentioned symmetry result for $r$-critical polygons. 
To that aim, it is convenient to reformulate more explicitly the shape derivative in the left hand side of  \eqref{f:area-constrained}. 
This has been done in \cite{BCT22}, but to make the paper self-contained we enclose a proof in the Appendix, see Lemma \ref{l:derivatives}.  
The outcome is the following: 
if $\{ \Om _ \e \}$ are obtained from $\Om$  respectively
by rotating the side $[A_i A _{i+1}]$ with respect to its midpoint $M _i$, or
by a parallel movement of such side  with respect to itself,   the stationarity condition 
\eqref{f:area-constrained} amounts to ask that, setting  $v_\Om (x) := \int _\Om h ( x-y  ) dy$, it holds 
\begin{eqnarray}
\int _{A_i}^ {M_i} \!\!\!\! v_\Om ( x) \,  |x M_i| \, dx - \int _{M_i}^ {A_{i+1} } \!\!\!\!   v_\Om ( x)  \,  |x M_i| \, dx  = 0   &  \label{f:angles}
\\ \noalign{\smallskip}
\int _{A_i}^ {A_{i+1}}v_\Om ( x) \, dx  = c \, \mathcal H ^ 1 (  [A_i A _{i+1}] ) \,. \qquad      &
& \label{f:sides}
\end{eqnarray}  
 The result below states that  
the validity of eqs \eqref{f:angles}-\eqref{f:sides} for $i = 1, \dots, N$ enforces symmetry,  provided the support of the kernel is  small enough. It can be viewed as a polygonal version of the Alexandrov-type symmetry recently proved in \cite[Corollary 7]{BF}. 
The proof  is obtained by using 
a reflection argument which is reminiscent of  \cite{BCT22, FV}.

\begin{theorem}[symmetry for $r$-critical polygons]\label{t:polygons} Let  $h$ be an admissible kernel  and   
let $\Om \in \mathcal P _N$ satisfy equations \eqref{f:angles}-\eqref{f:sides} for every $i = 1, \dots, N$. 
 When $N >3$,  assume further that $ {\rm spt} (h) \subseteq  {B _ r(0)}$, with $r$ such that 
\begin{equation}\label{f:card}   \partial \Om\cap  B _ r (x) \text{ is contained into a pair of consecutive sides of $\Om$}   \qquad \forall x\in \partial \Om \,.
\end{equation}  

Then $\Om$ is a regular $N$-gon.
\end{theorem}

\begin{figure}[ht]
		\centering
		\begin{tabular}{cc}
		\includegraphics[height=0.35\textwidth]{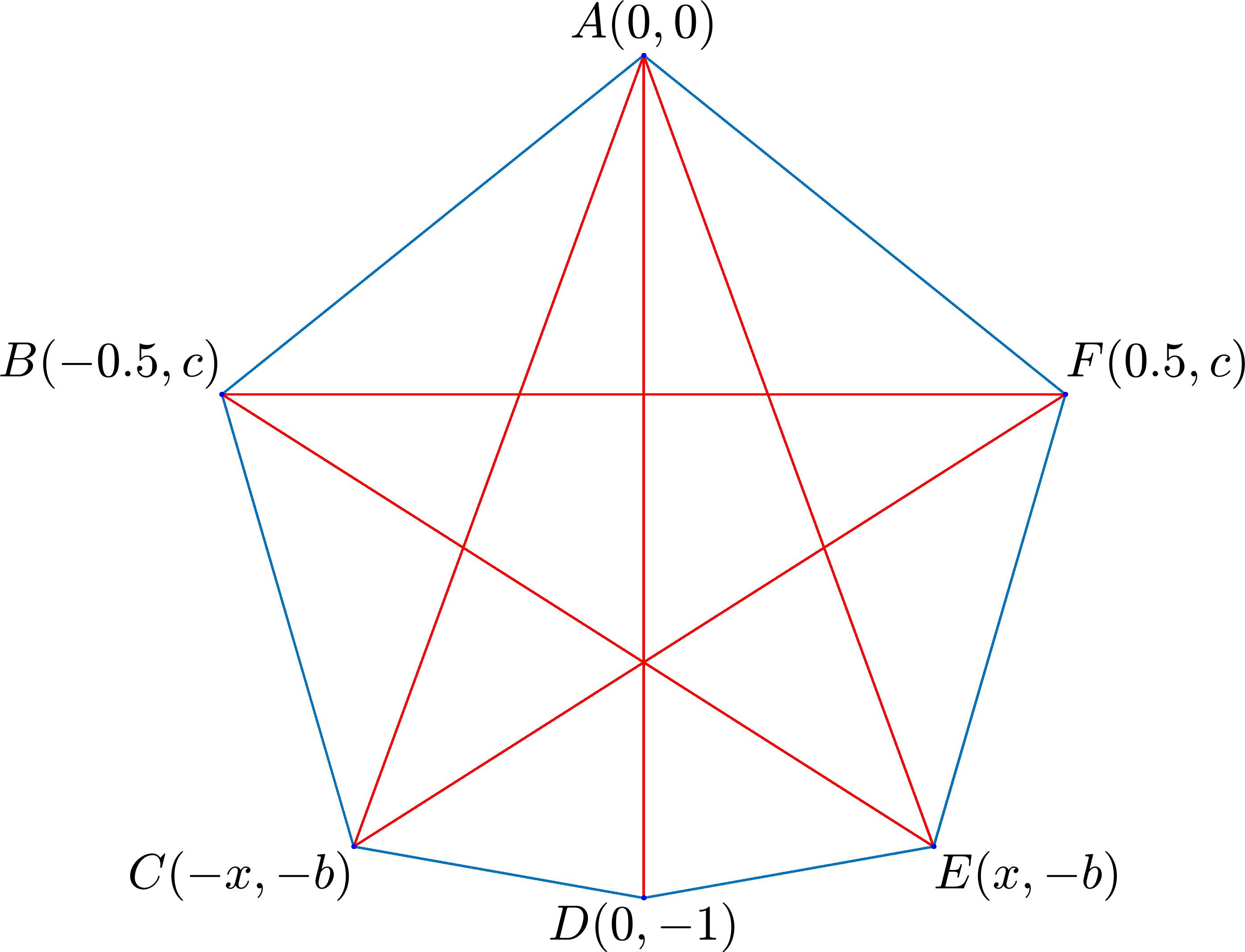}&
		\includegraphics[height=0.35\textwidth]{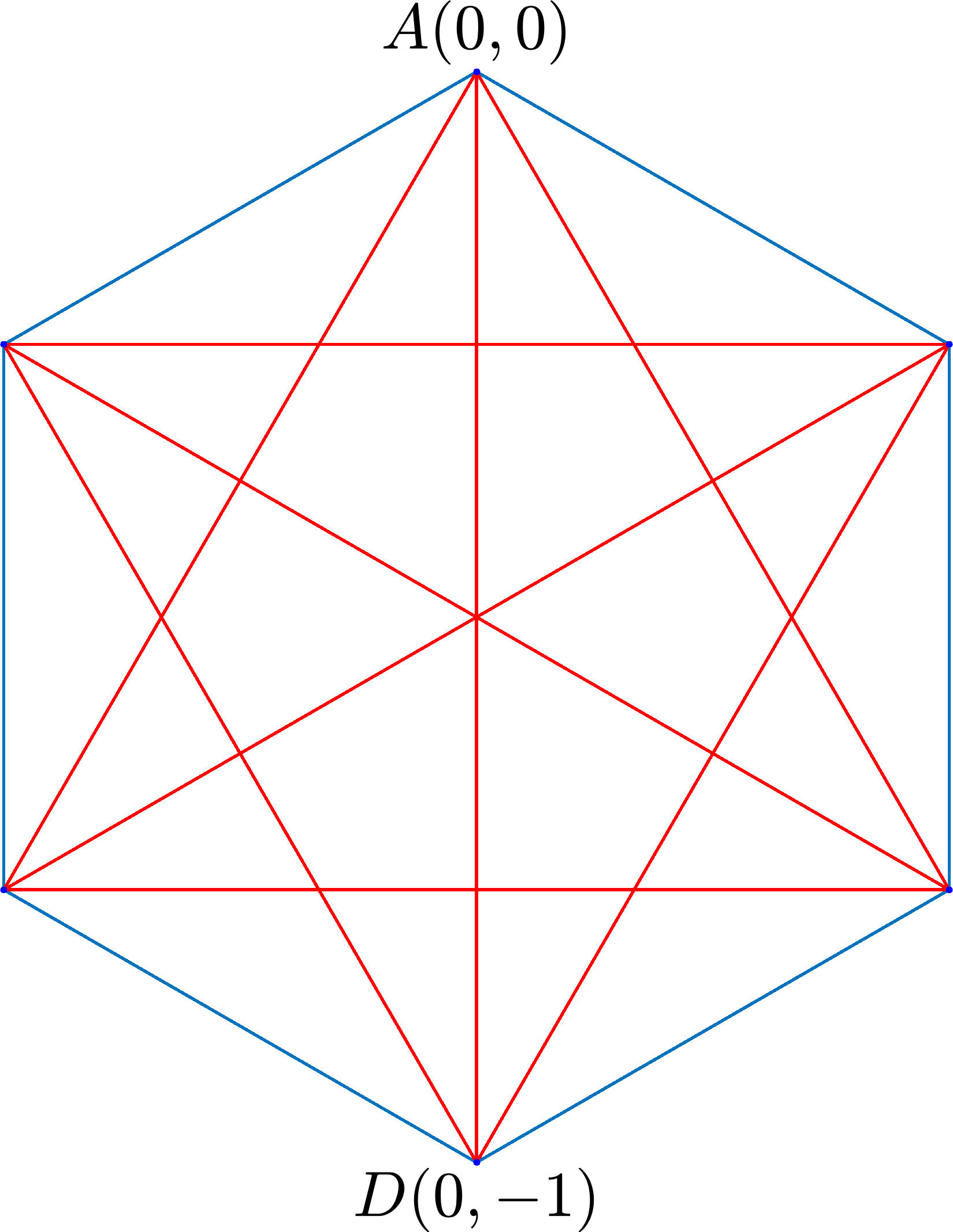} 
		\\
		Area = $0.674981$ & Area = $0.649519$
	\end{tabular}

		 \caption{The Graham hexagon (left) and the regular hexagon (right) with unit diameter.}   
		\label{fig:Graham}
	\end{figure}

\smallskip
The proof of Theorem \ref{t:false} is  obtained in a completely different way; indeed, it follows as a rather straightforward consequence of  a result by Reinhardt asserting that, when $N \geq 6$ is even, the regular $N$-gon is {\it not} a minimizer of the diameter under an area constraint
(while,  for $N$ odd, the regular $N$-gon is a minimizer), see \cite{R22, S58} and the expository paper \cite{Moss}.
 In the particular case $N = 6$, the optimal hexagon  $\Omega ^\sharp _6$  was found by  Graham \cite{graham, B61}, see Figure \ref{fig:Graham}.
 Its construction 
     can be done as follows.  First fix two points  $A(0,0)$ and $D(0,-1)$  at distance one, and then, denoting $c=d-b$ we determine the other four vertices  $B(-0.5,c)$, $C(-x,b)$,  $E(x,b)$, and  $F(0.5,c)$: 
taking $x$ as a parameter,  $b$ and $d$ are found from the relations
	$x^2+b^2=1$ and $(x+0.5)^2+d^2 = 1$.
	For $b = 0.939053346$ and $d = 0.536702650$, a numerical value for $x$ is given by
	$x = 0.343771453$. \footnote{The construction is taken from the MathWorld page
	\href{{https://mathworld.wolfram.com/GrahamsBiggestLittleHexagon.html}}{\nolinkurl{https://mathworld.wolfram.com/GrahamsBiggestLittleHexagon.html}}} 
 
 	Denoting by $H_G$ and $H_R$ suitable scalings of Graham and  regular hexagons, we have: 
 	\begin{itemize}
 		\item At {\it fixed diameter}, the area of the Graham hexagon is greater than the area of the regular one, the ratio of the areas being $\frac{\Area(H_G)}{\Area(H_R)} = 1.039201$.
 		\item At {\it fixed area}, the diameter of the Graham hexagon is smaller than the diameter of the regular one and the ratio of their diameters is $\frac{\diam(H_G)}{\diam(H_R)} = 0.980957$.
 	\end{itemize}

For  $N > 6 $ even,  the determination of a polygon
 minimizing the diameter under an area constraint is a challenging problem in discrete geometry, which remains, to the best of our knowledge, open. Let us just mention that, among equilateral polygons, it has been recently proved that the regular polygon is optimal \cite{A17}. Among arbitrary polygons, it has been proved in  \cite{FS07} that for every $N$ even the optimal polygon enjoys the following  property conjectured by Graham: its skeleton (namely the collection of diameters connecting any two vertices) is not an asterisk but consists of a $(N-1)$-cycle and one additional edge.  
Moreover, for $N = 8$, very accurate numerical solutions have been proposed in \cite{AHMX02} and \cite{Hmess}. We also refer to the recent paper \cite{AHS21}
for more recent symbolic calculations and for further bibliography.



\bigskip
We now come back to the  problem of minimizing the $h$-perimeter over $\mathcal P _N$.
We are going to 
focus on  some specific kernels which are not characteristic functions.  
In this respect, let us point out that the phenomenon of symmetry breaking  is not a prerogative of characteristic kernels.
 Indeed  from Theorem \ref{t:false} and keeping the same notation as in its statement, we easily obtain the following  
  
\begin{corollary}\label{c:corpotenze}  Let $N$ even, $N \geq 6$. 
\begin{itemize}
\item[(i)] If $h$ is a smooth admissible kernel close enough to $\chi _ {B _ r}(0)$, with $r \in [r'',   {\rm diam}\,  \Omega ^*_N) $, symmetry breaking occurs, i.e. 
 $P _h (\Om _N ^ \sharp) < P _h (\Om _N ^*)$. 
 
\smallskip 
\item[(ii)] If $h(x) = |x| ^k$, there exists $\overline k$ such that symmetry breaking occurs for every $k \geq \overline k$, i.e.
$P _h (\Om _N ^ \sharp) > P _h (\Om _N ^*)$. 
In particular, for $N=6$ we have $\overline k \leq 2832$. \end{itemize}
\end{corollary}  

 \EEE 
Regarding statement (ii) above recall that, since power-type kernels $h ( x) = |x| ^ k$ with $k>0$ are increasing, the corresponding inequalities \eqref{q:isop}
must be reversed.  Dealing with such  kernels, in view of Corollary \ref{c:corpotenze},  
the question becomes whether  the regular polygon  is a maximizer of the nonlocal perimeter at least for small $k$.
We show that the answer is affirmative in the particular cases  $k =2$ and $k= 4$:

 \begin{theorem}[symmetry for power-type kernels] \label{t:powers}  
 Let $h( x) = |x| ^ 2$ or $h  (x) = |x| ^ 4$. For every $N \geq 3$, 
 the regular poygon $\Om^*_N$ maximizes $P _h$ over $\mathcal P _N$. 
 In equivalent terms, we have  the following polygonal Riesz inequalities
   \begin{eqnarray}
  & \displaystyle \int _\Om \int _\Om  |x-y|^2  \, dx \, dy \geq  \int _{\Om_N^* }\!  \int _{\Om_N^*}  |x-y| ^2 \, dx \, dy \,, \qquad \forall \Om \in \mathcal P _N\,, \  N \geq 3 \, ;  & \label{f:power2} 
  \\ 
  & \displaystyle \int _\Om \int _\Om  |x-y|^4  \, dx \, dy \geq  \int _{\Om_N^* }\!  \int _{\Om_N^*}  |x-y| ^4 \, dx \, dy \,, \qquad \forall \Om \in \mathcal P _N\,, \ N \geq 3 \, . & \label{f:power4}  
    \end{eqnarray}  
 \end{theorem} 

The proof of  Theorem \ref{t:powers} relies on the idea to reduce the study of 
inequalities \ref{f:power2}-\eqref{f:power4} 
to the study of Hardy-Littlewood type polygonal inequalities.  More precisely, 
by writing explicitly the polynomials $|x-y| ^ 2$ and $|x-y| ^ 4$,  it turns out that the minimization of their
double integral   over $\Omega \times \Omega$ 
is equivalent to the minimization of the single integrals $\int_\Omega |x| ^ 2$ and  $\int_\Omega |x| ^ 4$. 
We are thus led to the following question: is it true that $\Omega ^*_N$ minimizes over $\mathcal P _N$ 
an integral  functional of the type $\int_\Omega |x| ^ k$?  More in general, for any admissible kernel $h$, 
one is led to investigate the following Hardy-Littlewood type maximization problem
\begin{equation}\label{p:hardy}
\max \Big \{ \int _\Om  h (x) \, dx \ :\ \Om \in \mathcal P _N \} \,.
\end{equation} 
As mentioned in the Introduction, in the restricted setting of convex polygons, 
problem \eqref{p:hardy} is mentioned as an open question by Fejes T\'oth in \cite{FT} when $h$ is a characteristic kernel. In this case, it amounts to solve the following purely
geometric problem:  find the polygon in $\mathcal P _N$ which maximizes the overlap with the ball $B_r (0)$.   Despite its elementary formulation,  the solution to
such geometric problem is far from being immediate, and it is also the heart of the matter in order to  solve  problem \eqref{p:hardy} for arbitrary kernels.  As in the proof of Theorem \ref{t:final}, the difficulty comes mainly from the fact that maximizing sequences of polygons may converge to  ``generalized polygons''  with self-intersections in their boundary. We overcome this difficulty via an  ad hoc 
geometric construction, allowing to reduce ourselves to deal with star-shaped polygons; once made this restriction,  we can take advantage 
of first order optimality conditions, which enable us to arrive at the regular $N$-gon.

\begin{theorem}[polygonal Hardy-Littlewood inequality]\label{t:hardy}  
Let  $h$ be an admissible kernel. 
For every $N \geq 3$, we have
\begin{equation}\label{f:hardy} 
 \int _\Om  h (x) \, dx  \leq \int _{\Om_N^*}  h ( x) \, dx \,, \qquad \forall \Om \in \mathcal P _N \,.
\end{equation}
\end{theorem}

Applying inequality \eqref{f:hardy} allows us to prove \eqref{f:power2}-\eqref{f:power4}, but the same strategy is not successful to obtain the analogous inequality
\begin{equation}\label{f:riesz-pot}
\int _\Om \int _\Om  |x-y|^k  \, dx \, dy \geq  \int _{\Om_N^* }\!  \int _{\Om_N^*}  |x-y| ^k \, dx \, dy \,, \qquad \forall \Om \in \mathcal P _N\,, \  N \geq 3
\end{equation}
for non-integers, or odd integers, or higher exponents $k$. 
We are just able to prove that \eqref{f:riesz-pot} continues to hold in some very  specific situations, that we gather in the statement below:

\begin{lemma}\label{l:specific} Inequality \eqref{f:riesz-pot} holds in the following cases: 
\begin{itemize}
\item[(i)] 
$k = 6$, $N = 8$,  under the restriction that $\Omega$ is convex and axisymmetric;  

\smallskip 
\item[(ii)] $k\ge 1$, $N\geq 3$,  under the restriction that $\Omega$ is  a linear image of $\Omega _N ^*$. 
\end{itemize} 
 \end{lemma} 

\medskip
\begin{remark} Theorem \ref{t:hardy} allows to extend the result in \cite{morbol} by Morgan and Bolton about the optimality of the hexagonal economic regions for the location problem to other kernels than the average distance, for instance power-type kernels.
\end{remark}

%
%
%
%
%
%
%
%
%

\medskip 
Clearly  our results raise many new questions,  some of which  may be very challenging.
 A short list is  given below. 
 
\bigskip 

{\bf Open problems} 

\begin{itemize}
\smallskip
\item [(A)] Characteristic kernels: Determine or estimate the radius $r'$ in Theorem \ref{t:final}. 

\medskip
\item [(B)] Power-type kernels: Determine for which 
values of $k$ the inequality  
\eqref{f:riesz-pot}
 holds. 
 \medskip
\item [(C)] Gaussian kernel $h(x)= e^{-|x|^2}$: Determine whether $\Omega ^*_N$ minimizes  the $h$-perimeter over $\mathcal P _N$   for any $N \ge 3$. 
(Alternatively, in terms of the heat content defined in \eqref{f:heatc}, does the inequality $Q_\Om(t) \le Q_{\Om_N^*} (t)$ $\forall \Om  \in {\mathcal P}_N$ hold 
for any $N \ge 3$ and every $t >0$?) 

\smallskip 
\item [(D)]  Arbitrary admissible kernels $h$: Determine whether $\Omega ^*_N$ minimizes  the $h$-perimeter over $\mathcal P _N$   for every $N$ odd ($N \geq 5$)
and whether, under some suitable assumptions on $h$,  the same holds for every $N$ even ($N \geq 6$). 

\smallskip
\item[(E)] More general kernels: Explore what happens also for kernels which are not locally integrable, 
but induce a  finite perimeter on the class of polygons,  as it is for instance the case for the  fractional kernel (see e.g.\ \cite[Corollary 1.2]{lomb}).

\end{itemize}

\EEE

\section{Numerical results about problems (B) and (C)}  \label{sec:numerics} 

In this section we bring some numerical evidence related to open problems (B) and (C).
The numerical results  are summarized below. Some of them (mainly on the local minimality) could be  turned into analytical ones  provided the approximations would be controlled and the numerical computations  certified. Let us point out that, as soon as the kernel is of polynomial type  {with sufficiently small degree}, the computations we perform are  {accurate up to rounding errors in double precision}. The computations use quadrature rules  {which are exact for low degree polynomials}. This is explained in Section \ref{sec:computstrat}.

\subsection {About problem (B)} \label{sec:subB}
We made multiple numerical optimizations  with randomized initialization, 
in order to minimize  over $\mathcal P _N$ the functional 
\begin{equation}\label{f:Jk} \int _\Omega \int _\Omega |x-y| ^ k \, dx \, dy. 
\end{equation}
We used the constrained optimization algorithm \texttt{interior-point} from the Matlab \texttt{fmincon} routine.  
 The computations were performed for $N \in \{5,6,7,8,9,10\}$, and  
 for $k \in \{6,8,10,12\}$. 
All simulations led  to the regular polygon.

\medskip 
Next,  in order to extract information about the local minimality of the regular $N$-gon, we looked at the 
sign of the eigenvalues of the Hessian matrix of the scale-invariant functional defined for all polygons with $N$ sides by 
 $$|\Omega|^{-(k+4)/2} \int _\Omega \int _\Omega |x-y| ^ k \, dx \, dy. $$ 

Needless to say,  since the above functional is invariant under rigid motions, several  zero eigenvalues must be expected, so that local minimality is gained as soon as the other eigenvalues are strictly positive.

We computed the Hessian matrix of our functional under vertices displacement, by using formula \eqref{eq:Hessian-F} in the Appendix and the classical Hessian formula for the area functional which can be found e.g.~in \cite[Section 2]{BB22}. We obtained $4$ eigenvalues equal to zero (corresponding to translations, rotations, and homotheties) and $2N-4$ eigenvalues which are strictly positive. 
	
The computations were performed for $N \in \{5,6,7,8,9,10\}$, and  
 for $k$ even, $k \leq 24$.

\medskip
\subsection {About problem (C)} \label{sec:subC}
 Differently from power-type kernels,  the  Gaussian kernel is no longer homogeneous under homotheties. 
Hence,  after rescaling, it not restrictive to consider the heat kernel  $h(x)= e^{- {|x|^2}/{t}}$ at different times, and work  
with polygons with fixed diameter equal to $1$ instead of polygons with fixed area $\pi$. 
Actually,  dealing with polygons with fixed diameter  turns out to be convenient in order to control the approximation  made when the power series expansion of the heat kernel  
\[ h(x,y) = \sum_{k=0}^\infty \frac{1}{k!} \left( -\frac{|x-y|^2}{t}\right)^k\]
is replaced by its partial sums\begin{equation}
h_Q(x,y) = \sum_{k=0}^Q \frac{1}{k!} \left( -\frac{|x-y|^2}{t}\right)^k\, .
\label{eq:Q-kernel}
\end{equation}
Notice that, in the  case $Q=0$, problem (C) becomes trivial since $h _ 0 = 1/t$, and also in the case $Q=1$ the regular $N$-gon maximizes $J_1(\Omega)$ for every $N$ thanks to Theorem  \ref{t:powers}. 

The idea is then to look at what happens for higher values of $Q$  such that the approximation of $h$ by $h _ Q$ is sufficiently good. 
If $x, y$ belong to a polygon with unit diameter and taking $t \geq 1$,   by the inequality 
$|e^z - \sum_{k=0}^Q \frac{z^k}{k!}|\leq \frac{|z|^{Q+1}}{(Q+1)!}$ holding for any  $z \in \C$ with $\text{Re} z<0$, we have
\[ |h(x,y)-h_Q(x,y)| \leq \frac{1}{(Q+1)!}.\]
In particular, for $Q=12$ (which in our computational strategy described in Section \ref{sec:computstrat} below  corresponds to a quadrature rule of order $24$) the inverse of $(Q+1) !$  is bounded from above by $1.6\times 10^{-10}$.  Consequently, for polygons with unit diameter  (having area at most $\pi$), the global numerical error done in  evaluating $J _{h _ Q}$ in place of $J _ h$ is  bounded from above by  $\pi  ^2 \times 1.6\times 10^{-10}$. Similar estimates yields global errors smaller than $10^{-7}$  when replacing    the gradient and the Hessian of  $J _ h$ by their analogues for $J _{h _ Q}$, according to the integral formulas in Appendix.

Then, we fix our attention on the functional $J _ { h _ Q}$ for $Q = 12$.  Clearly, working with such functional brings us back to a polynomial setting as in case of problem (B) discussed above,  with the difference   that now $J _ { h _ Q}$  is no longer homogeneous with respect to scalings. 
Hence we consider the Hessian matrix associated with the functional 
\begin{equation} \int _\Omega \int _\Omega h _ Q ( x, y) \, dx \, dy   -\ell_Q |\Omega|,
\label{eq:lagrange}
\end{equation}
being $\ell_Q$ a Lagrange multiplier chosen so that the regular $N$-gon under consideration is a critical point.  {The Lagrange multiplier $\ell_Q$ is the ratio of the norms of the gradient of  $\int_{\Omega}\int_{\Om} h_Q(x,y)dxdy$ and the gradient of the area; notice indeed that for the regular $N$-gon these gradients are collinear for symmetry reasons.}
The Hessian matrix is obtained again by using formula \eqref{eq:Hessian-F}.  We investigate the sign of its eigenvalues corresponding to 
 eigenvectors orthogonal to the gradient of the area (a space of dimension $2N-1$). 
 We obtained $3$ eigenvalues equal to zero (corresponding to translations and  rotations) and $2N-4$ eigenvalues which are strictly negative. 
	
 The computations where performed for $N \in \{ 5, 6, 7, 8, 9, 10\}$ and for a few choices of $t\in [1,100]$, including the endpoints. Numerically, we observe that the  Hessian eigenvalues vary monotonically with $t$: they are negative and have a decreasing absolute value as $t$ increases. Therefore, we conjecture that their sign remains  {negative} for all values of $t$ in the considered range.


Furthermore, the smallest absolute value of non-zero eigenvalues, which is obtained for $N=10$ and $t=100$, is larger than $10^{-4}$. The above discussion about the error estimates done when replacing the heat kernel by its polynomial approximation $h _ Q$ indicates that, 
also for the heat kernel, the regular $N$-gon is  a local maximizer  under area constraint.

These simulations motivate us to conjecture that  the regular $N$-gon is a local maximizer of the heat content for every $t\in \Bbb{R}$.

Let us also mention that the same computations above were made also for $Q \in [2,11]$: for $Q \in [2,5]$,  $N \in [5,10]$ and various choices of $t\in [1,100]$ an oscillatory behavior can be observed, namely the Hessian at the regular $N$-gon may have positive or negative eigenvalues;  however, for $Q\geq 6$ the behavior stabilizes and the non-zero eigenvalues become strictly  {negative}.

\medskip
\subsection {Computational strategy}\label{sec:computstrat}  
Let us now briefly  explain the strategy adopted for the computations in Sections \ref{sec:subB} and \ref{sec:subC}.
When $\Omega$ is a polygon, and $k$ is a positive even integer, functionals of the type  $\int_{\Omega} \int _\Omega |x-y|^k \, dx \, dy$
can be computed explicitly in terms of the coordinates of the vertices,  and the same assertion holds for the integrals involved in the  shape derivatives of such functionals.
However, the resulting expressions are difficult to interpret and implement.
Thus  we choose to adopt a different approach, based on quadrature rules. 
Given a $N$-gon $\Omega$, we split it into triangles $T_1,...,T_N$ (using an inner node) and we decompose the energy as
\[ \int_\Omega \int_ \Omega h(x,y)dxdy = \sum_{i,j=1}^N \int_{T_i}\!\int _{  T_j} h(x,y)dxdy\,,\]
so that we can focus on the computation of integrals  made over a product of triangles. 
A quadrature rule  for an integral over a triangle 
 $T$   is an approximation of the form
\begin{equation}
\int_T f(x)dx \approx \sum_{i=1}^M w_i f(P_i)\, , 
\label{eq:quad-formula}
\end{equation}
where $P_1,...,P_M$ are points in $T$ (expressed, for instance,  using barycentric coordinates in the triangle $T$), and
$w_1,...,w_M$ are the associated weights. 
A quadrature rule is said to be of order $k$ if the approximation \eqref{eq:quad-formula} is exact when $f$ is a polynomial of total degree at most equal to $k$. For any degree $k$, there exist quadrature rules of such degree, the number of quadrature points being increasing with respect to $k$. 
 
\begin{figure}[h]
	\centering
	\includegraphics[width=0.35\textwidth]{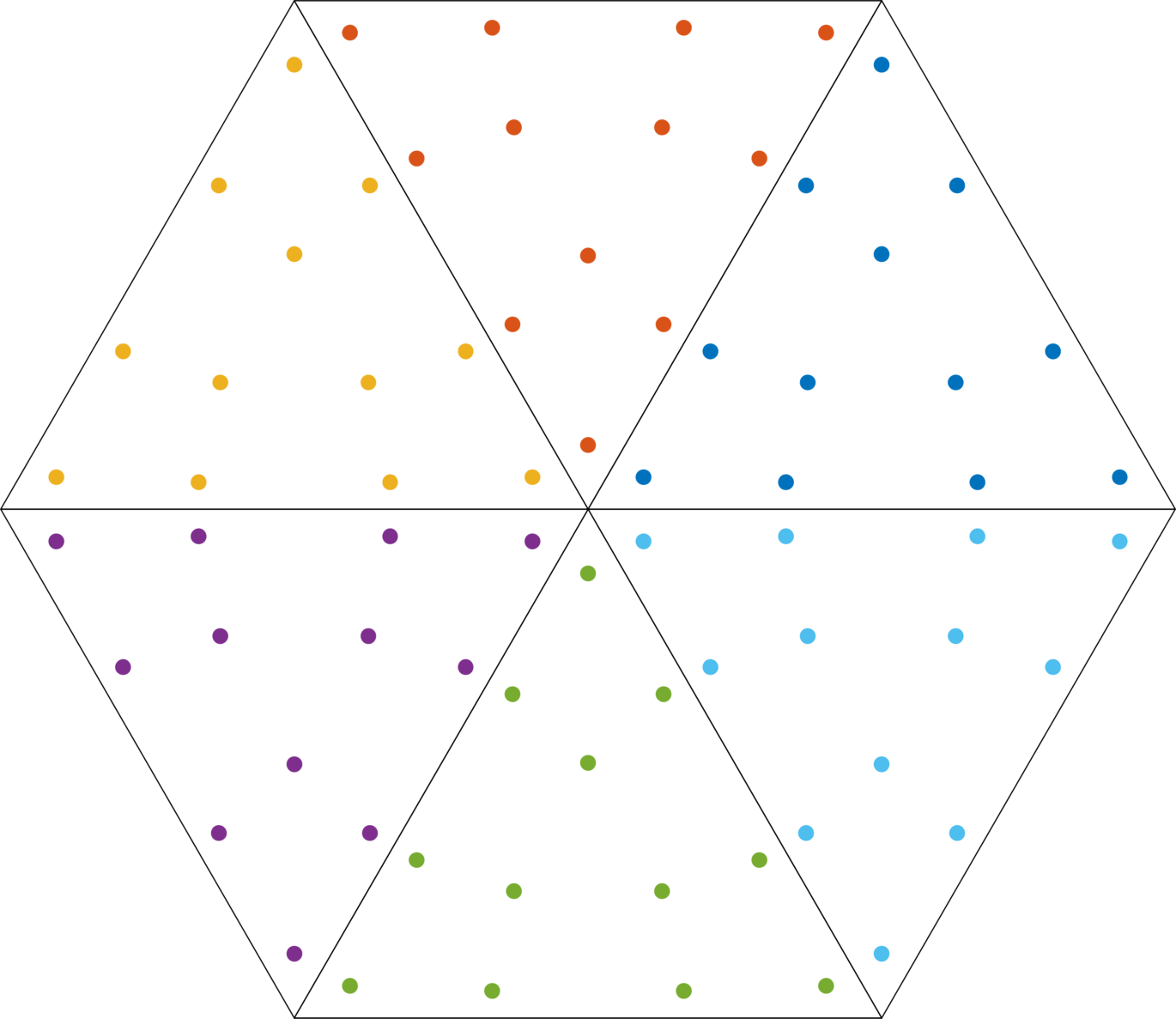}\quad 
	\includegraphics[width=0.35\textwidth]{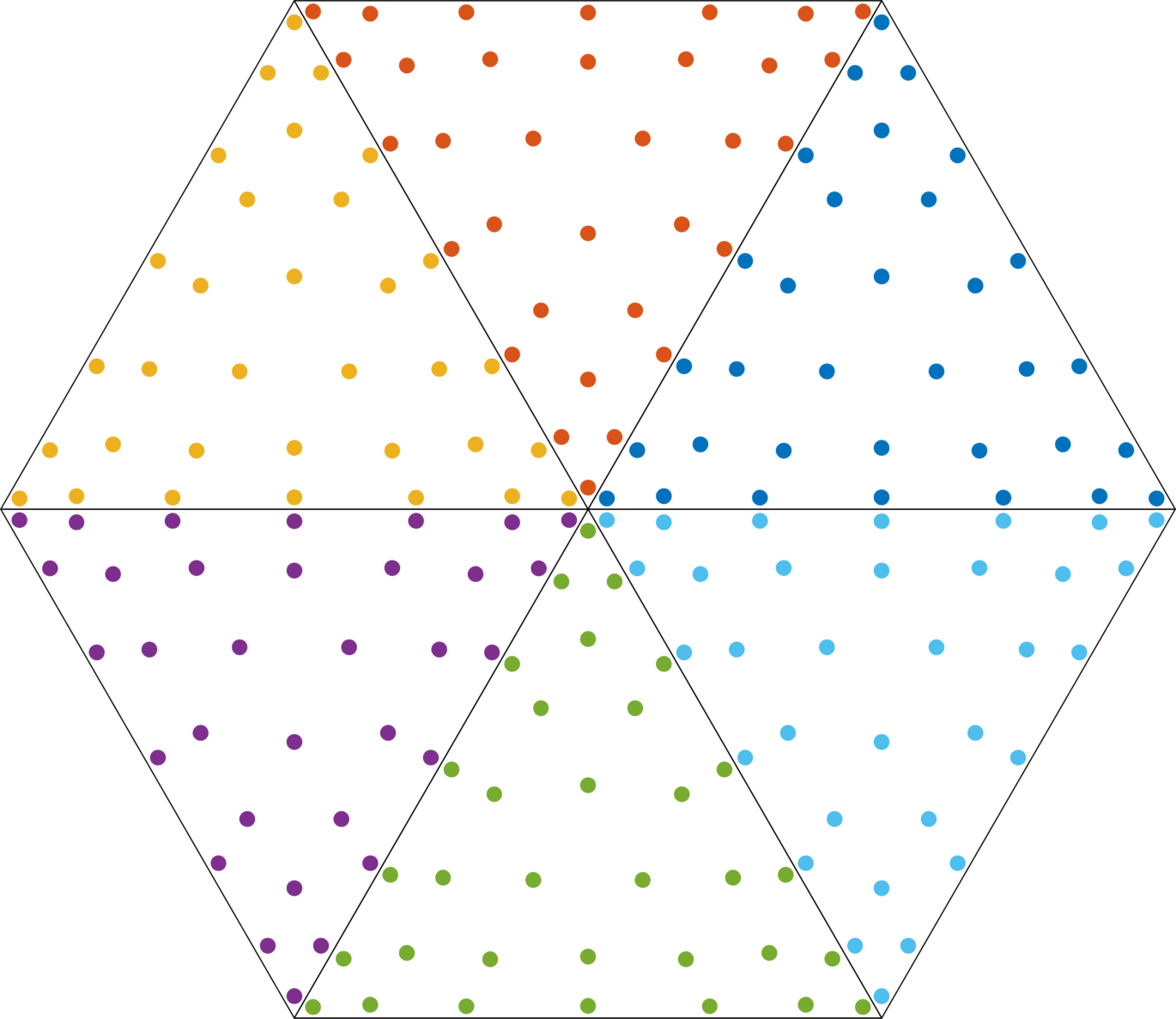}\quad
	\caption{Triangulation and corresponding quadrature points for a regular hexagon. The quadrature rules have degrees $6$ and $12$, respectively.}
	\label{fig:quad-poly} 
\end{figure}

An example of a triangulation and choice of quadrature points of degrees $6$ and $12$ for a regular hexagon  is shown in Figure \ref{fig:quad-poly}.

Handling quadratures rules involving double integrals is more complex, but relies on the same principles. In this case, given two triangles $T_1,T_2$ with corresponding quadrature points $P_1,...,P_M$, $Q_1,...,Q_M$ (having the same barycentric coordinates) and weights $w_1,...,w_M$,  we have
\begin{equation}
\int_{T_1\times T_2} h(x,y)dxdy \approx \sum_{i,j=1}^M w_iw_j h(P_i,Q_j);
\label{eq:quad-formula-prod}
\end{equation}
as above,  the quadrature rule is of order $k$  if  the approximation 
\eqref{eq:quad-formula-prod} is exact when 
$h(x,y)$ is a polynomial of total degree at most $k$. 

 In order to generate quadrature rules required in our computations,  we used the Matlab toolbox \href{https://mathworks.com/matlabcentral/fileexchange/72131-quadtriangle}{Quadtriangle} (accessed in November 2022). We used \emph{non product rules}, included in the referenced toolbox up to degree $25$.

 \section{Proof of Theorem \ref{t:final}}\label{sec:final} 
We proceed to prove the result first in the simplified setting of convex polygons and then in the general case.

 \subsection{Proof  of Theorem \ref{t:final} in the convex setting.}\label{sec:convex} 
 
\begin{figure}[ht]
	\includegraphics[width=0.4\textwidth]{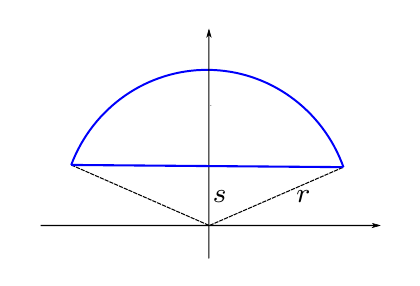}
	\caption{The set $\Delta _{r, s}$.}
	\label{fig:delta}   
\end{figure}

To prepare the proof, it is useful to introduce the set   $\Delta _{r, s}$ defined by
 $$\Delta _{r, s} := \Gamma^+_{0, s} \cap B _ r (0) \,,$$
 where
$$ \Gamma ^ + _{0, s}:= \Big \{ x=  (x_1, x_2)  \in \R ^ 2 \ :\ x_2 \geq s  \Big \} $$
 Clearly $\Delta_{r, s}$ is empty set for $s>r$, while 
for $s \in \ [0, r]$ it is given a circular segment of radius $r$ and  apothem $s$, see Figure \ref{fig:delta}. 

For $s \in [0, r]$,  the Lebesgue measure of $\Delta _{r, s}$ is immediately determined as
\begin{equation}\label{f:delta} |\Delta_ {r, s}| =    r ^ 2 \arccos \left ( \frac{s}{r} \right ) - s \sqrt {r ^ 2 - s ^ 2  } \,. \end{equation} 
Below we state  a simple geometric lemma which plays a key role in the proof; it provides a lower bound for the $r$-perimeter of a convex $N$-gon and an upper bound for the $r$-perimeter of a regular $N$-gon. They will be exploited in the limit of a vanishing radius.  
We focus on the case $N\ge 5$, since for $N=3,4$  the equality in \eqref{f:newsur}   is valid for every $r>0$.  
\begin{lemma}\label{l:ventagli} (i) Let $\Omega$ be a convex polygon in $ \mathcal P _N$. Assume that, for every side  $S_i$ of $\Omega$, denoting by $\ell_i$  its length and by $\theta _i , \theta _{i+1}$ its adjacent inner angles, it holds
\begin{equation}
\label{f:lunga1}
\ell _i  - {r}{\cot \Big ( \frac{\theta _i }{2} \Big ) } - {r}{\cot \Big ( \frac{\theta _{i +1 } }{2} \Big ) }  > 0 \qquad \forall i = 1, \dots, N \,.
\end{equation} 
Then
\begin{equation}\label{f:tutta}P _ r (\Omega)  \geq    \frac{2}{3}|\partial \Omega|  r ^ 3   - \frac{4}{3} \sum _{i = 1} ^ N {\cot \Big ( \frac{\theta _i }{2} \Big ) }  r ^ 4 =: \Phi _ r (\Omega) \,.\end{equation}
Moreover,  if \eqref{f:lunga1} fails for a family of sides (of cardinality at most $N-3$), the inequality \eqref{f:tutta} holds with $\Phi _ r (\Omega)$ replaced by $\Phi _ r (\widehat \Omega)$, where $\widehat\Omega$ is the convex polygon obtained  by eliminating any such side and prolonguing its two consecutive sides. 

(ii) Let $\Omega ^* _N$ denote a regular $N$-gon, with  $N \geq 5$.  Assume that,  denoting by $\ell$  the length   of its sides, it holds
\begin{equation}\label{f:lunga2}
\ell - 2 r >0\,.
\end{equation} 
 Then 
\begin{equation}\label{f:stimareg} P _ r (\Omega ^* _N) \leq \frac{2}{3} |\partial \Omega ^*_N|  r ^ 3 +  4N \Big (\pi - \frac{1}{3} \Big )   \,    r ^ 4\,.  
\end{equation} 

\end{lemma} 

\proof
 (i) Setting $\Om _ s:= \{ x \in \Om \, : \, {\rm dist} (x, \partial \Om) > s \}$, we have 
\begin{equation}\label{f:tubular}P _ r (\Omega) =  \int _\Om |\Om ^ c \cap B _ r ( x) |\, dx  = \int _ 0 ^ r  \, ds \int _{\partial \Om _ s}  |\Om ^ c \cap B _ r ( y) |\, d\mathcal H ^ 1 (y)\,.  \end{equation}
Thanks to assumption \eqref{f:lunga1}, for every $i = 1, \dots, N$ and every $s \in (0, r)$, $\partial \Om _s$ contains a segment $\Gamma_i$  of positive length, 
 made by points $y$ such that $\Om ^ c \cap B _ r ( y) $ intersects $\partial \Omega$ only along the side $S _i$ and  is congruent to $\Delta_{r,s}$.    
  For the  length of this segment we have the following lower bound: 
$$ |\Gamma_i |  \geq \ell _i  - {r}{\cot \Big ( \frac{\theta _i }{2} \Big ) } - {r}{\cot \Big ( \frac{\theta _{i +1 } }{2} \Big ) }\,,$$
and,   for points $y \in \Gamma_i$,    it holds
\begin{equation}\label{f:dueterzi} \int _ 0 ^ r |\Om ^ c \cap B _ r ( y) | \, ds =  \int_0 ^ r  |\Delta_{r,s}| \, ds = \frac{2}{3} r ^ 3 \,,
\end{equation}
where the last equality follows from \eqref{f:delta} and an elementary integration. 
Hence,  
$$P _ r (\Omega)  \geq  \sum _{i = 1  }^ N \Big ( \ell _i  - {r}{\cot \Big ( \frac{\theta _i }{2} \Big ) } - {r}{\cot \Big ( \frac{\theta _{i +1 } }{2} \Big ) }  \Big ) \frac{2}{3} \pi r ^ 3  
 =  \frac{2}{3}|\partial \Omega|  r ^ 3   - \frac{4}{3} \sum _{i = 1} ^ N {\cot \Big ( \frac{\theta _i }{2} \Big ) }  r ^ 4\,.$$ 


In case assumption \eqref{f:lunga1} fails for some index, we repeat the proof above with the following only modification:  in correspondence of any 
index for which \eqref{f:lunga1} is false, we remove that side from $\Omega$ and we consider the polygon $\widehat \Omega$ defined as in the statement. 
By construction, for every side $\widehat S_i$ of $\widehat \Omega$ and every $s \in ( 0, r)$, $\partial \Omega _s$ contains a segment $\widehat \Gamma _i$  
 of positive length (parallel to $\widehat S_i$),
 made by points $y$ such that $\Om ^ c \cap B _ r ( y) $  is congruent to $\Delta_{r, s}$.    
  For the  length of this segment we have now the following lower bound: 
$$ |\widehat \Gamma_i |  \geq  \widehat \ell _i  - {r}{\cot \Big ( \frac{\widehat \theta _i }{2} \Big ) } - {r}{\cot \Big ( \frac{ \widehat \theta _{i +1 } }{2} \Big ) }\,,$$
 where  $\widehat \ell _i$, $\widehat \theta _i$, and $\widehat \theta _{i+1}$ denote the length of $\widehat S_i$, and it adjacent angles. 
Summing over all the sides of $\widehat \Omega$, we find   the lower bound $P _ r (\Om) \geq \Phi _ r (\widehat \Om)$.

\medskip  
(ii)  We write the equality \eqref{f:tubular} for $\Omega ^ *_N$.  For every  $s \in (0, r)$,  the set of points $y$ at distance $s$ from $\partial \Omega ^ *_N$   
contains $N$ segments 
of positive length, bounded from above by the positive quantity $\ell - 2r$, 
made by points $y$ such that $\Om ^ c \cap B _ r ( y) $ is congruent to $\Delta_{s, r}$.  For points $y$ in such  segments, the equality  \eqref{f:dueterzi} holds.  
For points  $x\in \Omega$ such that $B _ r ( x) $ meets more than one side of $\partial \Omega ^*_N$, we simply estimate from above 
$|(\Om ^*_N ) ^ c \cap B _ r ( x) |$ by $\pi r ^ 2$.  The measure of these points is bounded from above by $4 N r ^ 2$. We end up with
$$P _ r (\Om ^*_N ) \leq  N ( \ell - 2 r) \frac{2}{3}  r ^ 3   + 4N  r ^ 2 ( \pi r ^ 2) \,.   $$ \qed

\bigskip
We are  now ready to prove Theorem \ref{t:final} for convex polygons.  We argue by contradiction. Assume the statement is false. 
Then, there exists an infinitesimal sequence of  radii $\{ r _k \}$ and a sequence of convex polygons $\{\Omega _k  \}\subset \mathcal P _N $ such that
\begin{equation}\label{f:absurd} P _ { r _k} (\Om_k ) < P _ { r _k } (\Om _N ^*)\,.
\end{equation}
Here and in the remaining of the proof, $\Om _N ^*$ denotes a regular $N$-gon of area $\pi$. 
By possibly passing to a subsequence and up to translations, the sequence of convex polygons $\{\Omega _k\}$ admits a limit  $\Omega _0$ in the Hausdorff complementary topology. 
There are two possibilities: either $\Omega _0 = \emptyset$, or  $\Omega _0 \neq \emptyset$. Let us show that both cases lead to a contradiction. 

\medskip
{\it Case 1)}: $\Omega _0 = \emptyset$. Let us consider the sequence of (possibly empty) convex polygons contained into  $\Omega _k$ defined by 
$$\omega _k := \Big \{ x \in \Omega _k \ :\ {\rm dist } ( x, \Omega_k ^ c ) \geq \frac{r_k}{2} \Big \}\,.$$ 
Up to a subsequence, we may distinguish to subcases: either $ \sup _k |\partial \omega _k|  < + \infty$, or 
 $|\partial \omega _k|  \to + \infty$. 
 
 \medskip 
 {\it Case 1a)}: $ \sup _k |\partial \omega _k|  < + \infty$. Up to a further subsequence,  either the convex polygons  $\omega _k$ are  empty, or they converge to a segment.  Anyhow, we have  that $|\Omega _k \setminus \omega _k| \to \pi$.  
 For every $x \in (\Omega _k \setminus \omega _k  )$, the intersection 
$  \Omega ^ c \cap B _ {r _k} ( x)$ contains a set congruent to $\Delta_{ r _k, \frac{r_k}{2}} $. Hence  
 \begin{equation}\label{f:uper} P _ {r_k} (\Om _k)  \geq 
  | \Omega _k \setminus \omega _k |  | \Delta _{  r _k, \frac{r_k}{2}} | = | \Omega _k \setminus \omega _k |  \Big ( \frac{\pi}{3} - \frac{\sqrt 3}{4} \Big ) r  _k ^ 2  \,,
  \end{equation}
  where in the last equality we have used \eqref{f:delta}.
  On the other hand, for $k$ sufficiently large assumption \eqref{f:lunga2} is fulfilled and hence the inequality \eqref{f:stimareg} in Lemma \ref{l:ventagli} holds with $r = r _k$. 
  In view of this fact, and since $|\Omega _k \setminus \omega _k| \to \pi$, the inequality \eqref{f:uper} contradicts \eqref{f:absurd} in the limit as $r _k \to 0$. 

\medskip 
 {\it Case 1b)}: $ |\partial \omega _k|  \to + \infty $.  We have
\begin{equation}\label{f:tubular2}P _ {r_k} (\Omega _k) \geq   |\partial \omega _k| \int _0 ^ { \frac{r_k}{2}}  |\Delta _{ r _k, s} |  \, ds  = 
\frac{ |\partial \omega _k| }{24}  (16+ 4 \pi  - 9  \sqrt 3  ) r_k^2 \,.
 \end{equation} 
Since \eqref{f:stimareg} holds with $r = r _k$, and since  $|\partial \omega _k|  \to + \infty $, the inequality \eqref{f:tubular2} contradicts \eqref{f:absurd} in the limit as $r _k \to 0$. 

\medskip
{\it Case 2)}: $\Omega _0 \neq \emptyset$. We apply Lemma \ref{l:ventagli} (i) to the sequence of convex polygons $\Omega _k$. 
We observe that, since $r _ k \to 0$,  for $k$ large enough assumption \eqref{f:lunga1} is certainly satisfied, except possibly  for certain indices 
corresponding to sides of infinitesimal length. 
Thus we have
$$P _ { r_k} (\Om _k) \geq \Phi _ {r_k} (\Om _k) \qquad (\text{or alternatively } P _ { r_k} (\Om _k) \geq \Phi _ {r_k} (\widehat \Om _k) )\,.$$ 
We observe that the coefficients of the polynomial function $\Om \mapsto \Phi _ r (\Om)$ only depend on the perimeter and on the inner angles of the polygon $\Om$ (see \eqref{f:tutta}); moreover, the  same holds for the polynomial function $\Om \mapsto \Phi _ r ( \widehat \Om)$, because the perimeter and the inner angles of $\widehat \Om$ can be easily expressed in terms of 
the perimeter and the inner angles of $\Om$. 
Now, since $\Omega_k$ converge to $\Omega _0$, the perimeter  and the inner angles of $\Omega _k$  converge respectively to the  perimeter and to the inner angles of $\Omega _0$. We conclude that, for $k$ sufficiently large,  the following lower bound holds: 
$$P _ {r_k}  (\Omega _k) \geq  \frac{2}{3}  |\partial \Omega| r_k  ^ 3   - C_0  r _k ^ 4  \,,$$ 
 where $C_0$ is a fixed constant independent of $k$. 
By combining the above lower bound with Lemma \ref{l:ventagli} (ii) (which applies since its assumption \eqref{f:lunga2} is satisfied for $k$ sufficiently large), we obtain that 
$$ \frac{2}{3}  |\partial \Omega| r_k  ^ 3   - C_0  r _k ^ 4  \leq  \frac{2}{3} |\partial \Omega ^*_N|  r  _k^ 3 +  4N \Big (\pi - \frac{1}{3} \Big )   \,    r _k^ 4 $$
and hence
$$\limsup _{k\to + \infty } |\partial \Omega _k  | \leq |\partial \Omega ^*_N|\,.$$   
By the  classical isoperimetric inequality for convex polygons, this implies that $\Omega _0 = \Omega ^* _N$. 
To conclude, we observe that it is not restrictive to assume that $\Omega _k$  is a minimizer of the $r_k$-perimeter over the class of convex polygons in $\mathcal P _N$ with area $\pi$.  (Notice that such a minimizer exists for any $r_k$ sufficiently small, because otherwise a maximizing sequence of polygons would degenerate, yielding a contradiction by the same arguments used to deal with  Case 1) above.) 
Then we have found a sequence of critical polygons  for the $r_k$-perimeter, converging to  $\Omega ^* _N$, and satisfying \eqref{f:absurd}. Since $r _k \to 0$, this contradicts Theorem \ref{t:polygons}. \qed

 \subsection{Proof  of Theorem \ref{t:final} in the general case.}\label{sec:general0} Also in this case we prepare the proof with a geometric lemma. 
 For every $\theta \in ( - \frac{\pi}{2}, \frac{\pi}{2})$ and every $s>0$, we set 
  $$\begin{array}{ll} 
  & \Gamma _{\theta, s}:= \Big \{ x=  (x_1, x_2)  \in \R ^ 2 \ :\  x_2  =  (\tan \theta) x_1 + s    \Big \} 
   \\  \noalign{\medskip} 
  & \Gamma ^ + _{\theta, s}:= \Big \{ x=  (x_1, x_2)  \in \R ^ 2 \ :\ x_2 \geq (\tan \theta) x_1 + s  \Big \} 
 \\  \noalign{\medskip} 
 &\Gamma^ - _{\theta, s}:= \Big \{x=  (x_1, x_2) \in \R ^ 2 \ :\ x_2 \leq (\tan \theta) x _1  + s  \Big \}  
 \end{array}
 $$

Given a family of straight lines $\{ \Gamma _{ \theta _i, s _i} \}$, for $i = 1, \dots, 2q-1$ which do not intersect each other in $B _ r (0)$, 
with   $\theta_i \in \big (\!-\!\frac{\pi}{2}, \frac{\pi}{2} \big )$ and  
  $   0 \le   s_1 <  \dots < s_{2q-1}$, consider the union of strips
\begin{equation}\label{f:sigma} 
\Sigma   : =   \Gamma   ^-_{\theta  _{1}  , s _{1}}  \cup \bigcup _{i = 1} ^ {q -1}  \big ( \Gamma   ^+_{\theta  _{2i}  , s _{2i} }   \cap   \Gamma  ^-_{\theta  _{2i +1}  
 s  _{2i+1}  } \big )  \,,
 \end{equation} 
see Figure \ref{fig:sigma}. 
Recalling that that $\Delta _{r, s}$ is the set defined at the beginning of Section \ref{sec:convex},   we prove the following estimate for the measure of $\Sigma ^c \cap B _ r (0)$:

\begin{figure} [h] 
\centering   
\includegraphics[width=0.55\textwidth]{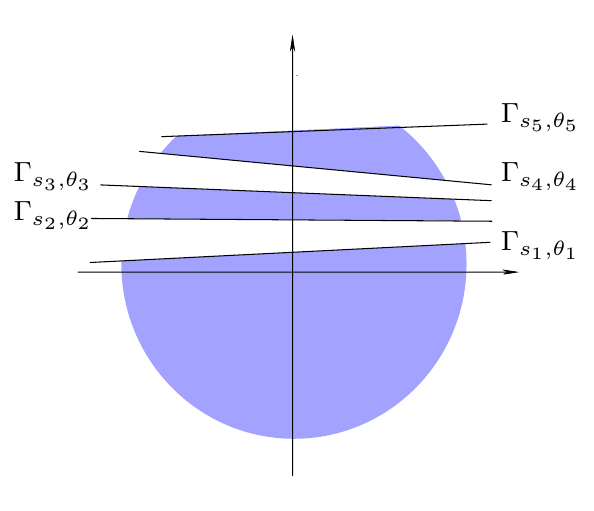}
\caption{The intersection of the set $\Sigma$ in \eqref{f:sigma} with $B _ r (0)$}
\label{fig:sigma}   
\end{figure}

\begin{lemma}\label{l:strips} Let $\Sigma$ be defined by \eqref{f:sigma}. Setting $\overline s:= \mathcal H ^ 1  \big (\Sigma \cap  (\{0 \} \times [0, +\infty))  \big )$, it holds 
\begin{equation}\label{f:millefoglie}
 |  \Sigma ^ c \cap B _r (0)|  \geq | \Delta _{r,  \overline s} | -  4 r ^ 2 \displaystyle  \sum _{  i = 1 }^{2 q -1}  |\theta_i|   \,. 
\end{equation}
\end{lemma}
\proof    In order to prove  \eqref{f:millefoglie}, it is not restrictive to assume that $\overline s \in [0, r]$, since otherwise it is trivially satisfied. We notice first that, 
for every $\theta \in \big (\!-\!\frac{\pi}{2}, \frac{\pi}{2} \big )$ and for every $s \in  [0, r] $, it holds
\begin{equation}\label{f:ventrot} 
|\Delta_{r, s} | -  4 |\theta| r ^ 2 \leq   | \Gamma^+_{\theta, s}\cap  B _ r (0)  |  \leq   |\Delta_{r, s} |  +  4 |\theta| r ^ 2\,.
\end{equation}
Namely,  for $\theta = 0$, \eqref{f:ventrot}  holds with equality signs, by the definition of $\Delta _{r, s}$. 
For $\theta \neq 0$, it holds since the  symmetric difference between  $\Gamma^+_{\theta, s}\cap  B _ r (0)$ and $\Delta _{r, s}$ is contained into the union of the two sets
$$B _ r (0) \cap (\Gamma ^+ _{\theta, s} \cap \Gamma ^-  _{0, s}  ) \quad \text{ and }  \quad B _ r (0) \cap (\Gamma ^+ _{0, s} \cap \Gamma ^-_ {\theta, s}  )  \,,$$
and each of these sets has measure bounded from above by $ \frac{ |\theta|} {2} (2r) ^ 2$. 
Now, the inequality \eqref{f:millefoglie} is a  consequence  of \eqref{f:ventrot} and of the monotonicity of   the map 
$$ [0, r]\ni s \mapsto {\mathcal H}^1 (\Gamma_{0,s} \cap B_r(0)),$$  
which implies that, for $0\le s\le  a<b$, it holds $|\Delta _{r,a}\setminus \Delta _{r,b}| \le |\Delta_{r,s} \setminus \Delta_{r, s+b-a}|$. 
\qed  

\begin{remark}\label{bufr23.06} 
In   Lemma \ref{l:strips},  inequality \eqref{f:millefoglie} remains trivially valid replacing $\Sigma$ by any  subset $\widetilde \Sigma \subseteq \Sigma$. We shall use this argument for subsets $\widetilde \Sigma $ of the type
\begin{equation}\label{f:sigma2} 
\widetilde \Sigma   : =   (\Gamma   ^+_{\theta  _{0}  , s _{0}}  \cap \Gamma   ^-_{\theta  _{1}  , s _{1}} ) \cup \bigcup _{i = 1} ^ {q -1}  \big ( \Gamma   ^+_{\theta  _{2i}  , s _{2i} }   \cap   \Gamma  ^-_{\theta  _{2i +1},  
 s  _{2i+1}  } \big )  \,,
 \end{equation} 
where $s_0 \le 0$ and  $\Gamma  _{\theta  _{0}  , s _{0}} $ does not intersect any  $\Gamma   _{\theta  _{i}  , s _{i}}$ with $i = 1, \dots, 2q-1$,  in $B_r(0)$.
\end{remark}

\bigskip
As a further preliminary, let us give the following definition, that was already used by the second and third authors in \cite{BF16}.

\begin{definition}\label{d:genpol} A {\it generalized polygon with at most $N$-sides} is the limit in the $H ^ c _{\rm loc}$ topology of a sequence $\{ \Om _n \}$ of classical polygons with at most $N$ sides (meant as {\it open} polygons) such that $\limsup _n | \Om _n | < + \infty$. 
\end{definition}

Recall that the convergence of $\{\Om _n\}$ to  $\Om$  in the $H ^ c _{\rm loc}$ topology means that, for every ball $B$, we have 
$\lim _n d_{H ^c} (\Om _n\cap B, \Om \cap B ) = 0$,   
$d_{H ^c}$ being the {\it Hausdorff complementary distance}, namely  
$$d_{H ^c} (\Om _n\cap B, \Om \cap B  ) := \sup _{x \in \R ^ 2} \big | {\rm dist} (x , (\Om _n\cap B) ^c) -
{\rm dist} (x ,  (\Om \cap B ) ^c)\big | \,,$$
(${\rm dist}$ stands for the Euclidean distance from  a closed set). 

As a consequence of well-known properties of such topology (see for instance \cite{BuBu, HP}), any generalized polygon  
is an open set of finite Lebesgue measure, which is simply connected (as its complement is connected), but possibly disconnected. 
Any connected component is  delimited by a finite number of line segments, which are pairwise joined at their endpoints to form a closed path, possibly containing self-intersections, given by points or line segments. 


\bigskip
We  now ready to prove Theorem \ref{t:final} for arbitrary polygons.  As in the convex case, we argue by contradiction, and we denote by 
$\Om _N ^*$ a regular $N$-gon or area $\pi$. 
If the statement is false  there exists an infinitesimal sequence of positive radii $\{ r _k \}$ and a sequence of polygons $\{\Omega _k  \}\subset \mathcal P _N $ such that
\begin{equation}\label{f:absurd2} P _ { r _k} (\Om_k ) < P _ { r _k } (\Om _N ^*)\,.
\end{equation}
To achieve the proof it is enough to show that 
\begin{equation}\label{f:convreg}\Omega _k \, \stackrel {H ^ c}    \to  \, \Omega ^*_N\,.
\end{equation}
Indeed,  the convergence \eqref{f:convreg} implies in particular that $\Omega _k$ is convex for $k$ large enough, which is a contradiction since we have already proved the statement for convex polygons. 

In order to prove \eqref{f:convreg}, we consider for every $k$ a triangulation of $\Omega _k$ made by $N-2$ disjoint (open) triangles $\{T ^ 1_k, \dots, T ^ {N-2}_k\}$, with vertices and sides belonging to the family of vertices and diagonals of $\Omega _k$, such that 
$$\overline \Omega _k = \overline T ^ 1_k \cup \dots \cup \overline  T _k ^ {N-2}\,.$$
Up to subsequences (which here and in the sequel are not relabeled),  there exist  sequences of vectors $\{ y _k ^j\}$ and sets $T^ j$ (which may be either a triangle or the empty set),  such that 
\begin{equation}\label{f:convtriang} T _k ^j - y ^j_k \, \stackrel {H ^ c}   \to \,  T ^ j \qquad \forall j = 1, \dots, N-2 \,.
\end{equation}  
For every $j = 1, \dots, N-2$, one of the following three situations occurs:
\begin{itemize}
\item[(a)]  $T ^ j = \emptyset$ and $|T ^ j _k | \to m ^ j >0$
\smallskip

\item[(b)] $T ^ j = \emptyset$ and $|T ^ j _k | \to 0$ 
\smallskip

\item[(c)] $T ^ j \neq \emptyset$ 
 \end{itemize} 
For convenience, we divide the remaining of the proof in three steps. 

\bigskip
{\it Step1: Situation (a) cannot occur.}   

Assume by contradiction that we are in situation (a) for some sequence $T _k ^ j$.  Hereafter, we omit for simplicity the index $j$.  
Then, up to a subsequence and to a rigid motion,  the vertices of $T _k$ are given by  
$$(-\ell _k, 0 ) \, , \quad (\ell _k, 0) \, , \quad (a_k, b _k), \quad \text{ with } b _ k >0 \,,
$$ 
where the horizontal side of length $2 \ell _k$ is the longest one and, 
as $k \to \infty$, 
$$\ell _k \to + \infty \quad \hbox { and } \quad \ell _k b _ k \to m >0 \,.$$ 
We divide the segment $[- \ell _k, \ell _k] $ into $2N$ equal segments, of length $\frac{\ell _k }{N}$. At least one of them, say up to a translation 
the segment $[- \frac{ \ell _k}{2N}, \frac{\ell _k}{2N}] $, is such that the half strip
$$S _k:= \big  [- \frac{ \ell _k}{2N}, \frac{\ell _k}{2N} \big ]  \times [0, +\infty)$$ 
does not contain any other vertex of $\Omega _k$. 
On the other hand, $S_k$ is crossed  by a certain number of sides of $\Omega _k$, including two sides of $T_k$.  
Thus, with the notation introduced at the beginning of Section \ref{sec:general0}, there exist angles $\theta ^ k _{ i}  \in \big (-\frac{\pi}{2}, \frac{\pi}{2} \big ) $
and positive numbers $s ^ k _i$, for $i = 0, \dots,   2 M _k   -1$  such that $0 = s ^ k _0 < s ^ k _1 < s ^k _2 < \dots s ^ k _{2 M _k  - 1 } $, and
$$ \begin{array}{ll} 
&  \displaystyle \Omega _ k \cap S _k = \bigcup _{i = 0} ^ {M_k  -1 }  \big ( \Gamma  ^+_{\theta _ k ^{2i}  , s _ k^{2i} }   \cap   \Gamma  ^-_{\theta _ k ^{2i +1}  , s _ k^{2i+1}  } \big )
\\ \noalign{\bigskip}
& \displaystyle T _ k \cap S _k =  \Gamma  ^+_{\theta _k ^{0}  , s _ k^{0} }   \cap   \Gamma  ^-_{\theta _ k ^{1}  , s _ k^{1}  } \,. 
\end{array}
$$ 
 Up to subsequences, for every $i =  0, \dots, M _k-1\EEE$, we have 
 $$\theta _k ^ {2i},  \theta _k ^ {2i+1} \to \theta ^i \in \big [ - \frac{\pi}{2}, \frac{\pi}{2} \big ]\, , \qquad s _k  ^ {2i},  s _k  ^ {2i+1},\to s ^i \in [ 0 , + \infty ]\,.$$ 
 
Moreover, by construction it holds $$0 = s ^ 0 \leq s ^ 1 \leq \dots \leq s ^{M_k-1} \,.$$ 
 Hence we may define   $q $  as  the largest index in $\{ 0, \dots , M _{k} -1  \}$ such that $s^q = 0$.  Notice that, as a consequence, we also have 
 \begin{equation}\label{f:convangles}
 \theta ^ 0 = \theta ^ 1 = \dots = \theta ^q = 0.
 \end{equation}

We set 
$$\Sigma _ k  : =  \Gamma  ^-_{\theta _ k ^ {1}  , s _ k^{1} } \cup  \bigcup _{i = 1} ^ {q }  \big ( \Gamma  ^+_{\theta _  k ^{2i}  , s _ k^{2i} }   \cap   \Gamma  ^-_{\theta _ k ^{2i +1}  
 s _ k^{2i+1}  } \big )  \,.$$ 
 Let us observe that, if $q <M_k-1$, then
 \begin{equation}\label{bufr23.01}
 \liminf_{k \ra +\infty} {\rm dist} \Big ( \Sigma_k\cap \big ( \big [ - \frac{\ell _k}{4N} , \frac{\ell _k}{4N} \big ] \times [ 0 , + \infty )  \big ), (S_k\cap \Om_k) \setminus \Sigma_k \Big ) >0 .
 \end{equation}
 Then, the idea is to adopt similar arguments as in  the convex case, just by using Lemma \ref{l:strips} in place of Lemma \ref{l:ventagli}.   More precisely, 
 we consider the set 
$$\omega  _k := \Big \{ x = ( x_1, x_2) \in \Sigma _k  \cap   \big ( \big [ - \frac{\ell _k}{4N} , \frac{\ell _k}{4N} \big ] \times [ 0 , + \infty )  \big )    \ :\   s_x \geq \frac{r_k}{2} \Big \}\,,$$ 
where
  $$s_x := \mathcal H ^ 1 \Big (  \Sigma _ k \cap \big ( \{ x_1 \} \times [ 0, + \infty) \big )  \Big ) \,.$$   
 In view of Lemma  \ref{l:strips}  and of the strict inequality \eqref{bufr23.01}, for $k$ large enough and for every $x \in \Sigma _k  \cap   \big ( \big [ - \frac{\ell _k}{4N} , \frac{\ell _k}{4N} \big ] \times [ 0 , + \infty )  \big )    $, it holds 
  $$|\Om_k^c \cap B_{r_k} (x)| \ge | \Delta _{r_k, s_x} | -  4 r_k ^ 2 \displaystyle  \sum _{  i = 0 }^{2 q +1}  |\theta_k^i|. $$ 
Then we follow the same proof as in Case 1) of Section \ref{sec:convex} to get a contradiction.
More precisely, 
denoting by $\gamma_k$ the projection of $\omega_k$ on the horizontal axis,  we distinguish the two cases $ \sup _k {\mathcal H}^1(\gamma_k)  < + \infty$, and 
 $  {\mathcal H}^1(\gamma_k)  \to + \infty$.   
 
Assume that $ \sup _k {\mathcal H}^1(\gamma_k)  < + \infty$. 
 This implies that $|\omega_k| \rightarrow 0$. 
 Thus, setting
 $$E_k :=  \Big ( T_k \cap  \big ( \big [ - \frac{\ell _k}{4N} , \frac{\ell _k}{4N} \big ] \times [ 0 , + \infty )  \big )  \Big  ) \setminus \omega _k \,,$$
  we have 
 \begin{equation}\label{f:misura} \liminf _k |E_k| \ge \liminf _k \frac{|T_k|}{16N^2}= \frac{m}{16 N ^ 2 } >0
 \end{equation} 
(where the first inequality holds by a proportion argument, which works since the side of length $2 \ell _k$ was assumed to be the longest one of $T _k$). 
Then we estimate $P_{r_k}(\Om_k)$ as follows:  
 $$P_{r_k}(\Om_k)\ge \int _{E_k}\Big ( | \Delta _{r_k, \frac{r_k}{2}} | -  4 r_k ^ 2 \displaystyle   \sum _{  i = 0 }^{2 q +1} |\theta_k ^i |\Big)dx  $$
In view of \eqref{f:misura}, recalling that $ | \Delta _{r_k, \frac{r_k}{2}} | = \Big ( \frac{\pi}{3} - \frac{\sqrt 3}{4} \Big )r_k^2$,  and taking into account that, 
by \eqref{f:convangles}, $\theta_i^k \ra 0$  for every $i= 0, \dots, 2q+1$, 
we conclude that
 \begin{equation}
 \label{f:ord3}
 \liminf_{k \to +\infty} \frac{P_{r_k}(\Om_k)}{r_k^3} =+\infty\,, 
 \end{equation} 
in contradiction with \eqref{f:stimareg} and \eqref{f:absurd2}. 
 
 If $  {\mathcal H}^1(\gamma_k)  \to + \infty$, then we estimate $P_{r_k}(\Om_k)$ as follows:
 $$P_{r_k}(\Om_k)  \ge  {\mathcal H}^1(\gamma_k) \frac{r_k}{2} \Big ( | \Delta _{r_k, \frac{r_k}{2}} | -  4 r_k ^ 2 \displaystyle   \sum _{  i = 0 }^{2 q +1} |\theta_k^i|\Big),$$
 so that \eqref{f:ord3} is again valid, in contradiction with \eqref{f:stimareg} and \eqref{f:absurd2}.

\bigskip
{\it Step 2: Identification of local concentrations.}   
Finally, only situations (b) and (c) can occur. Since any sequence of triangles in situation (b) does not affect the   limit of the sequence $\{\chi_{\Om_k}\}$ in $L ^ 1 _{loc} (\R ^2)$, in order to describe the geometry of local concentrations, we focus only on the sequences of triangles in situation  (c). For any pair of such sequences $\{T_k^i \}, \{ T_k^j \}$, we consider the corresponding sequences  of vectors 
$\{y_k^i \}, \{y_k^j\}$ such that \eqref{f:convtriang} holds, and we look at whether the distances  $\|y_k^i-y_k^j\|$  remain bounded or diverge as  $k \to +\infty$. This way we define an equivalence relation on the family of sequences of triangles in situation (c), which  splits them into a finite number $p$ of equivalence classes. By construction,  for $i = 1, \dots, p$, there exist  sequences of vectors $\{ \widetilde y _k ^ i \}$,  with $\|\widetilde y_k^i- \widetilde y_k^j\| \to + \infty$ for $i \neq j$,    such that 
\begin{equation}\label{f:tildey}
\Om_k-\widetilde y_k^i \ \stackrel{H^c_{loc}}{\longrightarrow}\  \Om^i _\sharp \,,
\end{equation} 
where 
$\Om ^1 _ \sharp, \dots, \Omega^ p_ \sharp$ are generalized polygons, with a total number of sides   not larger than $N$ and total area equal $\pi$, i.e., 
$$\sum_{i=1}^ p |\Om^ i_\sharp| =\pi\,.$$ 
Then we consider the open sets with polygonal boundary obtained as $ \Om_0  ^i :=  {\rm Int} (\overline {\Om^ i_ \sharp})$; we observe that 
\begin{equation}
 \label{bufr23.03}
 \sum _{i=1}^p |\partial  \Om^i _0|\ge |\partial \Om_N^*|,
  \end{equation} 
with equality if and  only if $p=1$ and $  \Om^1 _0=\Om_N^*$. 
Indeed, any open connected component $U$ of set $\Om^i_0$ has a boundary which is union of closed  polygonal lines, each one with at most $N$ edges. Then, by removing every bounded connected component of $\R^2 \setminus U$   and rescaling the set thus obtained by a factor less than $1$, it is possible to decrease the perimeter by preserving the area. Then the classical polygonal isoperimetric inequality ensures  that   $\sum _{i=1}^p |\partial  \Om^i _0|$ is not smaller than the sum of the perimeters of $p$ regular $N$-gons with total area $\pi$, and  \eqref{bufr23.03} follows from the sub-additivity of the map $\R^+\ni t \to \sqrt{t}$.

 \bigskip
{\it Step 3: We prove that }   
 \begin{equation}
 \label{bufr23.02}
 \liminf_{k \to +\infty} \frac   {P_{r_k}(\Om_k)}{r_k^3} \ge \frac 23 \sum_{i=1}^p   |\partial \Om_0^i|.
 \end{equation} 
The above lower bound, combined with the upper bound  inequality \eqref{f:stimareg} and with the assumption \eqref{f:absurd2},  will imply that \eqref{bufr23.03} holds with equality sign. This implies in particular that, for $k$ large enough, the sets $\Om_k$ must be convex. As we have seen in Section \ref{sec:convex}, this contradicts \eqref{f:absurd2}. 

Let us prove  \eqref{bufr23.02}. We fix an index $i \in \{1, \dots, p\}$ and we localize our estimates around the set $\Om ^i_0$; we may also  assume without loss of generality that the  corresponding 
vectors $\widetilde y ^i _k$ in \eqref{f:tildey} are equal to zero.  
Choosing $R_i>0$ such that $\overline {\Om^ i_\sharp}\subset B_{R_i}(0)$,  from \eqref{f:tildey} we have 
$$\Om_k \cap  B_{R_i}(0)\ \stackrel{H^c }{\longrightarrow}\ \Om^ i_ \sharp \,.$$  
Dropping the index $i$ for simplicity of notation, we have to show that 
\begin{equation}
 \label{bufr23.04}
 \liminf_{k \to +\infty} \frac {1}{r _k ^ {  3}  }   {\int_{B_{R}(0) \cap \Om_k} \!\int _{\Om_k^c }\chi_{B_{r_k}(0)} (x-y) dx dy}  \ge \frac 23   |\partial  \Om_0|.
 \end{equation}

We focus our analysis around a fixed side of $\Om_0$. Its endpoints are limit of vertices of $\Om_k$, but its interior as well may contain some accumulation points of vertices of $\Om_k$. These accumulation points divide  our side into several segments (at most $N$). We pick one of them, say $[- \ell, \ell]\times \{0\}$, with $\Om _0$ lying below the segment. From the $H^c$-convergence, and since the  open segment  $(- \ell, \ell)\times \{0\}$ does not contain any accumulation point of vertices of $\Omega _k$, if $\vps \in (0, \frac{\ell}{2})$ and $\delta >0$ are sufficiently small so that $[-\ell,  \ell]\times [-\delta, \delta] \subset B _ R (0)$, then inside the rectangle $[-\ell+\vps,  \ell-\vps]\times [-\delta, \delta]$  the structure of $\Om_k$ is similar to the one of the set $\Sigma$ in Lemma \ref{l:strips} (cf.\ \eqref{f:sigma}).

Precisely, we consider  the sides of $ \Om_k$
which intersect the rectangle $[-\ell +\vps, \ell-\vps]\times [-\delta, \delta]$,
and whose supporting lines $\Gamma_{\theta_k ^i, s_k^i}$ satisfy   $s_k ^i \ra 0$ as $k \to + \infty$. Assume that, as the index $i$ goes from $1$ to $2q+1$, 
 those lines are labelled from the bottom to the top. 
Choosing $\delta' >0$  such that $[-\ell+2\vps, \ell-2\vps]\times [-\delta, \delta']$ does not meet any other side of $\Om_k$, we 
 can locally represent $\Omega _k$ in $ [- \ell+2\vps, \ell-2\vps]\times [-\delta, \delta ']$ as the following union of strips: 
$$\Sigma_k= \Big ([- \ell+2\vps, \ell-2\vps]\times [-\delta, \delta']\Big ) \cap \Big (\Gamma^{-}_{\theta^{1}_k, s^ {1} _k}\cup \bigcup_{i=1}^q \Gamma^{+}_{\theta^{2i}_k, s^{2i} _k} \cap \Gamma^{-}_{\theta^ {2i+1}_k, s^{2i+1} _k}\Big ).$$
Then, using Lemma \ref{l:strips} (and Remark \ref{bufr23.06}), we get
$$ \begin{array}{ll} 
  &  \displaystyle \liminf_{k \to +\infty} \frac{1}{r_k ^ 3} {\int_{\Sigma_k}\int _{\Om_k^c }\chi_{B_{r_k}(0)} (x-y) dx dy}  \ge  \\ 
   \noalign{\bigskip}  
 &  \displaystyle \liminf_{k \to +\infty} \frac{1}{r_k^3}{\int_{- \ell+2\vps }^{\ell-2\vps} dx_1 \int _{\{x_2 \in \R \,: \, s_x\in [0,r_k]\}} |\Sigma_k^c \cap B_{r_k} (x) |dx_2} \ge  \\ 

 \noalign{\bigskip} 
 & \displaystyle \liminf_{k \to +\infty} \frac{1}{r_k^3}  {\int_{- \ell+2\vps }^{\ell-2\vps} dx_1 \int _0^{r_k} \big ( | \Delta _{r_k,  x_2 } | -  4 r_k ^ 2 \displaystyle   \sum _{  i = 0 }^{2 q +1} |\theta_k^i|\big)  dx_2} \ge 
 \\ \noalign {\bigskip} &  \displaystyle \frac 23   (2 \ell-4\vps).
 \end{array} 
$$
Inequality \eqref{bufr23.04} follows by repeating the above  argument around each side of $ \Om_0$ and letting $\vps \ra 0$.  

\qed 

\bigskip 
\section {Proofs of Theorem \ref{t:false} and of Corollary \ref{c:corpotenze}  } \label{sec:proof3}

\bigskip 
{\bf Proof of Theorem \ref{t:false}}.  For brevity, let us denote by $J _r$ the functional 
$$J _ r (\Omega) = \int _\Om \int _\Om \chi _{B _ r (0)}  (x-y) \, dx \, dy\,,$$ 
so that 
$$P _ r(\Omega) = |\Omega | |B _ r (0)| - J _ r (\Omega) = \pi ^ 2 r ^ 2 - J _ r (\Omega) \qquad \forall \Omega \in \mathcal P _N\,.$$ 
Clearly, for $\Om \in \mathcal P _N$, we have
$$ \begin{cases}J_ r (\Omega) =
\pi ^ 2 & \text{ if  } {\rm diam} \, \Om \leq r  
\\ 
J_ r (\Omega) <
\pi ^ 2 & \text{ if  } {\rm diam} \, \Om > r   
\end{cases}
$$ 
Setting $r'': = \min  \{ \rm diam \, \Om \, :\,  \Om \in \mathcal P _N \}$, 
by Reinhardt's Theorem \cite{R22}, for $N\geq 6$ even there exists a polygon $\Om^\sharp _N$, which is {\it not} a regular $N$-gon, such that 
$$r'' ={\rm diam }\,  \Om ^\sharp _N<  \rm diam \,  \Om ^* _N\,. $$
Therefore, for every $r \geq r'' $, we have
$$\pi ^ 2 = J _ r (\Om ^\sharp _N) > J _ r (\Om) \qquad \forall \Om \in \mathcal P _N \text{ with } {\rm diam}\, \Om > r  $$ 
(in particular, the above strict inequality holds for $\Om = \Om ^* _N$ if ${\rm diam} \, \Om ^* _N >r$). 
It follows that, for every $r\geq r''$, the maximum of $J _ r$ over $\mathcal P _N$ equals $\pi ^ 2$, 
or equivalently the minimum of $P _ r$ over  $\mathcal P _N$ equals $\pi ^ 2 (r ^ 2 -1)$, 
and they are attained at $\Om ^\sharp _N$.  \qed

%

\bigskip
{\bf Proof of Corollary \ref{c:corpotenze}}
  (i)  Let $r'' $ and $\Om ^\sharp _N$ as in the statement of Theorem \ref{t:false}. For every $ r \in [r'',  {\rm diam}\, \Om ^* _N )$, there exists a polygon $\Om^\sharp _N \in \mathcal P _N$ such that
$$J _ r (\Om ^ \sharp _N ) > J _ r (\Om^*_N)\,.$$
Then, it is enough to consider a sequence of 
non-negative and non-increasing radially symmetric kernels  $\{ h _ n \}$ in  $L ^ 1 _{\rm loc} (\R ^2)$
which converge increasingly to $\chi _ {B _r}(0)$. By the monotone convergence theorem, for $n$ large enough we have
$$J _ {h_n} (\Om ^ \sharp _N ) > J _ {h_n} (\Om^*_N)\,.$$

(ii)   Again, let $r''$ and $\Om ^ \sharp _N$ be as in the statement of Theorem \ref{t:false}. Set 
$\e:= [{\rm diam} (\Omega ^* _N ) - {\rm diam} (\Omega _N ^ \sharp)  ]/3$.  
If $[x_0, y _0]$ is a diameter of $\Omega ^*_N$, for $x \in B _ \e (x_0)$ and $y \in B _\e (y _0)$ we have 
$$|x-y| \geq |x_0 - y _0| - |x- x_0| - |y - y _0| \geq {\rm diam} (\Omega ^*_N ) - 2 \e = {\rm diam} (\Omega _N ^ \sharp) + \e =: r'' + \e \,.$$
Hence, 
\begin{equation}\label{f:doppio2}
\begin{array}{ll} 
\displaystyle \int_{\Omega ^*_N  } \int _{\Omega ^*_N}  |x-y|^k\, dx \, dy  & \displaystyle \geq  \int _{B _\e (x_0)\cap \Om ^*_N }\int _{B _\e (y_0)\cap \Om ^*_N }  |x-y|^k \, dx \, dy 
\\ \noalign{\bigskip} 
& \displaystyle \geq | B _\e (x_0)\cap \Om  ^*_N|\,  | B _\e (y_0)\cap \Om ^*_N |  (r'' + \e ) ^ k\,.   
\end{array}
\end{equation} 
On the other hand, we have
\begin{equation}\label{f:doppio1}  \int_{\Omega _N ^\sharp } \int _{\Omega _N ^\sharp}  |x-y|^k dx \, dy  \leq \pi^2 \, (r '') ^k
\end{equation} 
By comparing \eqref{f:doppio2} and \eqref{f:doppio1} we infer that, for $k$ large enough, 
$$ \int_{\Omega^*_N  } \int _{\Omega^*_N}  |x-y|^k\, dx \, dy  > \int_{\Omega _N ^\sharp } \int _{\Omega _N ^\sharp}  |x-y|^k \,.$$ 

We now examine in particular the case $N = 6$. Denoting by $H_R$ the regular hexagon with unit diameter and by $H_G$ the Graham hexagon with the same area, it follows that $\diam(H_G)=d\approx 0.980957$.  Then, a direct estimate using $|x-y|\leq d$ gives 
\[ \int_{H_G}\! \int _{ H_G} |x-y|^k dx \, dy \leq |H_G|^2 d^k.\]
On the other hand, by arguing as done above to obtain \eqref{f:doppio2}, we get
\[ \int_{H_R} \! \int_{H _ R}  |x-y|^k dx\, dy \geq 3|B_\varepsilon(x_0)\cap H_R||B_\varepsilon(y_0)\cap H_R|(d+\varepsilon)^k.\]
Observing that $|B_\varepsilon(x_0)\cap H_R|=\frac{\pi \varepsilon^2}{3}$, $\varepsilon=\frac{1-d}{3}$, and using the numerical value of $d$ shows that, for any $k \geq 2832$, $H _ G$ has a lower energy than  $H _ R$. 
\qed

\bigskip

\section{ Proof of Theorem \ref{t:polygons}}  \label{sec:proof4} 
We argue in two steps. 

\smallskip
{\it Step 1.} 
We claim that, 
if $\theta _i$ denotes  inner angle of $\Om$ at the vertex $A_i$, 
it holds $\theta _i = \theta _{i +1}$. 
We prove this claim by contradiction, via a reflection argument. Assume $\theta _{i+1} > \theta _i$. 
Let $H$ be the  {symmetry} axis of the side $[A_i,A_{i+1}]$, let $M_i$ be their intersection point and, for $x \in [A_i, M_i]$, let $x'$ be  its symmetric about $H$.  {Given $\omega$ in the plane we denote $\omega_\star$ its reflection about $H$.}
We have
\begin{equation}\label{f:sectors}
(\Om \cap B _ r (x))_\star \subseteq \Om \cap B _ r (x')\,, \text{ with strict inclusion for } x \text{ close to } A _{i}.
\end{equation} 
Indeed, the strict inclusion for $x$ close to $A_{i}$ readily follows from 
$$(\Om \cap B _ r (x))_\star = \Om_\star \cap B _ r (x') \qquad \text{ and } \qquad \theta _{i+1} > \theta _i \,.$$

We assert that the inclusion \eqref{f:sectors} remains valid for all $x \in  [A_i, M_i]$, namely that 
\begin{equation}\label{f:key} \Om_\star \cap B _ r (x')  \subseteq \Om \cap B _ r (x') \qquad \forall x \in  [A_i, M_i]\,.
\end{equation}
Once proved \eqref{f:key}, the  contradiction required to achieve the proof of Step 1 readily follows. Indeed, recalling that the inclusion becomes strict for $x$ close to $A_{i}$, we have 
$$\begin{array}{ll}
\displaystyle  \int _{A_i}^ {M_i } \!\!\!\!  v_\Om (x) \,  |x M_i| \, dx   
&<
\displaystyle  \int _{A_i}^ {M_i } \!\!\!\!   v_\Om (x')  \,  |x' M_i| \, dx 
\\ \noalign{\bigskip} 
& = 
\displaystyle  \int _{M_i}^ {A_{i+1}} \!\!\!\!  v_\Om (x')  \,  |x' M_i| \, dx'\,,
\end{array} 
$$
 {in contradiction with} \eqref{f:angles}. 

It remains to show \eqref{f:key}.  To that purpose, we distinguish the cases $N = 3$ and $N >3$.

\begin{figure} [h] 
\centering   
\includegraphics[height=0.39\textwidth]{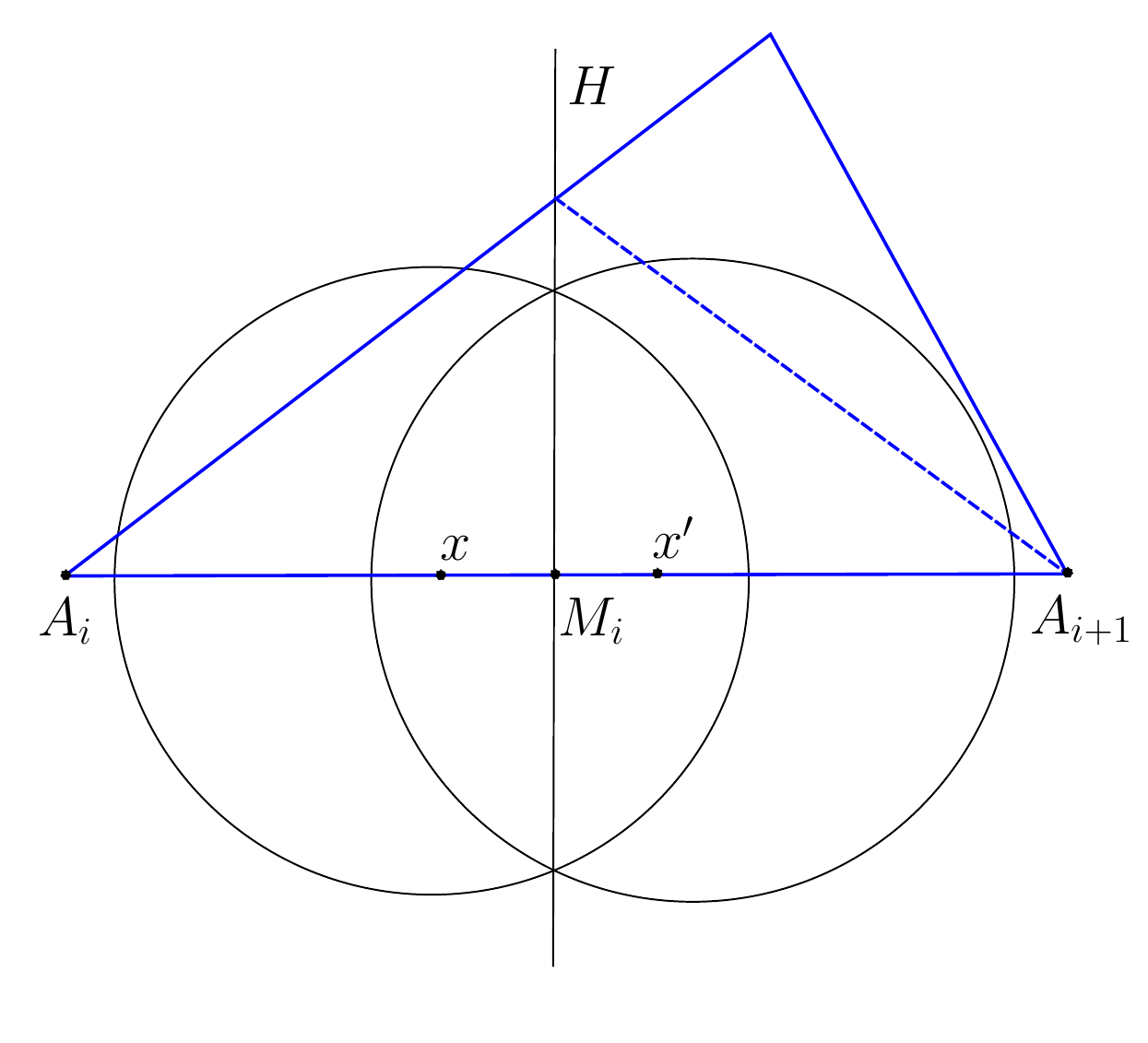}\quad
\includegraphics[height=0.39\textwidth]{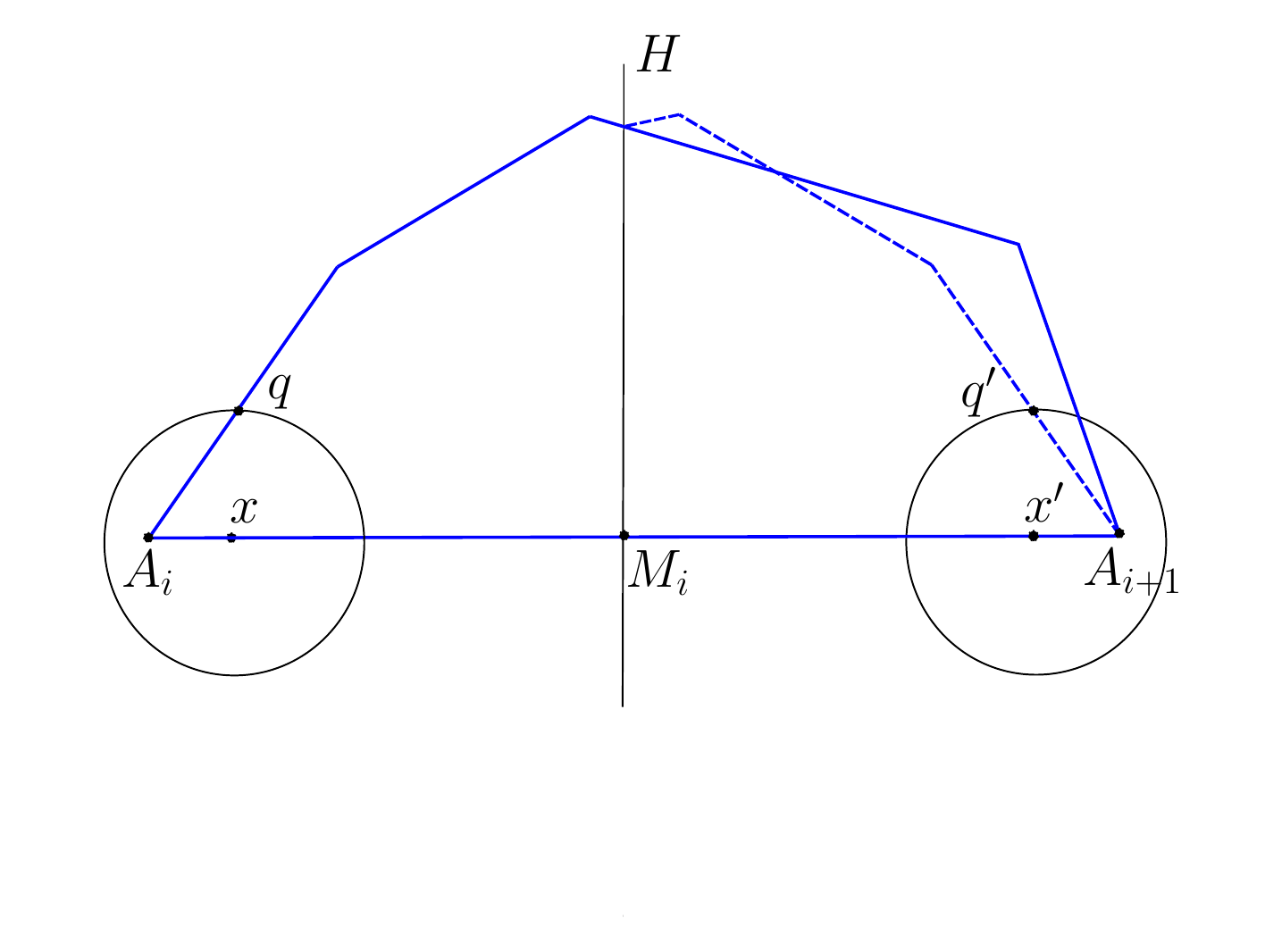}
\caption{The reflection argument for $N$-gons, with $N = 3$  (on the left) and  $N>3$ (on the right)}
\label{fig3-4}   
\end{figure}

\medskip
For triangles, \eqref{f:key}  is straightforward, since the inequality $\theta _{i+1} > \theta _i$ implies  that $\Om _ \star \subseteq \Om $, 
 see Figure \ref{fig3-4}, left.

For $N >3$, the inequality $\theta _{i+1} > \theta _i$ does not imply, in general, that $\Om _ \star \subseteq \Om $, see Figure \ref{fig3-4}, right. Nevertheless, thanks to assumption \eqref{f:card}, in order to check that \eqref{f:key} still holds, it is enough to show that,
for any 
$q \in [A_{i-1}, A _{i}] \cap   B _ r (x)$, 
if $q'$ denotes its symmetric about $H$, the trapezoid $\mathcal T$ with vertices $q, q', A_i, A _{i+1}$ is contained into  $\Om$. 
Again thanks to assumption \eqref{f:card}, the inclusion $\mathcal T \subset \Om$ holds true as soon as $q' \in \overline \Om$. This latter property is true because $q \in [A_{i-1}, A _i]$, and $\theta _{i+1} >\theta _ i$.

\bigskip
{\it Step 2.} 
We claim that, 
if $\ell _i$ denotes the length of the side $[A_i, A _{i+1}]$, 
it holds $\ell _i = \ell _{i +1}$.  
To prove this claim we are going to assume without loss of generality that $N \geq 4$. Indeed, if $\Om$ is a triangle, Step 1 proved above already gives that $\Om$ is equilateral. 

We point out that, by assumption \eqref{f:card}, it holds
\begin{equation} \label{f:length}
r \leq \min _{ i = 1, \dots, N } \frac{ \ell _i }{2} . 
\end{equation} 
Indeed,  in case $r >  \frac{ \ell _i }{2} $ for some index $i$, $\partial \Om \cap B _ r (M _i)$,  with $M _i$ mid-point of $[A_i, A _{i+1}]$ would not be contained into two consecutive sides of $\Om$, against \eqref{f:card}. 

Next we recall that, by Step 1, all the inner angles $\theta _i$ of $\Om$ are equal to a fixed angle $\theta$. Since $N \geq 4$, we have $\theta \geq \pi/2$. 
Therefore, points $x\in [A_i, A _{i+1}]$  such that $B _ r (x)$ intersects  another side of $\Om$ (which by \eqref{f:card} is necessarily a consecutive side) 
are points $x$ whose distance from $A_i$ or from $A_{i+1}$ does not exceed $r$. For such points,  
$v_\Om (x) $ is a function depending on  $|xA_i |$ (or $|xA_{i+1}|$), $r$, and $\theta $, but not on $\ell _i$. 

Taking  \eqref{f:length} into account, we conclude that 
for a suitable function $f$, it holds 
$$\int _{A_i}  ^ { A _{i+1}} v_\Om ( x)  \, dx = 
\frac{ 1}{2}  {\|h \| _{L^1 }} \big ( \ell _i - 2r)  + f ( r, \theta ) \qquad \forall i = 1, \dots, N\,, $$ 
 { where $\|h\|_{L^1}$ is the $L^1$ norm of $h$ in $\R^2$.} We now enforce condition \eqref{f:sides} to deduce 
$$\frac{ 1}{2} \|h \| _{L^1 } \big ( \ell _i - 2r)  + f ( r, \theta ) = c \ell _i  \qquad \forall i = 1, \dots, N\,. $$  
Clearly this system can be satisfied only if either all the $\ell _i$'s are equal, or $c = \frac{ 1}{2} \|h \| _{L^1 } $. 
But the latter equality cannot hold: indeed, by \eqref{f:sides}, $c$ is the integral mean over a side of $\Om$
of the function $v_\Om (x) $, and such function is always less than or equal to $\frac{ 1}{2} \|h \| _{L^1 }$, with strict inequality near the vertices. \qed

\section{Proof of Theorem \ref{t:powers} and Lemma \ref{l:specific}} \label{sec:proof2} 

{\bf Proof of Theorem \ref{t:powers}.}
Taking $h ( x) = |x| ^2$, for every $\Om \in \mathcal P _N$ we have
 $$ \begin{array}{ll}  \displaystyle \int_\Om\! \int_\Om |x-y| ^ 2 \, dxdy   & \displaystyle = \int_\Om\!\int_\Om \Big ( x_1 ^ 2 + x _ 2 ^ 2 + y _ 1 ^ 2 + y _ 2 ^ 2 - 2 x _ 1 y _ 1 - 2 x _ 2 y _2 \Big ) \, dx dy
\\ \noalign{\bigskip} 
 & \displaystyle = 2 \pi \int_\Om \big ( x_1 ^ 2 + x _ 2 ^ 2  \big  ) \, dx  - 2  \Big (  \int_\Om x _ 1  \, dx  \Big ) ^ 2  - 2 \Big ( \int _\Om x _ 2  \, dx  \Big ) ^ 2 \,.
 \end{array} 
$$ 
Since the energy is invariant by translations, it is not restrictive to assume that $\Om$ has its baricenter at the origin, and hence the result follows from Theorem \ref{t:hardy}.  

Taking $h ( x) = |x| ^4$, for every $\Om \in \mathcal P _N$ we have 
$$ \begin{array}{ll}  \displaystyle \int_\Om\! \int_\Om |x-y| ^ 4 \, dxdy   & \displaystyle = \int_\Om\! \int_\Om \Big ( x_1 ^ 2 + x _ 2 ^ 2 + y _ 1 ^ 2 + y _ 2 ^ 2 - 2 x _ 1 y _ 1 - 2 x _ 2 y _2 \Big ) \, dx dy
\\ \noalign{\bigskip} 
 & \displaystyle = \int_\Om \! \int_\Om \big [ ( x_1 ^ 2 + x _ 2 ^ 2 ) ^ 2 
 + ( y_1 ^ 2 +y  _ 2 ^ 2 ) ^ 2  + 4 (x_1 y_1 + x_2 y _ 2)  ^ 2  + 2 ( x_ 1 ^ 2 + x_ 2^2  )  (y _ 1 ^ 2 + y _ 2 ^ 2) 
\\ \noalign{\bigskip} 
 & \displaystyle \qquad \qquad  
  -4  ( x_ 1 y _ 1 + x _ 2 y _ 2) (x_1 ^ 2 + x _ 2 ^ 2)  
 -4  ( x_ 1 y _ 1 + x _ 2 y _ 2) (y_1 ^ 2 + y _ 2 ^ 2)  
   \big  ] \, dx dy
   \\ \noalign{\bigskip} 
   & \displaystyle =  2 \pi \int_\Om |x|  ^ 4  + 
6 \Big ( \int _\Om  x_ 1 ^ 2 \, dx \Big ) ^ 2  + 6 \Big ( \int _\Om  x_ 2 ^ 2 \, dx \Big ) ^ 2  
 \\ \noalign{\bigskip} 
      & \displaystyle + 8 \Big ( \int _\Om x_ 1 x_2 \, dx  \Big ) ^ 2 + 4 \Big ( \int _\Om  x_ 1 ^ 2 \, dx \Big )  \Big ( \int _\Om  x_ 2 ^ 2 \, dx \Big )  \, , 
 \end{array} 
$$
where the last equality holds up to assuming as done above that $\Om$ has its baricenter at the origin. 

Now we use the inequality
$6a^2+6b^2+4ab \geq 4(a+b)^2$, 
with equality if and only if $a=b$. Applying it for $a = \int_\Omega x_1^2, b = \int_\Omega x_2^2$,  we get
\begin{align*}
\int_{\Omega}\!  \int_{\Om}  |x-y|^4  \, dx dy & \geq  2 \pi  \int_\Omega |x|^4+ 4 \left(\int_\Omega |x|^2\right)^2
 + 8  \left(\int_\Omega x_1x_2\right)^2\\ 
& \geq 2 \pi \int_{\Omega_N^*} |x|^4+4\left(\int_{\Omega_N^*}|x|^2 \right)^2+ 8  \left(\int_{\Omega_N^*} x_1x_2\right)^2\\
& = \int_{\Omega_N^*  }\! \int _{ \Omega_N^*} |x-y|^4  \, dx dy \,,
\end{align*}
where the second inequality follows from Theorem \ref{t:hardy}, and last equality holds since 
\begin{equation}\label{f:formulas}
\int_{\Omega_N^*} x_1x_2=0 \qquad \text{  and} \qquad  \int_{\Omega_N^*}x_1^2=\int_{\Omega_N^*}x_2^2 \,.
\end{equation} 
To check these equalities, we decompose a regular $N$-gon  $\Omega_N^*$  centred at the origin into $N$ triangles $\{OA_1A_2,...,OA_NA_1\}$, with $A _i  = (p_i, q_i) =  (\cos(\alpha+i \frac{2\pi}{N}), \sin(\alpha+i \frac{2\pi}{N}))$.

Concerning the 
 first equality in \eqref{f:formulas}, on each triangle  we have 
 $$\int_{OA_iA_{i+1} } x_1 x_2  = \frac{1}{24}(p_i q_{i+1} -p_{i+1} q_i)(2p_i q_i+2p_{i+1}q_{i+1}+p_i q_{i+1}+p_{i+1}q_i) \,;$$ 
inserting the expressions of $(p_i, q_i)$ and summing over $i$ we get 
\[ \int_{\Omega_N^*} x_1 x_2 = \frac{1}{24}\sin \frac{2\pi}{N} \sum_{i=1}^N \left[\sin (2\alpha+ (i-1)\frac{4\pi}{N}) +\sin(2\alpha+  i \frac{4\pi}{N})+\sin(2\alpha+i \frac {2\pi}{N} )\right]=0 \]

Concerning the 
second equality in \eqref{f:formulas}, on each triangle we have 
\[ \int_{{OA_iA_{i+1} }} x_1^2 = \frac{1}{12} (p_i^2+p_ip_{i+1}+p_{i+1}^2)(p_iq_{i+1}-p_{i+1}q_i),\]
\[ \int_{{OA_iA_{i+1} }} x_2^2 = \frac{1}{12}  (q_i^2+q_iq_{i+1}+q_{i+1}^2)(p_iq_{i+1}-p_{i+1}q_i).\]

Inserting the expressions of $(p_i, q_i)$, we see that 
$(p_iq_{i+1}-p_{i+1}q_i) = \sin \frac{2\pi}{N}$ for every $i$, so that
\[ \int_{\Omega_N^*}  x_1^2 = \frac{1}{12}  \sin \frac{ 2\pi}{N} \sum _{i = 1} ^ N (p_i^2+p_iq_{i+1}+p_{i+1}^2) \]
\[ \int_{\Omega_N^*}  x_2^2 = \frac{1}{12}  \sin \frac{2\pi}{N} \sum _{i = 1} ^ N  (q_i^2+q_iq_{i+1}+q_{i+1}^2)  \]

Now,
\[ \sum_{i=1}^N (p_i^2+p_ip_{i+1}+p_{i+1}^2) = 2 \sum_{i=1}^N \cos^2 \big(\alpha+ i \frac{2\pi}{N}\big )+\sum_{i=1}^N \cos \big  (\alpha+i \frac{2\pi}{N}\big ) \cos(\alpha+(i+1)\frac{2\pi}{N})\]
\[ \sum_{i=1}^N (q_i^2+q_iq_{i+1}+q_{i+1}^2) = 2 \sum_{i=1}^N \sin^2 \big (\alpha+i \frac{2\pi}{N}\big )+\sum_{i=1}^N \sin\big  (\alpha+i \frac{2\pi}{N}\big ) \sin\big (\alpha+(i+1)\frac{2\pi}{N}\big )\]
Using the formulas
\[\cos^2 a = \frac{1+\cos 2a}{2}, \qquad 2\cos a \cos b = \cos (a+b)+\cos(a-b)\]
\[\sin^2 a = \frac{1-\cos 2a }{2}, \qquad 2\sin a \sin b = \cos (a-b)-\cos(a+b)\]
we conclude that
\[\sum_{i=1}^N \cos^2 \big (\alpha+i \frac{2\pi}{N}\big ) = \sum_{i=1}^N \sin^2 \big (\alpha+i \frac{2\pi}{N}\big )  = \frac{N}{2}\] 
\[\sum_{i=1}^N\cos\big  (\alpha+i \frac{2\pi}{N}\big ) \cos\big (\alpha+(i+1)\frac{2\pi}{N}\big ) =\sum_{i=1}^N \sin \big (\alpha+i \frac{2\pi}{N}\big ) \sin\big (\alpha+(i+1)\frac{2\pi}{N}\big ) =  N\cos  \frac{2\pi}{N}.\]

\qed 

\bigskip 
{\bf Proof of Lemma \ref{l:specific}}. 
 
(i) If $\Omega \in \mathcal P _8$ is axially symmetric, 
some elementary computations give
$$\begin{array}{ll} \displaystyle \int_\Om \!\int_\Om |x-y| ^ 6 \, dxdy   & \displaystyle = 2  \pi \int _\Om |x| ^ 6 \, dx + 18 \int _\Om |x| ^ 2 \, dx \int _\Om |x | ^ 4 \, dx \, 
\\ \noalign{\bigskip} 
& \displaystyle
+ 12  \int_\Om (x_ 1 ^ 2 - x _ 2 ^ 2) \, dx 
 \int_\Om (x_ 1 ^ 4 - x _ 2 ^ 4)   \, dx 
\\ \noalign{\bigskip} 
& \displaystyle
\geq 2 \pi  \int _{\Om _8 ^*} |x| ^ 6 \, dx + 18 \int _{\Om_8^*} |x| ^ 2 \, dx \int _{\Om_8^*} |x | ^ 4 \, dx
\\ \noalign{\bigskip} 
& \displaystyle
+ 12  \int_\Om (x_ 1 ^ 2 - x _ 2 ^ 2) \, dx 
 \int_\Om (x_ 1 ^ 4 - x _ 2 ^ 4)   \, dx \,,
   \end{array}$$ 
where the inequality is obtained by invoking as usual Theorem \ref{t:hardy}. Since
$$ \int_{\Om_8 ^*}  (x_ 1 ^ 2 - x _ 2 ^ 2) \, dx 
 \int_{\Om_8 ^*}  (x_ 1 ^ 4 - x _ 2 ^ 4) = 0 \, , $$
 to conclude the proof it is enough to show that 
$$\int_\Om (x_ 1 ^ 2 - x _ 2 ^ 2) \, dx 
 \int_\Om (x_ 1 ^ 4 - x _ 2 ^ 4)   \, dx \geq 0 \,.$$
 To that aim,  we exploit the assumption that $\Om$ is an axisymmetric convex octagon. 
Assuming without loss of generality that $\Om$ has four vertices at the points $(\pm1, 0)$, $(0, \pm 1)$, 
 and one in the region $\{ x _1>0, \ x_2 > 0, \ x _ 2 \geq x _ 1, \ x _ 2 \geq 1 - x _ 1\} $, we have that $\Om$ contains an axisymmetric convex octagon $Q$, 
 which has still four vertices at the points $(\pm1, 0)$, $(0, \pm 1)$ and one on the straight line $x_1 = x _2$. Then
$$\int_\Om (x_ 1 ^ 2 - x _ 2 ^ 2) \, dx 
 \int_\Om (x_ 1 ^ 4 - x _ 2 ^ 4)   \, dx \geq \int_Q (x_ 1 ^ 2 - x _ 2 ^ 2) \, dx 
 \int_Q (x_ 1 ^ 4 - x _ 2 ^ 4)   \, dx  = 0  \,.$$

(ii) We can assume without loss of generality that $k>1$. Indeed,  once the inequality is proved for such $k$, one can pass to the limit as $k \ra 1$. 
The proof is inspired from \cite{Laug21}. In the remaining of this proof, the functional $J _h$ with $h ( x) = |x| ^ k$ will be denoted for brevity by $J _k$, i.e. we set 
$$J _ k ( \Om):= \int _ \Om \!\int _\Om | x-y| ^ k \, dx \, dx\,.  $$
Our target is to show that, for every real $2 \times 2$ volume preserving real matrix $M$, it holds 
$$J _ k ( M( \Om_N^*)) \geq J _ k (\Om _N ^*) \,. $$
We have $M = ASB$, with $S$ diagonal and $A, B \in O (2)$. Since $J _ k (\Om ) = J _k ( A (\Om))$ for any $\Om \in \mathcal P _N$ and any $A \in O (2)$, it is not restrictive
to assume that $A = {\rm Id}$. Moreover, by considering the regular polygon $B ( \Om _N ^*)$ in place of $ \Om _N ^*$, we may assume that also $B = {\rm Id}$. 
We are thus reduced to show that 
$$J _ k ( S( \Om_N^*)) \geq J _ k (\Om _N ^*) \,,\qquad \text{ where } S = {\rm diag } \{ \sigma \,, \,  \sigma ^ { -1} \}, \ \sigma \in \R ^+ \,. $$

We consider the family of polygons $\Om _ t := S _ t ( \Om_N^*) $, with $S _ t = {\rm diag } \{ \sigma^ t \,, \,  \sigma ^ { -t} \}$, so that
$\Om _ 0 =  \Om_N^*$, and $\Om _ 1 = S( \Om_N^*)$.  We  claim that
$J _ k  ( \Om_N^*  )  > J _k ( \Om _ t)$ for every  $t \in (0, 1]$, or equivalently that 
the map  $g_k (t):= J _k ( \Om _ t) $  satisfies 
$
g _k( t) >g_k ( 0)$ for every  $t \in (0, 1]$.
Indeed,  let us show that 
\begin{equation}\label{f:derig} 
g_k' (0) = 0\, ,  \qquad  g_k ''(t) \geq 0 \text{ on } (0, 1)\,.
\end{equation} 

By  arguing as in \cite[Lemma 4.3 and Lemma 4.4]{Laug21},  we are allowed to differentiate under the sign of integral. 
Setting $\Phi _k ( r ) = r ^ k$, and writing for brevity $r$ in place of $|S_ t ( x) - S _ t ( y)|$, by direct computations we have 
$$g_k' (t)   = \int _ \Om \int _ \Om \frac{\partial }{\partial t }  
\Phi _k ( |S _ t x  - S _ t y | ) \, dx \, dy \,, \qquad  g_k'' (t)   = \int _ \Om \int _ \Om \frac{\partial ^2}{\partial t ^ 2}  
\Phi _k ( |S _ t x  - S _ t y | ) \, dx \, dy $$ 
where

$$\begin{array}{ll} \displaystyle \frac{\partial }{\partial t }  
\Phi _k ( |S _ t x  - S _ t y | )  
 & \displaystyle= \Phi' ( r)  \ \frac{S_t ( x) - S_ t (y)}{|S_ t ( x) - S_ t (y) | } \cdot ( \dot S _ t ( x) - \dot S _ t (y))  
 \\ \noalign{\bigskip} 
 & \displaystyle= k r ^ { k-1}  \ \frac{S_t ( x) - S_ t (y)}{|S_ t ( x) - S_ t (y) | } \cdot ( \dot S _ t ( x) - \dot S _ t (y)) 
 \end{array}  
$$ 
and

$$\begin{array}{ll} \displaystyle \frac{\partial ^2}{\partial t ^ 2}  
\Phi _k ( |S _ t x  - S _ t y | )  
& \displaystyle = \Big ( \Phi _ k '' (r) - \frac{\Phi _k ' ( r) } {r} \Big ) \Big ( \frac{S_t ( x) - S_ t (y)}{|S_ t ( x) - S_ t (y) | } \cdot ( \dot S _ t ( x) - \dot S _ t (y))  \Big ) ^2
 \\ \noalign{\bigskip} 
& \displaystyle 
+ \frac{\Phi_k ' ( r) }{r} \Big ( |\dot S_t (x) -\dot S _ t (y) | ^ 2 + ( S_ t(x) - S _ t (y)) \cdot ( \ddot S _ t ( x) - \ddot S _ t (y))  \Big ) 
 \\ \noalign{\bigskip} 
& \displaystyle =  k ^ 2 r ^ { k-2}  \Big ( \frac{S_t ( x) - S_ t (y)}{|S_ t ( x) - S_ t (y) | } \cdot ( \dot S _ t ( x) - \dot S _ t (y))  \Big ) ^2
 \\ \noalign{\bigskip} 
& \displaystyle 
+ k r ^ { k-2}  \Big ( |\dot S_t (x) -\dot S _ t (y) | ^ 2 + (\log \sigma ) ^ 2  |  S_ t(x) - S _ t (y) | ^ 2  \Big ) \,.
\end{array} 
$$

Hence conditions \eqref{f:derig} are satisfied (we exploit here the assumption $k>1$). 

\qed

\bigskip

\section {Proof of Theorem \ref{t:hardy}}\label{sec:proof1} 

We consider the maximization problem
\begin{equation}\label{f:pb_r} 
 \max \Big \{ \int _\Om  h ( x) \, dx  \ :\  \Om \in \mathcal P _N \Big \} \,,  \end{equation}  
and we proceed as follows: 

\smallskip
-- In Section \ref{sec1} we prove that it is not restrictive to take $h= \chi _{B_r (0)}$ (cf.\ Proposition \ref{p:reduction}).

\smallskip

-- In Section \ref{sec2}  we prove that  problem \eqref{f:pb_r} admits a solution, which is  a classical  star-shaped  polygon (cf.\ Proposition \ref{p:existence}). 

\smallskip

-- In Section \ref{sec3}  we prove that, among classical star-shaped polygons, the regular $N$-gon is optimal (cf.\ Proposition \ref{p:regular}).

\smallskip
The validity of Theorem \ref{t:hardy} follows at once by combining Propositions \ref{p:reduction}, \ref{p:existence}, and \ref{p:regular}. 

\smallskip
In the sequel, we write for brevity $B _r$ in place of $B _ r (0)$.

\subsection{Reduction to the characteristic kernel}\label{sec1} 

\begin{proposition}\label{p:reduction} 
If a polygon solves problem \eqref{f:pb_r} when $h = \chi _ {B _ r}$ for all $r>0$, then it solves problem    \eqref{f:pb_r} for every admissible kernel $h$. 
%
%
  \end{proposition}

\proof  {Assume $\Omega^*$ solves \eqref{f:pb_r} for every $r>0$. Consider a sequence of radii $0<r_1<...<r_k$ and some positive reals $d_1,...,d_k$. Then by optimality of $\Omega^*$ for every radius $r_i$, $i=1,...,k$ we have 
\[ \int_{\Omega} \sum_{i=1}^k d_i \chi_{B_{r_i}} \leq \int_{\Omega^*} \sum_{i=1}^k d_i \chi_{B_{r_i}}\]

Since every radially decreasing step function can be written in the form $\sum_{i=1}^k d_i \chi_{B_{r_i}}$ it follows that $\Omega^*$ solves \eqref{f:pb_r} for radially decreasing step functions.

Every radially decreasing function $h$ can be written as the limit of an increasing sequence of radial step functions $\{h_n\}_{n\geq 1}$. Passing to the limit in the inequalities
\[ \int_\Omega h_n(x)dx \leq \int_{\Omega^*} h_n(x)dx,\]
shows that $\Omega^*$ solves \eqref{f:pb_r} for arbitrary admissible kernels. 
}
\qed

\subsection{Existence and reduction to  star-shaped polygons}\label{sec2}
  In view of Proposition \ref{p:reduction}, we are going to focus our attention on the maximization problem
\begin{equation}\label{f:pb_rr}
\max  \Big \{ \mathcal E (\Om): = | \Omega \cap B _r| \ : \ \Om \text{ is a polygon with $N$ sides with }    |\Om| \leq \pi \Big \}\,.
\end{equation}

\begin{proposition}\label{p:existence}
Problem \eqref{f:pb_rr} 
 admits a solution. Moreover, every solution is   star-shaped.
\end{proposition} 

As a preliminary remark, let us 
observe that there are some ranges for the value of $r$ for which
problem \eqref{f:pb_rr}  can be elementarily solved.  More specifically, 
let $\Om _N ^{r, {circ}}$  and $\Om _N ^{r, {in}}$   denote respectively
 the regular $N$-gons  circumscribed and inscribed to $B _r$.   Then: 
\begin{itemize} 
\item[--]   If $\pi \geq | \Om _N ^{r, {circ}} |$,   
 the maximum in \eqref{f:pb_rr}  is equal to $\pi r ^2$, and 
 it is attained either at infinitely many admissible polygons 
 among which the regular $N$-gon of area $\pi$  (if the inequality is strict),  
 or uniquely at $\Om _N ^{r, {circ}}$ (if the inequality holds with equality sign). 

 \item[--]  If $\pi \leq | \Om _N ^{r, {in}} |$,   
 the maximum in \eqref{f:pb_rr}  is equal to $\pi$, and 
 it is attained either at infinitely many admissible polygons 
 among which the regular $N$-gon of area $\pi$  (if the inequality is strict),  
 or uniquely at $\Om _N ^{r, {in}}$  (if the inequality holds with equality sign).  
  \end{itemize} 

Thus, in the remaining of this section we always tacitly assume that  $r$ is chosen so that 
\begin{equation}\label{f:assu} | \Om _N ^{r, {in}} | <  \pi <  | \Om _N ^{r, {circ}} | 
\end{equation}
 (this is just for definiteness, as our proof below works also in the `trivial' cases left). 
 
\smallskip
The proof of Proposition \ref{p:existence}  requires as a key ingredient  a geometric construction  that we state separately in the next lemma, along with its application in our problem in the subsequent remark.
%
\begin{lemma}\label{l:geometric}  
Let $\Om $ be a polygon with  $N$ sides such that $|\Omega \cap B _ r| >0 $ and $|\Omega \setminus B _ r|>0$. Then there exists another polygon $\Om '$ with $N$ sides, 
which is star-shaped and satisfies the inclusions
 \begin{equation}\label{f:inclusions}  
  (\Omega ' \cap B _ r) \supseteq (\Omega \cap B _r) \qquad \text{ and } \qquad ( \Omega' \setminus B _ r) \subseteq  (\Omega \setminus B _ r )\,.
  \end{equation} \end{lemma}

\begin{remark}\label{r:procedure} Let $\Om$ and $\Om'$ be  polygons as in Lemma \ref{l:geometric}.  We claim that, starting from $\Omega'$, it is easy to construct another star-shaped polygon  $\widetilde \Om$, still having $N$ sides, such that 
\begin{equation}\label{f:tilde}  |\widetilde \Om| \leq |\Om| \qquad \text{ and } \qquad |\widetilde \Om \cap B _ r| \geq |\Om \cap B _ r| \,.\end{equation}   
Indeed, from the first inclusion in \eqref{f:inclusions} we have
$\E  ( \Omega '   )   = \E (\Omega ) + \delta$,  with  $\delta >0$.  
Compare then the areas of $\Omega' $ and $\Omega$.  In case $|\Omega ' | \leq |\Omega|$, we simply define
 $\widetilde \Omega := \Omega  ' $. 
 In the case left, namely when $|\Omega '  | > |\Omega|$, we define $\widetilde \Omega $ as the polygon homothetic to $\Omega' $ which as the same area as $\Omega$, namely we take
 $\widetilde \Omega  :=  ({|\Omega |}/{|\Omega' |} ) ^{1/2}  \Omega'$.
 Clearly,  such $\widetilde \Omega $ is a still a star-shaped polygon with $N$ sides, and  it is easy to check that that  its energy is not less than the energy of $\Omega $. Actually we have: 
 $$  \E (\widetilde \Omega  ) \geq \frac{|\Omega |}{|\Omega' |} \E  ( \Omega '  )   =\frac{|\Omega  |}{|\Omega' |}   \big ( \E (\Omega  ) + \delta \big ) 
 = \frac{\E(\Omega  ) + |\Omega \setminus B _ r|}{  \E (\Omega  ) + \delta + |\Omega'   \setminus B _r|}   \big ( \E (\Omega  ) + \delta \big )  \geq \E (\Omega ) \,,
   $$ where the last inequality follows after an immediate computation as a consequence of the two inequalities
 $\delta  >0$ and  $ |\Omega  ' \setminus B _ r| \leq |\Omega  \setminus   B _ r |$, 
the latter holding by the second inclusion in \eqref{f:inclusions}. 
\end{remark}

Let us assume for a moment that Lemma \ref{l:geometric} holds true, and let us show how  Proposition \ref{p:existence} follows.  Let  $\{ \Omega _n \}$ be a maximizing sequence for problem \eqref{f:pb_rr}.
Clearly, up to a subsequence each polygon $\Omega _n$ has a non-negligible intersection both with $B _r$ and with its complement (recall we are assuming \eqref{f:assu}). Then,   
 for every $n$,   we denote by $\Om_n'$ the star-shaped polygon given by Lemma \ref{l:geometric}. 
Proceeding as in Remark \ref{r:procedure}, 
we obtain a new polygon $\widetilde \Om_n$ with its star-shapedness centre inside  $B_r$, which satisfies \eqref{f:tilde}. Thus, $\{ \widetilde \Om _n \}$ is still a maximizing sequence.     \EEE 
%

Now, by 
the compactness and lower semicontinuity properties of the Hausdorff complementary topology \cite[Section 2]{HP},  
up to 
a subsequence,  we may assume that 
\begin{equation}\label{locconv} {\widetilde \Om_n}  \sr{H ^ c _{\rm loc}}{\lra} \Om\,
\end{equation}
for some $\Om$ which is a generalized polygon with $N$ sides, according to Definition \ref{d:genpol}, having area at most $\pi$.
We observe that the perimeter of the sets $( \widetilde \Om_n \cap B_r)$ is uniformly bounded  from above (since all the polygons $ \widetilde \Omega_n$ have a fixed number $N$ of sides). Hence, by the compact embedding of  $BV ( B _ r)$ into $L ^ 1 ( B _ r)$, the energy $\E( \widetilde \Om_n) $ converges to $\E ( \Omega ) $. 

Notice carefully that the limit generalized polygon $\Om$ is still star-shaped (this is precisely the  scope reached through the  modification of the sequence $\Om_n$ into the sequence $\widetilde \Om_n$).

It may still occur that $\partial \Omega$ contains some self-intersections, but they can only be  contact segments between two consecutive sides. 
Hence,  it is enough to remove any such contact segment, in order to transform $\Om$ 
into a classical polygon, which will be a solution to problem \eqref{f:pb_rr}.  \qed

\bigskip  We now turn to the most delicate part of the proof, namely  the geometric construction in Lemma \ref{l:geometric}.

\bigskip
{\bf Proof of Lemma \ref{l:geometric}}. 
Let $\Om  \in   {\mathcal P _N}$ be a polygon  as in the assumptions of the Lemma. 
To prove the statement, we can further assume with no loss of generality 
that  $\Omega$ has no side tangent to $B_r$ and 
no vertex in $\partial B_r$. Indeed, if this is not the case, 
once the Lemma is proved for polygons with no side tangent to $B _r$ and no vertex in $\partial B _ r$, 
we can approximate $\Omega$  (in the Hausdorff complementary topology) by a sequence of polygons $\{\Omega _n\}$ satisfying such additional conditions, 
and apply the Lemma to each $\Omega _n$: we find a sequence of polygons $\{\Omega ' _n\}$, whose limit polygon $\Omega'$ 
(which exists up to passing to a subsequence  and is still star-shaped) does the job for $\Omega$. 

Thus, let $ \Omega \in   {\mathcal P _N}$ be a polygon as in the assumptions of the Lemma, which 
in addition has no side tangent to $B _ r$, and no vertex in $\partial B _r$.  
For the sake of clearness, we give first the construction of the polygon $\Om'$  in a simplified situation, 
namely when the intersection between $\Omega$ and $\partial B _r$ consists precisely of $N$ arcs of circle, and then 
 we proceed in the general case.

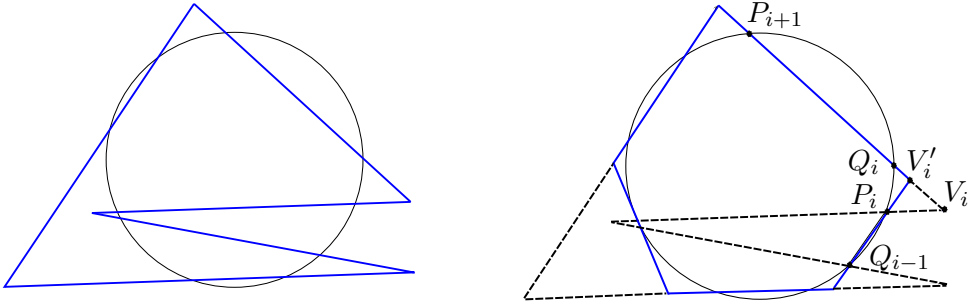
\begin{figure} [h] 
\centering   
\def\svgwidth{7cm}   
\begingroup%
  \makeatletter%
  \providecommand\color[2][]{%
    \errmessage{(Inkscape) Color is used for the text in Inkscape, but the package 'color.sty' is not loaded}%
    \renewcommand\color[2][]{}%
  }%
  \providecommand\transparent[1]{%
    \errmessage{(Inkscape) Transparency is used (non-zero) for the text in Inkscape, but the package 'transparent.sty' is not loaded}%
    \renewcommand\transparent[1]{}%
  }%
  \providecommand\rotatebox[2]{#2}%
  \ifx\svgwidth\undefined%
    \setlength{\unitlength}{265.53bp}%
    \ifx\svgscale\undefined%
      \relax%
    \else%
      \setlength{\unitlength}{\unitlength * \real{\svgscale}}%
    \fi%
  \else%
    \setlength{\unitlength}{\svgwidth}%
  \fi%
  \global\let\svgwidth\undefined%
  \global\let\svgscale\undefined%
  \makeatother%
  \begin{picture}
  (4,1)%
    \put(0.1,0){ \includegraphics[height=4.5cm]{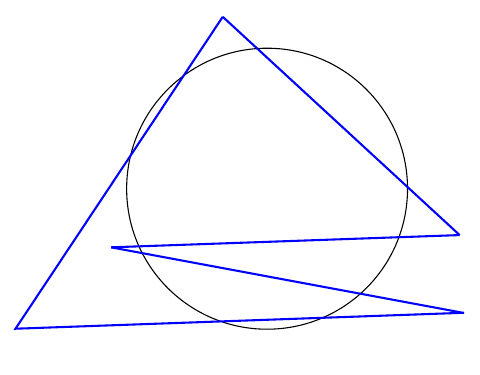} \hskip 1 cm  \includegraphics[height=4.5cm]{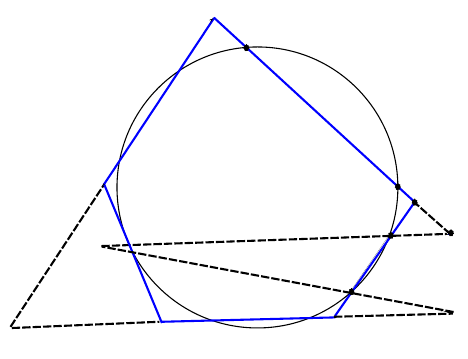}} 
     \put(1.54, 0.58){\color[rgb]{0,0,0}\makebox(0,0)[lb]{\smash{$P_{i+1}$}}}%
    \put(1.735, 0.235) {\color[rgb]{0,0,0}\makebox(0,0)[lb]{\smash{$P_i$}}}%
    \put(1.91, 0.24) {\color[rgb]{0,0,0}\makebox(0,0)[lb]{\smash{$V_i$}}}%
    \put(1.84, 0.3) {\color[rgb]{0,0,0}\makebox(0,0)[lb]{\smash{$V'_i$}}}%
    \put(1.73, 0.3) {\color[rgb]{0,0,0}\makebox(0,0)[lb]{\smash{$Q_i$}}}%
        \put(1.77, 0.12) {\color[rgb]{0,0,0}\makebox(0,0)[lb]{\smash{$Q_{i-1}$}}}%
         \end{picture}%
\endgroup%
\caption{Costruction of the polygon $\Omega'$ in Lemma \ref{l:geometric}, case when each side has  endpoints outside $\overline B _ r$ and meets  $B _r$ }
\label{fig:hardysimple}   
\end{figure}

\smallskip
$\bullet$ {\it Case when each side of $\Omega$ has both its endpoints outside $\overline B _ r$, and intersects  $B _r$}.  
 (Equivalently,  the intersection between $\Omega$ and $\partial B _r$ consists precisely of $N$ arcs of circle.) Starting from a fixed endpoint of such an arc, say $P _1$, and following a  counter-clockwise oriented parametrization of $\partial B _r$,   name these arcs  $\wideparen{P_i Q_i}$,  for  $i = 1, \dots, N$; none of these arcs is degenerated into a point, since by assumption no side of $\Omega$ is tangent to $B _r$, see Figure \ref{fig:hardysimple}, left.  
 
 For $i = 1, \dots N$, let $\gamma _i$ be the straight line through $Q_i$ and $P _{i+1}$,  with the convention $P _{N+1} = P _ 1$, and let $\pi _i$ be  the  (open) half-plane determined by $\gamma _i$ which contains all the points $Q_k, P_j$  for $k \neq i$ and $j \neq {i+1}$.  
 
 \begin{equation}\label{f:easy} \Omega'  := \bigcap _{ i = 1, \dots, N} \pi _i \,.
 \end{equation} 
By construction, $\Omega'$ is a classical convex polygon in $\mathcal P _N$, see Figure \ref{fig:hardysimple}, right.

\vskip -2.5cm 
 
 The first inclusion in  \eqref{f:inclusions} is satisfied because any circular segment delimited by the arc $\wideparen{P_{i+1} Q_i}$ and the segment $P_{i+1} Q_i$ cannot intersect $\Om$ (otherwise the arc $\wideparen{P_{i+1} Q_i}$ would be crossed by a side of $\Om$). 

The second inclusion in \eqref{f:inclusions} is satisfied because,   
denoting by $V_i$ the common vertex of the two consecutive sides of $\Om$ containing $P_i$ and $Q_i$, and by $V'_i$ the intersection of the straight lines $\gamma _i$ (through $Q_i$ and $P _{i+1}$) and $\gamma _{i-1}$  
(through $Q_{i-1}$ and $P_i$), it holds 
$${    \Delta \, Q_i V_i' P_i } \ \subseteq \  { \Delta \, Q_i V_i P_i } \,.  $$  
This is due to the fact that   the point $V'_i$ belongs to both the half-plane determined by the straight line trhough $V_i$ and $P_i$ and containing $Q_i$, and  the half-plane determined by  $V_i$ and $Q_i$ and containing $P_i$.    
%
  
\smallskip
$\bullet$ {\it General case}.
 Consider the intersection between $\Omega$ and $\partial B _ r$. Such intersection consists now of  $M$ arcs of circle, with $M \leq N$. Starting from a fixed endpoint of such an arc, say $P _1$, and following a   counter-clockwise   parametrization of $\partial B _r$,   name these arcs  $\wideparen{P_i Q_i}$,  for  $i = 1, \dots, M$.  
  Notice that, if we  equip $\partial \Omega$   with an oriented parametrization such that $\Om$ lies on the left of each side, then,  at every point $Q_i$, the side of $\Om$ passing through $Q_i$ is entering into $B_r$.  
 
  For $i = 1, \dots M$, let $\pi _ i$ be the half-planes defined as above.   Let also $\{A _1, \dots, A _k\}$ denote the (possibly empty) family of vertices of $\Omega$ lying inside $B _r$ (recall that by assumption no vertex of $\Omega$ lies on $\partial B _r)$.   

We point out that, if one would define $\Omega '$ as in \eqref{f:easy}, none of the two inclusions in \eqref{f:easy} would be in general satisfied: the former due to the possible presence of vertices inside $B _r$, the latter due to the possible presence of sides exterior to $B _r$ (in both cases, with possible self-intersections occuring in $\partial \Omega$). 

For this reason, the definition of $\Omega'$ is more involved: 
we are going to construct it as the union of two sets, denoted by $\Omega ' _{in}$ and $\Omega '  _{out}$, which lie respectively inside and outside $B _r$. Such sets are ``curvilinear'' polygons, whose boundaries do not contain self-intersections, and consist in a finite number of  ``sides'', meant as arcs of circle lying on $\partial B _r$ or line segments (in case of $\Omega' _{in}$, the segments lie inside $B _r$, while in case of $\Omega ' _{out} $ they lie outside). The closures of the two sets $\Omega ' _{in}$ and $\Omega '  _{out}$ intersect precisely at the $M$ arcs $\wideparen{P_i Q_i}$, so that the set
  \begin{equation}\label{f:defgen}\Omega ' := \Omega' _{in} \cup \Omega ' _{out} \cup  \big \{  \wideparen{P_i Q_i} \ :\ i = 1, \dots, M \} \, \end{equation}
turns out to be a classical polygon, which by construction will be a star-shaped one. 

Let us specify how $ \Omega' _{in}$ and  $\Omega ' _{out}$ are defined. We set
 $$\Omega' _{in} := {\rm conv} \Big ( \bigcap _{ i = 1, \dots, M} \pi _i   , A_1, \dots , A _k \Big )\cap B _r  \,.$$ 
 where ${\rm conv}$ denotes the convex envelope.  By 
   construction, specifically thanks to the presence of the vertices $\{A _1, \dots, A _k\}$ in the above definition, we have
 \begin{equation}\label{f:hard1}
 \Omega' _{in}  \supseteq ( \Omega \cap B _ r)\,.
 \end{equation}

 In order to define $\Omega' _{out}$,  let us choose a point in the interior of $\Omega'_{in}$, say $x_0$, which does not belong to any of the straight lines supporting the edges of $\Om$.
 
 Let us denote by $\mathcal F$ the family of  all the straight lines supporting the non-circular edges of  $\Omega'_{in}$.
   For every $i = 1, \dots, M$, we introduce the set
  $$ \Delta_{i}   := \Big \{ x + t \overrightarrow {x_0 x}   \ :\ x \in  \wideparen{P_i Q_i}  , \ t \in (0, \lambda (x) ) \Big \}$$ 
where  
  $$\lambda ( x)  := \inf \Big \{ t >0 \ :\ x + t\overrightarrow {x_0 x}  \in (\partial \Omega \cup \mathcal F )\Big \}\,. $$

 We define 
 $$ \Omega' _{out} :=   \bigcup_{i=1}^M\Delta_{i}  \,.$$ 
 
 \vskip - 1 cm 
   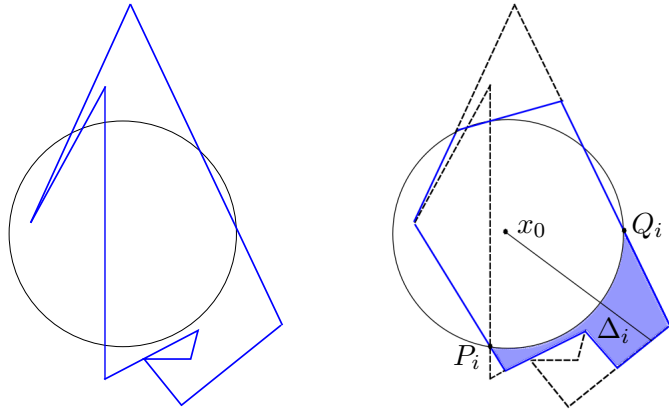
\begin{figure} [h] 
\centering   
\def\svgwidth{7cm}   
\begingroup%
  \makeatletter%
  \providecommand\color[2][]{%
    \errmessage{(Inkscape) Color is used for the text in Inkscape, but the package 'color.sty' is not loaded}%
    \renewcommand\color[2][]{}%
  }%
  \providecommand\transparent[1]{%
    \errmessage{(Inkscape) Transparency is used (non-zero) for the text in Inkscape, but the package 'transparent.sty' is not loaded}%
    \renewcommand\transparent[1]{}%
  }%
  \providecommand\rotatebox[2]{#2}%
  \ifx\svgwidth\undefined%
    \setlength{\unitlength}{265.53bp}%
    \ifx\svgscale\undefined%
      \relax%
    \else%
      \setlength{\unitlength}{\unitlength * \real{\svgscale}}%
    \fi%
  \else%
    \setlength{\unitlength}{\svgwidth}%
  \fi%
  \global\let\svgwidth\undefined%
  \global\let\svgscale\undefined%
  \makeatother%
  \begin{picture}
  (4,1)%
    \put(0.3,0){ \includegraphics[height=5.5cm]{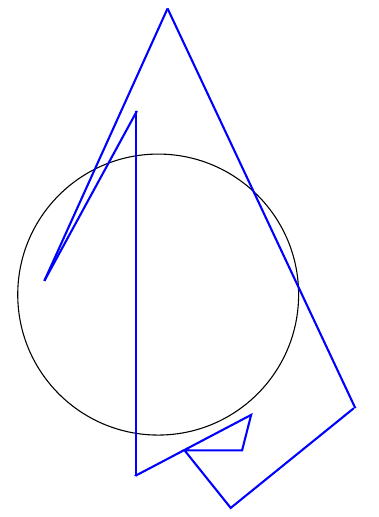} \hskip 1 cm  \includegraphics[height=5.5cm]{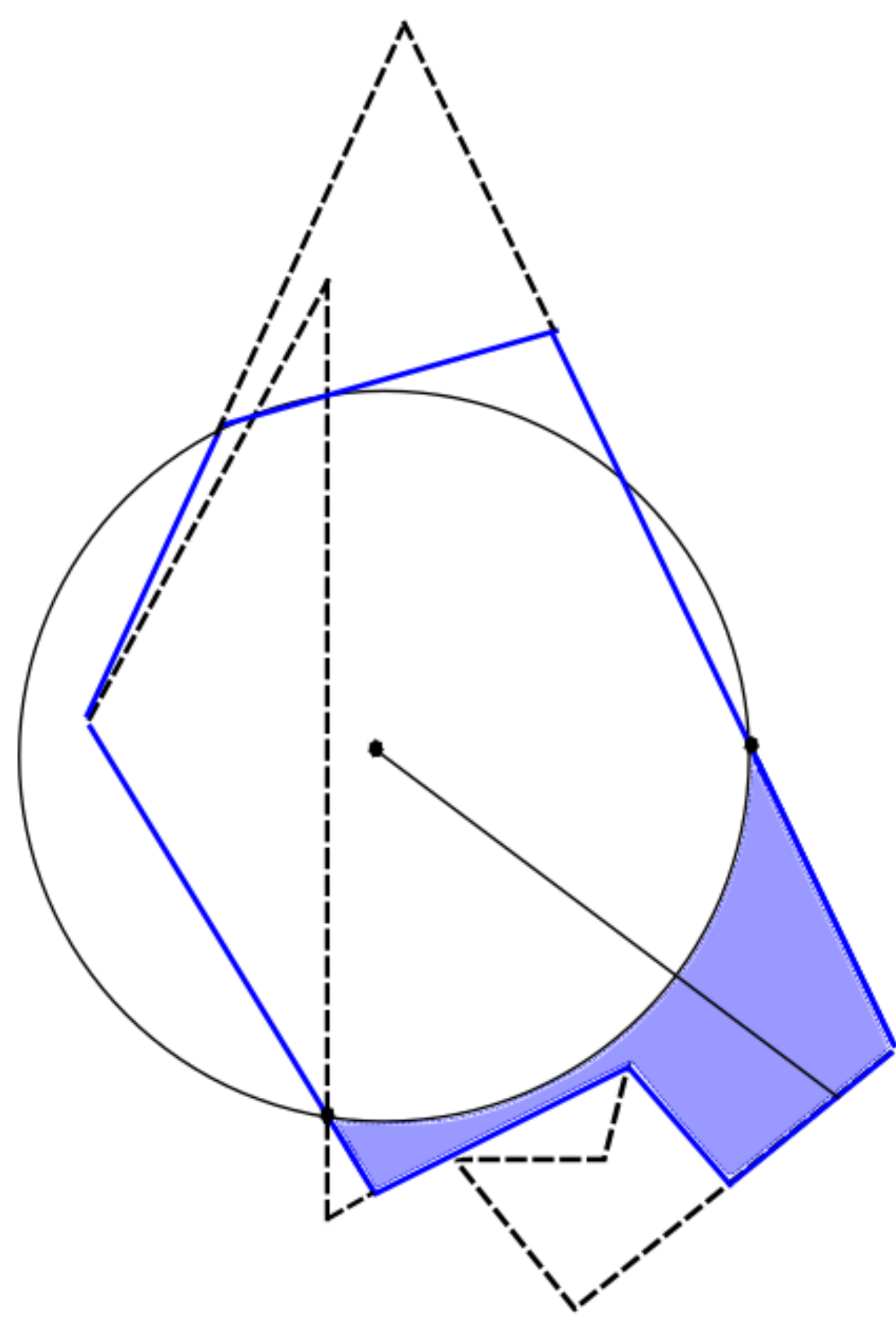}} 
    \put(1.18 , 0.09) {\color[rgb]{0,0,0}\makebox(0,0)[lb]{\smash{$P_i$}}}%
    \put(1.515, 0.34) {\color[rgb]{0,0,0}\makebox(0,0)[lb]{\smash{$Q_i$}}}%
    \put(1.3, 0.34) {\color[rgb]{0,0,0}\makebox(0,0)[lb]{\smash{$x_0$}}}%
       \put(1.45, 0.14) {\color[rgb]{0,0,0}\makebox(0,0)[lb]{\smash{$\Delta_i$}}}%
         \end{picture}%
\endgroup%
\caption{Costruction of the polygon $\Omega'$ in Lemma \ref{l:geometric}, general case}
\label{fig:hardyhard}
\end{figure}

 By  construction we have
 \begin{equation}\label{f:hard2} \Omega' _{out}  \subseteq ( \Omega \setminus B _ r)\,.
 \end{equation}
In view of \eqref{f:hard1}, and \eqref{f:hard2}, the set $\Omega'$ defined in \eqref{f:defgen} is a classical  star-shaped  polygon which satisifes both inclusions in \eqref{f:inclusions} as required.

 \smallskip
 To conclude, it remains to show that  $\Omega' \in \mathcal P _N$, namely that our procedure respects the constraint on the number of sides. We set
 $$\begin{array}{ll} 
 N _{in }& := \text{ number of sides of $\Omega$ which intersect $B _r$ } 
 \\ 
 N _{out }& := \text{ number of sides of $\Omega$ which do not intersect $B _r$\,,} 
 \end{array}
 $$
so that $ N _{in }+ N _{out }  = N$. Denoting by $\mathcal S'$ the family of sides of $\Omega'$, 
we have 
$$\mathcal S' = \mathcal S' _{in} \cup \mathcal S' _{out}\,,$$
  where $  \mathcal S' _{in}$ and $\mathcal S' _{out}$ are  the families of non-circular sides respectively of $\Omega' _{{in}}$ and of $\Omega' _{out}$.  To prove that $\Omega' \in \mathcal P _N$ (namely that ${\rm card} ( \mathcal S' ) \leq N$), we are going to show that
  \begin{equation}\label{f:counting} {\rm card} ( \mathcal S' _{in}  ) \leq N _{in} \qquad \text{ and } \qquad {\rm card} ( \mathcal S' _{out} \setminus \mathcal S' _{in}  ) \leq N _{out} \,.
  \end{equation}

 \medskip
{\it  (i) Counting inside. }
For every $i = 1, \dots, M$, let us denote by $A_i^1, \dots, A_i^{j_i}$ the vertices of $\Omega'_{in}$ lying in the interior of the circular segment delimited by the arc $\wideparen{Q_i P_{i+1}}$ and the line segment $Q_i P_{i+1}$.  Then the number of (non-circular) sides of $\Omega' _{in}$ which join $Q_i$ to $P_{i+1}$ is exactly $1+j_i$. 
To any such side, we can associate a side $\Om$ which intersects $B _r$, in the following way: 

\begin{itemize}
\item[--] to the side starting at $Q_i$, we associate the only side of $\Om$ passing through $Q_i$, which is entering into $B_r$;
\item[--]  to the side starting at $A_i ^k$, we associate the side of $\Om$ which starts at $A _i ^k$ in the positive orientation of $\partial \Omega$. 
\end{itemize} 

Since this association is injective, the first inequality in \eqref{f:counting} holds true. 

 \medskip
 {\it  (ii) Counting outside.}  
  As above, to any side in $\mathcal S' _{out} \setminus \mathcal S' _{in} $ we can associate a side of $\Om$, which in this case does not intersect $B _r$. Specifically,   we distinguish two cases: a side in $\mathcal S' _{out} \setminus \mathcal S' _{in} $ is  
  either a side of $\Om$ which does not intersect $B _r$ or a newly created one. 
  
  \begin{itemize}
  \item[--]
 to  a side  which does not intersect $B _r$, we associate the corresponding side of $\Omega$; note that the association is injective, because a side of $\Omega$  cannot be simultaneously in the boundary of two distinct sets $\Delta _{i}$'s. 
 \item[--]  to a newly created side, we associate a side of $\Omega$ which  is not part of the boundary of $\Delta _{S,i}$ in the following way: 
 a newly created side occurs when some point $A:= x + \lambda ( x) \overrightarrow  {x_0 x} $ is a vertex of $\Om$ and, for $\vps >0$ sufficiently small, we still have $x + (\lambda ( x)+\vps) \overrightarrow  {x_0 x} \in \Om$. This means that $A$ is an endpoint of a side of $\Omega$ which  is not part of the boundary of $\Delta _{i}$: this is precisely  the side we associate to the newly created one (again, with an injective association). 
   \end{itemize} 
   
   Thus also the second equality in \eqref{f:counting} holds true, and our proof is achived. \qed 
    \EEE

\subsection{Optimality of the regular $N$-gon}\label{sec3} 

\begin{proposition}\label{p:regular} 
Problem \eqref{f:pb_rr} 
 is solved by the regular $N$-gon of area $\pi$.
\end{proposition} 

\proof Let $\Om$ be a solution to problem \eqref{f:pb_rr}, which exists by Proposition \ref{p:existence}.
 As in the proof of Proposition \ref{p:existence},  
we are going to assume   with no loss of generality that  the inequalities
 \eqref{f:assu} are satisfied.  
We are going to prove the result through several claims.  We stress that in each of these claims we estimate the variation
of the area and of the ``energy'' $\E (\Om):= |\Om \cap B _ r|$ when $\Omega$ is perturbed by some kind of deformation,  
preserving the number of sides: this strategy is allowed precisely by the crucial information $\Omega$ is a {\it classical} polygon, whose boundary does not contain self-intersections. 

\smallskip
$\bullet$ {\it No side entirely outside $B _ r$}. Indeed, assume that $\Omega$ has a side $S$ which does not intersect $B_r$. By the right inequality in \eqref{f:assu}, there exists another side $S'$ which intersects $B _r$. 
We move simultaneously $S$  and $S'$, both in a  parallell way to themselves, respectively towards the interior and towards the exterior of $\Omega$:
 the area is preserved while the energy increases, contradicting optimality.  

\smallskip
$\bullet$ {\it No side with one vertex in $B _r$ and one vertex outside $\overline B _r$}. Indeed, assume that $\Omega$ as 
such a side $S$.  We perform a rotation  of $S$ around its mid-point, in such way that the vertex
of $S$ which lies inside $ B _ r$ moves towards the exterior of $\Omega$: 
 the area is preserved at first order, while the energy increases, again contradicting optimality.

\smallskip 
$\bullet$ {\it No vertex  in $B _r$}.
Assume by contradiction that some side $S$ of  $\Omega$ has an endpoint in $B _r$.
 By the previous claim,  we know that
the other endpoint of  $S$ cannot lie outside $\overline B _ r$, hence it lies either in $B _r$ or on $\partial B _ r$. 
On the other hand, by the left inequality in \eqref{f:assu}, we can exclude that all vertices of $\Omega$ lie in $B _r$.  
We deduce that necessarily $\Omega$ contains 
a chain of consecutive sides, all entirely contained into $\overline B_r$, such that  
the first and the last sides in the chain have one vertex in $\partial B _ r$ and the other one in $B _r$, while all the intermediate sides in the chain
are entirely contained into $B _ r$. 
Now, we can move to $\partial B _r$ all the vertices of the chain lying in $B _r$ 
so to construct another polygon $\Omega' \in \mathcal P _N$  which satisfies the inclusions \eqref{f:inclusions}.  
Starting from this polygon $\Omega'$ and arguing as in Remark \ref{r:procedure}, 
we find another polygon $\widetilde \Omega \in \mathcal P _N$, with $|\widetilde \Omega| \leq \pi$,  which has a strictly larger energy, contradicting optimality.


\smallskip
$\bullet$ {\it $\Om$ is inscribed into a circle, concentric with $B _r$, of radius $>r$.} 
By the previous item, we may associate with each side of $\Omega$ a chord of $B _r$, given by its intersection with ${B _r}$ (a priori possibly coinciding with the side itself). 
We perform a rotation of a fixed arbitrary side around its mid point: the first order optimality conditions  yield that
the mid-point of any side must coincide with the mid-point of the chord associated with it  (apply Lemma \ref{l:derivatives},  eq.~\eqref{f:r2} with $h = \chi _{B _r}$ and eq.~\eqref{f:r1} from the Appendix in Section \ref{sec:appendix}). We infer that all the vertices of $\Om$ have the same distance from the center of $B _r$, namely that $\Om$ is inscribed into a circle concentric with $B _r$. By the left inequality in \eqref{f:assu}, this circle has radius strictly larger than $r$.

\smallskip
$\bullet$ {\it $\Om$ is a regular polygon.}  We make a simultaneous parallel movement of two different sides in such a way to preserve the area of $\Omega$. 
Denoting by $\ell _i$ the lengths of the sides of $\Om$ and by $c_i$ the lengths of the corresponding chords (obtained by intersecting the sides with $B _ r$), 
the first order optimality conditions yield  $$\frac{c_i}{\ell _i}= \frac{c_j}{\ell _j} \qquad \text{ for every }  i \neq j \,$$
   (apply Lemma \ref{l:derivatives},  eq.~\eqref{f:p2} with $h = \chi _{B _r}$ and eq.~\eqref{f:p1}).  
 Combined with the previous item, this yields that  
 $\ell  _i =  \ell  _j$ for every $i,j$, and hence $\Om$ is a regular polygon. 

\smallskip
$\bullet$ {\it $\Om$ is the  regular $N$-gon of area $\pi$.}     
We already know from the previous steps that $\Om$ is a regular polygon, with number of vertices at most $N$ and area at most $\pi$.

If $\Om$ has a number of sides strictly less than $N$, we can add a side just by ``cutting'' a corner which lies outside $B _ r$.  In this case we obtain a new optimal polygon with an edge not intersecting $B_r$, which contradicts the first step of this proof.  Hence, the optimal polygon is a regular $N$-gon centred with $B_r$. Then, 
since regular $N$-gons centred with $B _r$ are monotone by inclusions, the optimal polygon must have the maximal admissible area, i.e.~area equal to $\pi$.    \qed

\bigskip
\section{Appendix: first and second order shape derivatives}\label{sec:appendix}

 \subsection{General formulas}\label{sec:general}
When $h$ is an integrand of class $\mathcal C ^ 2$, integral energies  on $\R^d$ such as
$$\mathcal E _h (\Om):=  \int _\Om  h (x) \, dx\, , \qquad  \text{ or }  \qquad J _h (\Om):=  \int _\Om \int _\Om h ( x-y) \, dx \, dy\,.$$ 
are twice differentiable with respect to domain perturbations. 
More precisely,  given a Lipschitz velocity field $\theta  :\R ^ d \to \R ^ d$, let $\Phi _   t ^ \theta( x)$ denote a one parameter family of diffeomorphisms from $\R ^ d$ into itself with initial velocity $\theta$, i.e.\ $\Phi _ t ( x)  =  x + \theta(x)   t + o ( t)$. 

The first and second order Fr\'echet shape derivatives of $\mathcal E _h$,  meant respectively as
$$\begin{array}{ll}
(\mathcal E _h ) ' _\theta ( \Omega) & \displaystyle  = \lim _ {t \to 0} \frac{ \mathcal E _h(\Phi _ t ^ \theta (\Omega)  ) - \mathcal E _h (\Omega)  }{t}  \\ \noalign{ \bigskip } (\mathcal E _h ) '' _{\theta, \xi }   ( \Omega) & \displaystyle  = \lim _ {t \to 0} \frac{(  \mathcal E' _h )  _\theta (\Phi _ t ^ \xi  (\Omega)  ) - (  \mathcal E' _h )  _\theta  (\Omega)  }{t}   \,,
\end{array}
$$ 
exist, and   their computation as stated in the next lemma is classical, see \cite[Theorem 5.2.2, eq.\ (5.11)]{HP}, 
\cite[Section 2]{laurain2ndDeriv}.

\begin{lemma}\label{p:variations} Assume $h$ is of class $\mathcal C ^ 2$. Then, for any  open bounded 
Lipschitz domain $\Omega \subset \R ^ d$ and any Lipschitz deformation $\theta , \xi :\R ^ d \to \R ^ d$, it holds 
\begin{equation}\label{f:ell} (\mathcal E _h ) ' _\theta ( \Omega)  = \int_\Omega \dv ( h \theta ) 
\end{equation}
and  
\begin{equation}\label{f:b}  \begin{array} {ll} \displaystyle (\mathcal E _h ) '' _{\theta, \xi }   (\Omega)  =  &\displaystyle  \int _\Omega\Big  [ 
 \nabla ^ 2 h \theta \cdot \xi + \nabla h \cdot ( \theta  \dv \xi + \xi \dv \theta) 
\\  \noalign{\medskip}   &  \quad + h \big  ( \dv \theta  \dv \xi - \frac{1}{2} (\nabla ^ T \theta  : \nabla \xi + \nabla ^ T \xi : \nabla \theta ) \big  ) 
\Big ] \,. \end{array} 
\end{equation}  
\end{lemma} 

Since the above result is valid in every space dimension, it applies in particular to energies of the type $J _h$ on $\Omega \times \Omega$. In that case, we need to consider ``doubled'' vector fields of the form $\Theta = (\theta(x),\theta(y)), \Xi = (\xi(x),\xi(y))$, see Section \ref{sec:hessian} for more details.

 \subsection{First order shape derivative under rotation/ parallel movement of a side}

We give hereafter the expressions of the first order shape derivatives for the energies 
$$\mathcal E _h (\Om):=  \int _\Om  h (x) \, dx\,  \qquad \text{ and }  \qquad J _h (\Om):=  \int _\Om \int _\Om h ( x-y) \, dx \, dy\,,$$ 
when  a polygon $\Omega$ is perturbed by  two distinct relevant deformations, preserving the number of sides,  
which have been previously  considered in   \cite{BF16} (see also \cite{BCT22, FV}).  We enclose their definitions to make the presentation self-contained.  
Let $\Om$ be a fixed polygon, and let  $S$ be a fixed side of $\Om$, with consecutive sides $S_1$ and $S_2$. 

\smallskip
(i)  The polygons $\Om _\e$ are obtained from $\Om$ by {\it rotation of the side $S$ around its mid-point} if they are obtained from $\Om$ by keeping  the other sides are fixed, and replacing  the three sides $(S, S_1, S_2)$ by the new sides $(S^ \e, S_1 ^ \e, S _ 2 ^ \e)$ described as follows 

\begin{itemize}
\item[--]  $S ^ \e$ lies on the straight-line obtained by rotating of an oriented angle $\e$, around the mid-point of $S$, the straight-line containing $S$; 

\smallskip
\item[--]  $S ^ \e _1$ and $S ^ \e _2$ lie on the same straight-line containing respectively $S_1$ and $S_2$;  

\smallskip
\item[--]   the lengths of $S ^ \e$, $S ^ \e _1$ and $S ^ \e _2$, are chosen so that the three sides are consecutive.
\end{itemize} 

\smallskip
(ii) The polygons $\Om _\e$ are obtained from $\Om$ by  {\it parallel movemenf of the side $S$} if they are obtained from $\Om$ by keeping  the other sides are fixed, and replacing  the three sides $(S, S_1, S_2)$ by the new sides $(S^ \e, S_1 ^ \e, S _ 2 ^ \e)$ described as follows

\begin{itemize} 
\item[--] $S ^ \e$ lies on the straight-line parallel to $S$ having signed distance $\e$ from $S$;

\smallskip
\item[--]  $S ^ \e _1$ and $S ^ \e _2$ lie on the same straight-line containing respectively $S_1$ and $S_2$;  

\smallskip
\item[--]  the lengths of $S ^ \e$, $S ^ \e _1$ and $S ^ \e _2$, are chosen so that the three sides are consecutive.
\end{itemize} 
 
\vskip -1cm
\begin{figure}[h]
\centering
\includegraphics[height=0.3\textwidth]{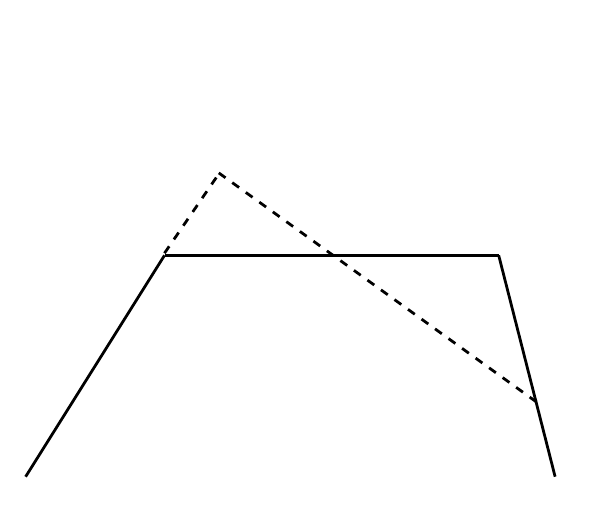} \hskip 1 cm \includegraphics[height=0.3\textwidth]{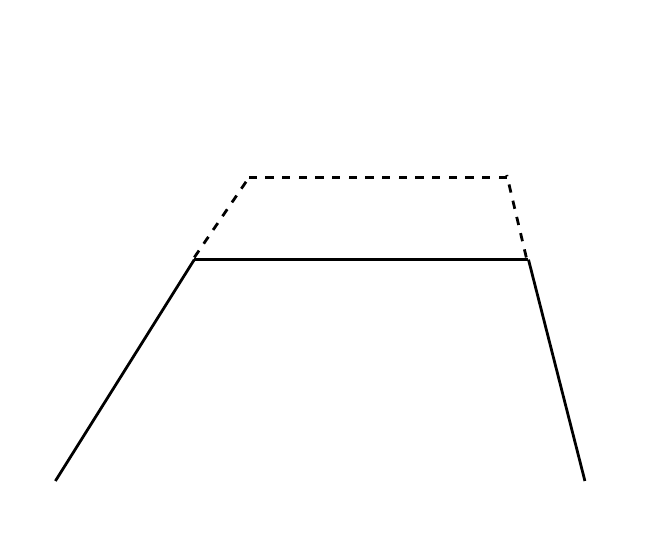} 
\label{fig:0}
\caption{Rotation around mid-point (left) and parallel movement (right).} 
\end{figure}

It follows from Lemma \ref{p:variations} that 
\begin{eqnarray}
& \displaystyle \frac{d}{d \e}  {\mathcal E} _ h (\Omega _\e) \Big | _{ \e = 0 } =  
\int _{\partial \Omega} h (x) \,  X(x)  \cdot \nu _\Omega (x)  \, d \mathcal H ^ {1} (x) & \label{f:stat0}  
\\ \noalign{\medskip}  
 & \displaystyle \frac{d}{d \e}  J _ h (\Omega _\e) \Big | _{ \e = 0 } =  2
\int _{\partial \Omega}v_\Omega (x) \,  X  (x) \cdot \nu _\Omega   (x)   \, d \mathcal H ^ {1} (x) \,, & \label{f:stat1} \end{eqnarray}

where $v_\Om (x) = \int _\Om h ( x-y  ) dy$,     $X = \frac{d} {d {\e}} \phi _{\e}  \big | _{ \e= 0}$ is the initial velocity of the deformation,  and $\nu _\Omega$ is the unit outward normal defined $\mathcal H ^ {1}$-a.e.\ on $\partial \Omega$.

 \bigskip

Now, some elementary geometric considerations show that,  in cases (i) and (ii),  the normal component of the velocity field $X$ is given respectively by
\begin{equation}\label{f:formula2} 
X(x)\cdot \nu _\Om (x)=\begin{cases}|xM| \quad \text{if}\ x\in[A_1,M]\\
-|xM| \quad \text{if}\ x\in[M,A_2]
\end{cases}, \qquad X(x)\cdot \nu _\Om(x)=1 \quad \forall x\in[A_1,A_2]\, ,
 \end{equation} 
where $A_1, A_2$ are the endpoints of $S$ and $M$ is its midpoint.

We end up with the following

\begin{lemma}\label{l:derivatives} 

(i) If $\Omega _\e$ are obtained   from $\Omega$ by  rotation of the side $S$ around its mid-point, it holds
\begin{eqnarray}
& \displaystyle \frac{d}{d \e}  {\mathcal E} _ h (\Omega _\e) \Big | _{ \e = 0 } =  
\int _{A_1} ^M  h (x)   |xM|   \, d \mathcal H ^ {1}(x) - \int _{M} ^{A_2}  h(x)    |xM|   \, d \mathcal H ^ {1} (x)  & \label{f:r2}  
\\ \noalign{\medskip}  
 & \displaystyle \frac{d}{d \e}  J _ h (\Omega _\e) \Big | _{ \e = 0 } =  2
\Big [  
\int _{A_1} ^M  v_\Omega (x)   |xM|   \, d \mathcal H ^ {1} (x) - \int _{M} ^{A_2}  v_\Omega (x)    |xM|   \, d \mathcal H ^ {1} (x) \Big ] \,. & \label{f:r3} \end{eqnarray} 
In particular, it follows from \eqref{f:r2} taking $h \equiv 1$ that
\begin{equation}\label{f:r1}
 \frac{d}{d \e}  |\Omega _\e| \Big | _{ \e = 0 } =  
0\,. \end{equation} 

(ii) If $\Omega _\e$ are obtained   from $\Omega$ by   parallel movement of the side $S$ with respect to itself, it holds
\begin{eqnarray}
& \displaystyle \frac{d}{d \e}  {\mathcal E} _ h (\Omega _\e) \Big | _{ \e = 0 } =  
\int _{S} h (x)    \, d \mathcal H ^ {1} (x)  & \label{f:p2}  
\\ \noalign{\medskip}  
 & \displaystyle \frac{d}{d \e}  J _ h (\Omega _\e) \Big | _{ \e = 0 } =  2
\int _S  v_\Omega (x)      \, d \mathcal H ^ {1} (x)  \,. & \label{f:p3} \end{eqnarray} 
In particular, it follows from \eqref{f:r2} taking $h \equiv 1$ that
\begin{equation}\label{f:p1}
 \frac{d}{d \e}  |\Omega _\e| \Big | _{ \e = 0 } =  
 \mathcal H ^ 1 (S) 
\,. \end{equation} 

\end{lemma}

\qed

\subsection{Gradient and Hessian under vertices displacement}\label{sec:hessian} 

\begin{figure}
	\centering 
	\includegraphics[width=0.35\textwidth]{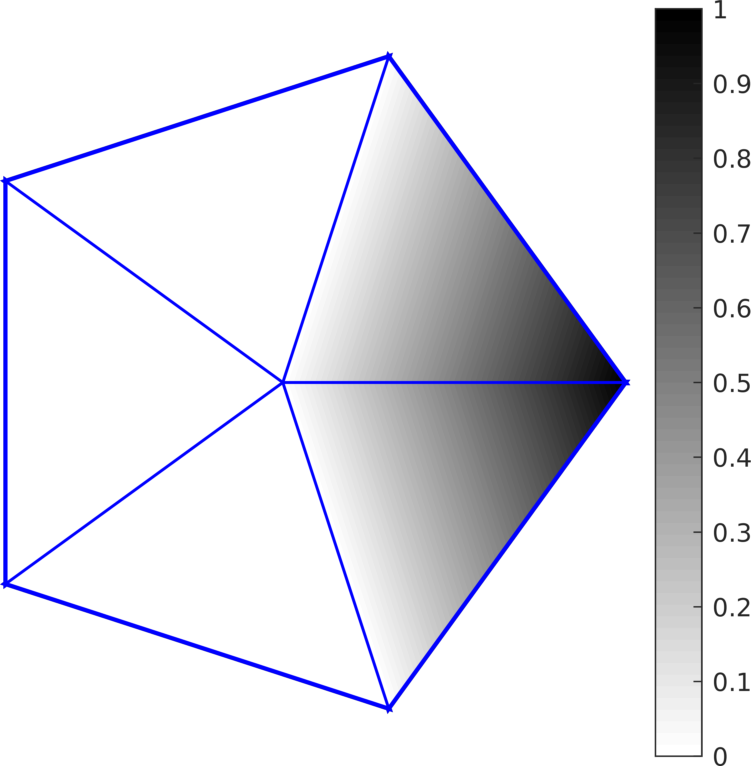}
	\caption{Triangulation of the regular polygon used for constructing polygonal deformations together with the piecewise affine function $\varphi_1$.}
	\label{fig:tripoly}
\end{figure}

Let  $\Omega$ be a $N$-gon  with vertices $A_1,..., A_N$, and let 
$\mathcal T=(T_i)_{i=1}^{M}$ be a triangulation of $\Omega$ such that the edges of $\Omega$ are edges of some triangles $T_i$. An example is shown in Figure \ref{fig:tripoly}. Following \cite{laurain2ndDeriv, BB22},  let $\varphi_i$ denote the piece-wise affine function on the triangulation $\mathcal T$ such that $\varphi_i(A_j) = \delta_{ij}$. 

Given vectors $\theta_1,\dots, \theta_N$ and $\xi _ 1, \dots, \xi _ N$ which perturb the vertices $A_1,...,A_N$, respectively, consider the double perturbation fields
$$ \Theta(x,y) = \sum_{i=1}^{N}\begin{pmatrix}\theta_i\varphi_i(x) \\ \theta_i \varphi_i(y)  \end{pmatrix}\in \Bbb{R}^2\times \Bbb{R}^2\, , \qquad 
\Xi(x,y) = \sum_{i=1}^{N}\begin{pmatrix}\xi_i\varphi_i(x) \\ \xi_i \varphi_i(y)  \end{pmatrix}\in \Bbb{R}^2\times \Bbb{R}^2.$$ 
From the general results recalled in Section \ref{sec:general}, we have
$$\begin{array}{ll} 
& \displaystyle J' _ \Theta (\Omega) = \int_{\Omega}\!  \int_{ \Omega} \nabla h \cdot \Theta + h  \di \Theta  
\\  \noalign{\medskip} 
 & \displaystyle J'' _{\Theta,\Xi} (\Omega) = \int_{\Omega} \!  \int _ {\Omega}    \nabla^2h\Theta \cdot \Xi + \nabla h \cdot(\Theta \di \Xi + \Xi \di \Theta)  + h(\di \Theta \di \Xi- \nabla\Theta^T : \nabla \Xi ) 
\end{array}
$$

Now, the different terms appearing in the above equalities can be made explicit by 
using the expressions of $\Theta, \Xi$   (the computations are similar to those in \cite[Appendix A]{laurain2ndDeriv}). 
 We have the following formulas: 
$$\di \Theta = \sum_{i=1}^N \theta_i \cdot \varphi_i(x) + \sum_{i=1}^N  \theta_i\cdot \varphi_i(y)\,;$$

\begin{align*}
\di \Xi \di \Theta & = \sum_{i,j=1}^{N} (\theta_i \cdot \nabla \varphi_i(x))(\xi_j\cdot \nabla \varphi_j(x))+ \sum_{i,j=1}^{N} (\theta_i \cdot \nabla \varphi_i(x))(\xi_j\cdot \nabla \varphi_j(y))\\& 
+ \sum_{i,j=1}^{N} (\theta_i \cdot \nabla \varphi_i(y))(\xi_j \cdot \nabla \varphi_j(x))+ \sum_{i,j=1}^{N} (\theta_i \cdot \nabla \varphi_i(y))(\xi_j \cdot\nabla \varphi_j(y))\\
& = \sum_{i,j=1}^{N}\theta_i \cdot \Big(\nabla \varphi_i(x)\otimes \nabla \varphi_j(x)+\nabla \varphi_i(x)\otimes \nabla \varphi_j(y)\\
&+\nabla \varphi_i(y)\otimes \nabla \varphi_j(x)+\nabla \varphi_i(y)\otimes \nabla \varphi_j(y)\Big) \xi_j\,;
\end{align*}

\begin{align*}
\nabla \Theta^T : \nabla \Xi & = \sum_{i,j=0}^{n-1}\theta_i \cdot (\nabla \varphi_j(x)\otimes \nabla \varphi_i(x)+\nabla \varphi_j(y)\otimes \nabla \varphi_i(y)) \xi_j\,;
\end{align*}

\begin{align*}
\nabla h \cdot(\Theta \di \Xi + \Xi \di \Theta)= \sum_{i,j=1}^N \theta_i \cdot \Big( & \varphi_i(x) ( \nabla_x h\otimes \nabla \varphi_j(x))+\varphi_j(x) (
\nabla \varphi_i(x)\otimes  \nabla_x h) \\+&\varphi_i(x) ( \nabla_x h\otimes \nabla \varphi_j(y))+\varphi_j(y) (
\nabla \varphi_i(x)\otimes  \nabla_y h)  \\
+& \varphi_i(y) ( \nabla_y h\otimes \nabla \varphi_j(x))+\varphi_j(x) (
\nabla \varphi_i(y)\otimes  \nabla_x h) \\
+& \varphi_i(y) ( \nabla_y h\otimes \nabla \varphi_j(y))+\varphi_j(y) (
\nabla \varphi_i(y)\otimes  \nabla_y h)\Big)\xi_j\,;
\end{align*}

\begin{align*}
\nabla ^2 h(x,y) \Theta \cdot \Xi = &\sum_{i,j=1}^N \theta_i\cdot \Big(   \varphi_i(x)\varphi_j(x)\nabla^2_{xx}h(x,y)+ \varphi_i(x)\varphi_j(y)\nabla^2_{xy}h(x,y)\\
+&\varphi_i(y)\varphi_j(x)\nabla^2_{yx}h(x,y)+ \varphi_i(y)\varphi_j(y)\nabla^2_{yy}h(x,y))dxdy\Big)\xi_i \,.
\end{align*} 

Inserting these formulas into the expressions of 
$J' _ \Theta (\Omega)$  and $J'' _{\Theta,\Xi} (\Omega)$,  
 we find that
 $$J' _ \Theta (\Omega) 
 = B \cdot \theta = \sum _{i = 1 }^ N  B _ i \cdot \theta _i      \qquad \text{ and } \qquad 
  J'' _{\Theta,\Xi} (\Omega)  =  \theta ^ T \cdot M \xi  = \sum _{i , j= 1} ^ N \theta _i \cdot  M _ {ij}   \xi _j\,, 
 $$
 where the vector $B = (B_i) _{ i = 1, \dots, N} $ and the matrix $ (M _{ij})_{i,j= 1, \dots, N} $, representing respectively 
  the gradient and the Hessian of $J$ with respect to the vertices, are given by
  
  \begin{equation}  B _ i =  \int_{\Omega} \int _{\Omega} (\varphi_i(x)\nabla_x h+\varphi_i(y)\nabla_y h+   h(x,y)(\nabla \varphi_i(x)+\nabla \varphi_i(y)))dxdy \label{eq:grad-F}
\end{equation}

and 

\begin{align}
M_{ij}  = \int_{\Omega}\int _{\Omega}h(x,y)(&\nabla \varphi_i(x)\otimes \nabla \varphi_j(x)-\nabla \varphi_j(x)\otimes \nabla \varphi_i(x)\notag\\
+&\nabla \varphi_i(x)\otimes \nabla \varphi_j(y)+\nabla \varphi_i(y)\otimes \nabla \varphi_j(x)\notag\\
+&\nabla \varphi_i(y)\otimes \nabla \varphi_j(y)-\nabla \varphi_j(y)\otimes \nabla \varphi_i(y))dxdy\notag\\
+\int_{\Omega} \int_{ \Omega}( &\varphi_i(x) ( \nabla_x h\otimes \nabla \varphi_j(x))+\varphi_j(x) (
\nabla \varphi_i(x)\otimes  \nabla_x h)\notag \\+&\varphi_i(x) ( \nabla_x h\otimes \nabla \varphi_j(y))+\varphi_j(y) (
\nabla \varphi_i(x)\otimes  \nabla_y h) \label{eq:Hessian-F}  \\ 
+& \varphi_i(y) ( \nabla_y h\otimes \nabla \varphi_j(x))+\varphi_j(x) (
\nabla \varphi_i(y)\otimes  \nabla_x h) \notag \\
+& \varphi_i(y) ( \nabla_y h\otimes \nabla \varphi_j(y))+\varphi_j(y) (
\nabla \varphi_i(y)\otimes  \nabla_y h))dxdy\notag \\
+\int_{\Omega} \int_{ \Omega} (&\varphi_i(x)\varphi_j(x)\nabla^2_{xx}h(x,y)+ \varphi_i(x)\varphi_j(y)\nabla^2_{xy}h(x,y) \notag \\
+&\varphi_i(y)\varphi_j(x)\nabla^2_{yx}h(x,y)+ \varphi_i(y)\varphi_j(y)\nabla^2_{yy}h(x,y))dxdy\,. \notag
\end{align}

\bigskip \bigskip  

{\it Note: On behalf of all authors, the corresponding author states that there is no conflict of interest. }
\bibliographystyle{mybst}
\bibliography{References}

\def\cprime{$'$}
\providecommand{\bysame}{\leavevmode\hbox to3em{\hrulefill}\thinspace}
\providecommand{\MR}{\relax\ifhmode\unskip\space\fi MR }
\providecommand{\MRhref}[2]{%
  \href{http://www.ams.org/mathscinet-getitem?mr=#1}{#2}
}
\providecommand{\href}[2]{#2}
\begin{thebibliography}{10}

\bibitem{A17}
{C.} Audet, \emph{Maximal area of equilateral small polygons}, Amer. Math.
  Monthly \textbf{124} (2017), no.~2, 175--178.

\bibitem{AHMX02}
{C.} Audet, {P.} Hansen, {F.} Messine, and {J.} Xiong, \emph{The largest small
  octagon}, J. Combin. Theory Ser. A \textbf{98} (2002), no.~1, 46--59.

\bibitem{AHS21}
{C.} Audet, {P.} Hansen, and {D.} Svrtan, \emph{Using symbolic calculations to
  determine largest small polygons}, J. Global Optim. \textbf{81} (2021),
  no.~1, 261--268.

\bibitem{B61}
{H.} Bieri, \emph{Ungel\"oste Probleme: Zweiter Nachtrag zu Nr. 12}, Elem.
  Math. \textbf{16} (1961), 105--106.

\bibitem{BB22}
{B.} Bogosel and {D.} Bucur, \emph{On the Polygonal Faber-Krahn Inequality},
  Arxiv arXiv:2203.16409.

\bibitem{BCT22}
{M.} Bonacini, {R.} Cristoferi, and {I.} Topaloglu, \emph{Riesz-type
  inequalities and overdetermined problems for triangles and quadrilaterals},
  J. Geom. Anal. \textbf{32} (2022), no.~2, Paper No. 48, 31.

\bibitem{BBM01}
{J.} Bourgain, {H.} Brezis, and {P.} Mironescu, \emph{Another look at {S}obolev
  spaces}, Optimal control and partial differential equations, IOS, Amsterdam,
  2001, pp.~439--455.

\bibitem{BLL}
{H.J.} Brascamp, {E.H.} Lieb, and {J.M.} Luttinger, \emph{A general
  rearrangement inequality for multiple integrals}, J. Functional Analysis
  \textbf{17} (1974), 227--237.

\bibitem{BuBu}
D.~Bucur and G.~Buttazzo, Variational methods in shape optimization problems,
  Progress in Nonlinear Differential Equations and their Applications, 65,
  Birkh\"auser Boston Inc., Boston, MA, 2005.

\bibitem{BF16}
{D.} Bucur and {I.} Fragal\`a, \emph{A {F}aber-{K}rahn inequality for the
  {C}heeger constant of {$N$}-gons}, J. Geom. Anal. \textbf{26} (2016), no.~1,
  88--117.

\bibitem{BF21}
{D.} Bucur and {I.} Fragal\`a, \emph{Symmetry results for variational energies
  on convex polygons}, ESAIM Control Optim. Calc. Var. \textbf{27} (2021),
  Paper No. 3, 16.

\bibitem{BF}
D.~Bucur and I.~Fragal\`a, \emph{Rigidity for measurable sets}, Adv. Math.
  \textbf{414} (2023).

\bibitem{BFVV}
{D.} Bucur, {I.} Fragal\`a, {B.} Velichkov, and {G.} Verzini, \emph{On the
  honeycomb conjecture for a class of minimal convex partitions}, Trans. Amer.
  Math. Soc. \textbf{370} (2018), no.~10, 7149--7179.

\bibitem{CaffLin}
{L. A.} Caffarelli and {F. H.} Lin, \emph{An optimal partition problem for
  eigenvalues}, J. Sci. Comput. \textbf{31} (2007), no.~1-2, 5--18.

\bibitem{CS08}
{L. A.} Caffarelli and {P.E.} Souganidis, \emph{A rate of convergence for
  monotone finite difference approximations to fully nonlinear, uniformly
  elliptic {PDE}s}, Comm. Pure Appl. Math. \textbf{61} (2008), no.~1, 1--17.

\bibitem{CRS10}
{L.A} Caffarelli, {J.M.} Roquejoffre, and {O.} Savin, \emph{Nonlocal minimal
  surfaces}, Comm. Pure Appl. Math. \textbf{63} (2010), no.~9, 1111--1144.

\bibitem{Ca01}
{T.} Carroll, \emph{Old and new on the bass note, the torsion function and the
  hyperbolic metric}, Irish Math. Soc. Bull. (2001), no.~47, 41--65.
  \MR{1880191}

\bibitem{CDNV}
{A.} Cesaroni, {S.} Dipierro, {M.} Novaga, and {E.} Valdinoci, \emph{Minimizers
  for nonlocal perimeters of {M}inkowski type}, Calc. Var. Partial Differential
  Equations \textbf{57} (2018), no.~2, Paper No. 64, 40.

\bibitem{CN18}
{A.} Cesaroni and {M.} Novaga, \emph{The isoperimetric problem for nonlocal
  perimeters}, Discrete Contin. Dyn. Syst. Ser. S \textbf{11} (2018), no.~3,
  425--440.

\bibitem{degiorgi}
{E.} De~Giorgi, \emph{Sulla propriet\`a isoperimetrica dell'ipersfera, nella
  classe degli insiemi aventi frontiera orientata di misura finita}, Atti
  Accad. Naz. Lincei Mem. Cl. Sci. Fis. Mat. Natur. Sez. Ia (8) \textbf{5}
  (1958), 33--44.

\bibitem{FT}
{G.} Fejes~T\'{o}th, \emph{On the intersection of a convex disc and a polygon},
  Acta Math. Acad. Sci. Hungar. \textbf{29} (1977), no.~1-2, 149--153.

\bibitem{FFMMM}
{A.} Figalli, {N.} Fusco, {F.} Maggi, {V.} Millot, and {M.} Morini,
  \emph{Isoperimetry and stability properties of balls with respect to nonlocal
  energies}, Comm. Math. Phys. \textbf{336} (2015), no.~1, 441--507.

\bibitem{FS07}
{J.} Foster and {T.} Szabo, \emph{Diameter graphs of polygons and the proof of
  a conjecture of {G}raham}, J. Combin. Theory Ser. A \textbf{114} (2007),
  no.~8, 1515--1525.

\bibitem{FV}
{I.} Fragal\`a and {B.} Velichkov, \emph{Serrin-type theorems for triangles},
  Proc. Amer. Math. Soc. \textbf{147} (2019), no.~4, 1615--1626.

\bibitem{FLS08}
{R. L.} Frank, {E.H.} Lieb, and {R.} Seiringer, \emph{Hardy-{L}ieb-{T}hirring
  inequalities for fractional {S}chr\"{o}dinger operators}, J. Amer. Math. Soc.
  \textbf{21} (2008).

\bibitem{FS08}
{R. L.} Frank and {R.} Seiringer, \emph{Non-linear ground state representations
  and sharp {H}ardy inequalities}, J. Funct. Anal. \textbf{255} (2008), no.~12,
  3407--3430.

\bibitem{graham}
{R. L.} Graham, \emph{The largest small hexagon}, J. Combinatorial Theory Ser.
  A \textbf{18} (1975), 165--170.

\bibitem{HLP}
{G.H.} Hardy, {J.E.} Littlewood, and {G.} P\'{o}lya, Inequalities, Cambridge
  Mathematical Library, Cambridge University Press, Cambridge, 1988, Reprint of
  the 1952 edition.

\bibitem{Hmess}
{D.} Henrion and {F.} Messine, \emph{Finding largest small polygons with
  {G}lopti{P}oly}, J. Global Optim. \textbf{56} (2013), no.~3, 1017--1028.

\bibitem{HdG}
{A.} Henrot, Shape optimization and spectral theory, De Gruyter, 2017.

\bibitem{HP}
{A.} Henrot and {M.} Pierre, Shape variation and optimization, EMS Tracts in
  Mathematics, vol.~28, European Mathematical Society (EMS), Z\"{u}rich, 2018.

\bibitem{Laug21}
{R.S.} Laugesen, \emph{Minimizing capacity among linear images of rotationally
  invariant conductors}, Anal. Math. Phys. \textbf{12} (2022), no.~1, Paper No.
  21, 25.

\bibitem{laurain2ndDeriv}
{A.} Laurain, \emph{Distributed and boundary expressions of first and second
  order shape derivatives in nonsmooth domains}, J. Math. Pures Appl. (9)
  \textbf{134} (2020), 328--368. \MR{4053037}

\bibitem{lomb}
{L.} Lombardini, \emph{Fractional perimeters from a fractal perspective}, Adv.
  Nonlinear Stud. \textbf{19} (2019), no.~1, 165--196.

\bibitem{L14}
{M.} Ludwig, \emph{Anisotropic fractional perimeters}, J. Differential Geom.
  \textbf{96} (2014), no.~1, 77--93.

\bibitem{MRT}
{J.M.} Maz\'on, {J.D.} Rossi, and {J.J.} Toledo, Nonlocal perimeter, curvature
  and minimal surfaces for measurable sets, Frontiers in Mathematics,
  Birkh\"{a}user/Springer, Cham, 2019.

\bibitem{morbol}
{F.} Morgan and {R.} Bolton, \emph{Hexagonal economic regions solve the
  location problem}, Amer. Math. Monthly \textbf{109} (2002), no.~2, 165--172.

\bibitem{Moss}
{M.J.} Mossinghoff, \emph{A \$1 problem}, Amer. Math. Monthly \textbf{113}
  (2006), no.~5, 385--402.

\bibitem{PS}
G.~P\'{o}lya and G.~Szeg\"{o}, Isoperimetric {I}nequalities in {M}athematical
  {P}hysics, Annals of Mathematics Studies, no. 27, Princeton University Press,
  Princeton, N. J., 1951.

\bibitem{Pr04}
{M.} Preunkert, \emph{A semigroup version of the isoperimetric inequality},
  Semigroup Forum \textbf{68} (2004), no.~2, 233--245.

\bibitem{R22}
{K.} Reinhardt, \emph{Extremale polygone gegebene Durchmessers}, Jahresber.
  Deutsch. Math.-Verein \textbf{31} (1922), 251--270.

\bibitem{riesz}
{F.} Riesz, \emph{Sur une in\'egalit\'e int\'egrale}, J. London Math. Soc.
  \textbf{5} (1930), no.~3, 162--168.

\bibitem{S58}
{J.J} Sch\"affer, \emph{Ungel\"oste Probleme: Nachtrag zu Nr. 12}, Elem. Math.
  \textbf{13} (1958), 85--86.

\bibitem{SZ}
{A.Y.} Solynin and {V. A.} Zalgaller, \emph{An isoperimetric inequality for
  logarithmic capacity of polygons}, Ann. of Math. (2) \textbf{159} (2004),
  no.~1, 277--303.

\bibitem{BeGi16}
{M.} van~den Berg and {K.} Gittins, \emph{On the heat content of a polygon}, J.
  Geom. Anal. \textbf{26} (2016), no.~3, 2231--2264.

\bibitem{BeSr90}
{M.} van~den Berg and {S.} Srisatkunarajah, \emph{Heat flow and {B}rownian
  motion for a region in {${\bf R}^2$} with a polygonal boundary}, Probab.
  Theory Related Fields \textbf{86} (1990), no.~1, 41--52.

\end{thebibliography}

\end{document}